\documentclass{article}
\usepackage[utf8]{inputenc}
\usepackage{amsmath}
\usepackage{amssymb}
\usepackage{amsthm}
\usepackage{color}
\usepackage{enumitem}
\usepackage[T1]{fontenc}
\usepackage{fullpage}
\usepackage{microtype}
\usepackage{hyperref}
\usepackage{tikz}

\usepackage{subcaption}
\usepackage{soul}

\newtheorem{theorem}{Theorem}[section]
\newtheorem{corollary}[theorem]{Corollary}
\newtheorem{lemma}[theorem]{Lemma}
\newtheorem{definition}[theorem]{Definition}
\newtheorem{proposition}[theorem]{Proposition}
\newtheorem{remark}[theorem]{Remark}
\newtheorem{example}[theorem]{Example}

\title{Real Analytic Multi-parameter Singular Radon Transforms: necessity of the Stein-Street condition}
\author{Lingxiao Zhang}
\date{}

\begin{document}
\maketitle

\begin{abstract}
We study operators of the form
$$
Tf(x)= \psi(x) \int f(\gamma_t(x))K(t)\,dt,
$$
where $\gamma_t(x)$ is a real analytic function of $(t,x)$ mapping from a neighborhood of $(0,0)$ in $\mathbb{R}^N \times \mathbb{R}^n$ into $\mathbb{R}^n$ satisfying $\gamma_0(x)\equiv x$, $\psi(x) \in C_c^\infty(\mathbb{R}^n)$, and $K(t)$ is a ``multi-parameter singular kernel'' with compact support in $\mathbb{R}^N$; for example when $K(t)$ is a product singular kernel. The celebrated work of Christ, Nagel, Stein, and Wainger studied such operators with smooth $\gamma_t(x)$, in the single-parameter case when $K(t)$ is a Calder\'on-Zygmund kernel. Street and Stein generalized their work to the multi-parameter case, and gave sufficient conditions for the $L^p$-boundedness of such operators. This paper shows that when $\gamma_t(x)$ is real analytic, the sufficient conditions of Street and Stein are also necessary for the $L^p$-boundedness of $T$, for all such kernels $K$.
\end{abstract}

\section{Introduction}\label{lyd}
We study operators of the form
\begin{equation}\label{yf}
Tf(x)= \psi(x) \int f(\gamma_t(x))K(t)\,dt,
\end{equation}
where $\gamma_t(x)$ is a real analytic function of $(t,x)$ mapping from a neighborhood of $(0,0)$ in $\mathbb{R}^N \times \mathbb{R}^n$ into $\mathbb{R}^n$ satisfying $\gamma_0(x)\equiv x$, $\psi(x) \in C_c^\infty(\mathbb{R}^n)$, and $K(t)$ is a ``multi-parameter singular kernel'' with compact support in $\mathbb{R}^N$; for example when $K(t)$ is a product singular kernel. The celebrated work of Christ, Nagel, Stein, and Wainger \cite{CNSW} studied such operators with smooth $\gamma_t(x)$, in the single-parameter case when $K(t)$ is a Calder\'on-Zygmund kernel.

The goal is to find necessary and sufficient conditions on the real analytic $\gamma_t(x)$, such that for all $K(t)$ in a given class of ``multi-parameter singular kernels'' supported in a sufficiently small neighborhood of $0$ in $\mathbb{R}^N$, and for all smooth $\psi(x)$ supported in a sufficiently small neighborhood of $0$ in $\mathbb{R}^n$, the corresponding operators $T$ in (\ref{yf}) are bounded on $L^p(\mathbb{R}^n)$ for all $1<p<\infty$. (In fact, this is equivalent to the boundedness on $L^p(\mathbb{R}^n)$ for some $1<p<\infty$; see Theorem \ref{lieri}).

Street and Stein \cite{L2, LP, ANALYTIC} gave sufficient conditions on a $C^\infty$ $\gamma_t(x)$, for all such operators to be bounded on $L^p(\mathbb{R}^n)$ for all $1<p<\infty$, (and in particular, in the single-parameter case, such operators with real analytic $\gamma_t(x)$ are always bounded). This paper proves these sufficient conditions are also necessary when $\gamma_t(x)$ is real analytic.

For simplicity, we introduce in this section the necessary and sufficient condition in the special case that $K$ is a product kernel with compact support in $\mathbb{R}^2$. We suppose $\gamma_t(x)$ is real analytic. We say a distribution $K(s,t)$ is a product kernel on $\mathbb{R}^2$, provided $K(s,t)$ agrees with a $C^\infty$ function on $\mathbb{R}^2$ away from $s=0$ and $t=0$, satisfying the following differential inequalities: for each $\alpha=(\alpha_1, \alpha_2) \in \mathbb{N}^2$, and for $s \neq 0, t \neq 0$,
$$
|\partial_{s}^{\alpha_1} \partial_{t}^{\alpha_2} K(s,t)| \leq C_\alpha |s|^{-(1+\alpha_1)} |t|^{-(1+\alpha_2)}, \text{~for some constant~} C_\alpha>0,
$$
and $K(s,t)$ satisfies certain ``cancellation conditions'' (see Definition \ref{productkernel}). For example, $K(s,t)= \frac{1}{st}$.

Note that for every small $(s,t)$, $\gamma_{(s,t)}(\cdot )$ is a local diffeomorphism near $0$ in $\mathbb{R}^n$. We let $\gamma_{(s,t)}^{-1}(\cdot)$ denote the inverse maps. We define
$$
W((s,t),x) = \frac{d}{d\epsilon}\Big|_{\epsilon =1} \gamma_{(\epsilon s, \epsilon t)} \circ \gamma_{(s,t)}^{-1}(x)\in T_x \mathbb{R}^n.
$$
Write the Taylor series of $W$ in $(s,t)$:
\begin{equation}
W((s,t),x) \sim \sum_{|\alpha|>0}(s,t)^\alpha X_\alpha,
\end{equation}
where $\alpha=(\alpha_1, \alpha_2) \in \mathbb{N}^2$ denotes a multi-index, $(s,t)^\alpha = s^{\alpha_1} \cdot t^{\alpha_2}$, and the $X_\alpha$ are real analytic vector fields on a neighborhood of $0$ in $\mathbb{R}^n$. We say $\alpha$ is a pure power if $\alpha$ is nonzero only in one component; otherwise we say $\alpha$ is a nonpure power. Let
\begin{align*}
&\mathcal{P}=\{X_\alpha: \alpha \text{~is a pure power}\},\\
&\mathcal{N}=\{X_\alpha: \alpha \text{~is a nonpure power}\},
\end{align*}
and let $\mathcal{L}(\mathcal{P})$ be the smallest set containing $\mathcal{P}$ and closed under Lie brackets.

A necessary condition on the real analytic $\gamma_{(s,t)}(x)$, for all operators $T$ associated with smooth functions $\psi$ and product kernels $K$, both with sufficiently small supports, to be bounded is that every element of $\mathcal{N}$ is spanned by $\mathcal{L}(\mathcal{P})$ with bounded coefficients on some neighborhood of $0$. However, this condition is not sufficient. For every $\delta =(\delta_1, \delta_2) \in (0,1]^2$, denote $\delta^\alpha = \delta_1^{\alpha_1}\cdot \delta_2^{\alpha_2}$, define
\begin{align*}
&\mathcal{P}_\delta = \{\delta^\alpha X_\alpha: \alpha \text{ is a pure power}\},\\
&\mathcal{N}_\delta = \{\delta^\alpha X_\alpha: \alpha \text{ is a nonpure power}\},
\end{align*}
and let $\mathcal{L}(\mathcal{P}_\delta)$ be the smallest set containing $\mathcal{P}_\delta$ and closed under Lie brackets. The necessary and sufficient condition is that every element of $\mathcal{N}_\delta$ is spanned by $\mathcal{L}(\mathcal{P}_\delta)$ with bounded coefficients in a neighborhood of $0$ ``uniformly in $\delta$''; see Theorem \ref{fy}. For instance, see Theorem \ref{illuminate} for this uniformity condition in the $\mathbb{R}^1$ polynomial case.

For an example of $\gamma_{(s,t)}(x)$, not satisfying the above necessary condition, consider $\gamma_{(s,t)}(x)= x+st$ defined on $\mathbb{R}^2 \times \mathbb{R}$. We have $W((s,t),x) = 2st \partial_x$, and therefore $2\partial_x \in \mathcal{N}$, but $\mathcal{L}(\mathcal{P})$ only consists of zero vector fields. In fact, there exist product kernels with arbitrarily small supports such that the corresponding operators are unbounded on $L^p(\mathbb{R})$ for all $1< p<\infty$ (this idea first dates back in \cite{17}, see Section 17.5 of \cite{L2} for a proof).

There have been many papers which give sufficient conditions for the $L^p$ boundedness of operators of the form (\ref{yf}). The first paper (Fabes \cite{FABES}) on the operator (\ref{yf}) was about the Hilbert transform on the parabola, where $\gamma_t(x_1, x_2) = (x_1 -t, x_2 -t^2)$; this is a single-parameter singular Radon transform. After this first paper, many authors studied the $L^p$ boundedness of single-parameter singular Radon transforms. This culminated in the important paper by Christ, Nagel, Stein, and Wainger \cite{CNSW}; see that paper for a detailed history of the subject. More recently, Street and Stein \cite{L2, LP, ANALYTIC} generalized the operators to the multi-parameter case and studied the $L^p$ boundedness; even in the single-parameter case, their results are more general than those in \cite{CNSW}. Street \cite{BRIAN17} later studied the boundedness in Sobolev spaces of the ``multi-parameter'' singular Radon transforms. Some other directly related works include \cite{GREENBLATT1} and \cite{SW03} on the single-parameter case, \cite{GREENBLATT5}, \cite{GREENBLATT3}, and \cite{GREENBLATT4} on the $L^p$ Sobolev regularity, \cite{DOUBLE} and \cite{12} which gave uniform estimates, and \cite{CHKY08} which studied such operators along flat surfaces.

There are few papers on the necessary and sufficient conditions. Carbery, Wainger, and Wright \cite{CWW00} gave a necessary and sufficient condition for (local) double Hilbert transforms along polynomial surfaces in $\mathbb{R}^3$ to be bounded on $L^p (1<p<\infty)$. See also Patel \cite{d=3} and Patel \cite{d=1} for extensions of this result to global Hilbert transforms and to polynomials in $\mathbb{R}^1$, respectively. Carbery, Wainger, and Wright (Remark 1.4 in \cite{CWW06}) gave a necessary and sufficient condition for double Hilbert transforms along real analytic surfaces in $\mathbb{R}^3$ to be bounded on $L^p (1<p<\infty)$. See also Pramanik and Yang \cite{8} for an extension to real analytic surfaces in $\mathbb{R}^{d+2}$. Carbery, Wainger, and Wright \cite{9} gave necessary and sufficient conditions for triple Hilbert transforms along polynomial surfaces in $\mathbb{R}^4$ to be bounded on $L^2$. See also Cho, Hong, Kim, and Yang \cite{TRIPLE} for an extension to the $L^p$ boundedness. Kim \cite{KIM} gave necessary and sufficient conditions for multiple Hilbert transforms along polynomial surfaces to be bounded on $L^p (1<p < \infty)$.
However, these papers are different from ours, in that they fix the kernel $K$, and that they only consider translation-invariant operators.

Note our conditions are necessary only when operators associated to all kernels are bounded, and are not necessary when we consider a particular kernel. Of course when $K=0$, no conditions are needed. More interestingly, Vitturi and Wright (Theorem 3.1 of \cite{12}) gave a necessary and sufficient condition for the $L^2$ boundedness of such operators on the Heisenberg group $\mathbb{H}^1$ with a multi-parameter singular kernel in terms of the $L^2$ boundedness of such operators on $\mathbb{R}^2$ with the same kernel and with certain truncations. We refer to this paper for an extensive bibliography of other works on translation-invariant singular Radon transforms. Besides, in the case when $K(s,t) = \frac{1}{st}$, and when the operator is along a polynomial surface $(s,t, P(s,t))$ in $\mathbb{H}^1$, Theorem 1.2 of \cite{12} gave a necessary and sufficient condition for $L^2(\mathbb{H}^1)$ boundedness of such operators uniformly in the coefficients of $P$, in terms of the powers of $(s,t)$ appearing in the polynomial $P$.

Also note that by a change of variable, for any $k\in \mathbb{N}$ odd, any real analytic function $P(s,t)$, and any small supported product kernel $K(s,t)$, the $L^p$ boundedness of the operator in $\mathbb{H}^1$:
$$
f \mapsto \int_{\mathbb{R}^2} f\big( (x,y,z) \cdot (s,t,P(s^{1/k},t))^{-1} \big) K(s,t)\,ds\,dt,
$$
is equivalent to the $L^p(\mathbb{H}^1)$ boundedness of
$$
f \mapsto \int_{\mathbb{R}^2} f\big( (x,y,z) \cdot (s^k,t,P(s,t))^{-1} \big) \widetilde K(s,t)\,ds\,dt,
$$
where $\widetilde K$ is another small supported product kernel. Therefore Theorem 3.1 of \cite{12} also holds for certain non-real analytic cases. To be more precise, Theorem 3.1 of \cite{12} immediately implies that for odd $k$ and real analytic $P(s,t)$, the $L^2$ boundedness of the operator along $(s,t, P(s^{1/k}, t))$ in $\mathbb{H}^1$ with a small supported product kernel is equivalent to the $L^2$ boundedness of the operator along $(t, P_k(s^{1/k}, t))$ in $\mathbb{R}^2$ with the same kernel and with certain truncations, where the real analytic function $P_k(s,t)$ is obtained from $P(s,t)$ by removing nonpure terms $s^{\alpha_1}t^{\alpha_2}$ with $\alpha_1\geq k$. Likewise, by a change of variables, current work also extends to certain non-real analytic cases: our conditions on a ``fractional surface'' $(P_1(s^{1/k}, t), P_2(s^{1/k}, t),P_3(s^{1/k}, t))$ in $\mathbb{H}^1$ with odd $k$ and real analytic $P_1, P_2$, and $P_3$, are still necessary and sufficient for the $L^p$ boundedness of operators along the surface with all small supported product kernels.

The main tool we use is a quantitative and uniform version of the classical theorem of Frobenius due to Street, which is Theorem 4.1 in \cite{CARNOT} (see Theorem \ref{16}). This is related to Carnot-Carath\'eodory geometry. The earliest such result is by Nagel, Stein, and Wainger \cite{NSW85}, which is a special case of Theorem 4.1 in \cite{CARNOT}, in the single-parameter case. Tao and Wright (Section 4 of \cite{TW03}) later developed this into a uniform version in the multi-parameter case, which is also a special case of Theorem 4.1 in \cite{CARNOT}. The results of \cite{CARNOT} were later strengthened and generalized by Stovall and Street \cite{BETSY, COORDINATES, STREET20}. Carnot-Carath\'eodory geometry has long played an important role in these operators. This was first laid out in \cite{CNSW}. See also e.g. \cite{L2}, \cite{LP}, \cite{ANALYTIC}, \cite{BRIAN17}, \cite{BETSY11}, \cite{STOVALL}, \cite{GREENBLATT}, \cite{TW03}, \cite{GREENBLATT6}, \cite{ENDPOINT}, \cite{SPARSE}, \cite{GREENBLATT2}, \cite{CHRISTSTREET}.

In Section \ref{raj}, we will state the main results of this paper, which are the necessary and sufficient conditions. The main goal of this paper is to show these conditions are necessary, the sufficiency having been proved by Street and Stein in Theorem 5.1 of \cite{LP}. In Section \ref{real}, we will explain the main results and the basic ideas of the proof in the special case of the real line $\mathbb{R}^1$. In Section \ref{sntz}, we will give an outline of the proof in the special case of the Heisenberg group $\mathbb{H}^1$, to motivate our proof in the general real analytic case. In Section \ref{boruozhang}, we will cite some properties of real analytic functions. In Section \ref{shaolinchangquan}, we will state the main result in a more complete form, and we will divide it into proving three main propositions: Propositions \ref{1}, \ref{2}, and \ref{3}. In Section \ref{nianhuazhi}, we will prove or cite several lemmas or theorems, in preparation for our proof. The main parts of this paper are Sections \ref{A}, \ref{xj}, and \ref{as}, which are the proofs for Propositions \ref{1}, \ref{2}, \ref{3}, respectively. We put the proofs of some lemmas used in the above three sections in Appendix \ref{jinzhangzhi}.

\textbf{Acknowledgements.} The author would like to thank her advisor Brian Street for introducing this problem, for his guidance and illuminating discussions. The author would also like to thank the anonymous referees for many helpful comments and suggestions.

\section{The main results}\label{raj}
We have informally described in Section \ref{lyd} a necessary and sufficient condition in the product kernel case. To state precisely the necessary and sufficient conditions in the general case, we start with several definitions.

We first make precise the notion of a ``multi-parameter singular kernel'' with compact support in $\mathbb{R}^N$. Let $\nu$ be the number of parameters. Besides $N, \nu$, we need a set of vectors $e=\{e_1, \ldots, e_N\} \subseteq [0,\infty)^\nu$ to determine the multi-parameter dilations for a class of kernels. We also write $e=\{e_i^\mu: 1 \leq i \leq N, 1\leq \mu \leq \nu\}$. We require that
\begin{equation}\label{e}
\text{for each~} 1\leq i \leq N, e_i^\mu \neq 0 \text{~for some~} 1\leq \mu \leq \nu,\quad \text{for each~} 1\leq \mu \leq \nu, e_i^\mu \neq 0 \text{~for some~} 1\leq i \leq N.
\end{equation}
Using $e$, we define multi-parameter dilations on $\mathbb{R}^N$ as follows: for $t=(t_1, \ldots, t_N) \in \mathbb{R}^N$ and $\delta \in [0, \infty)^\nu$,
$$
\delta t = (\delta^{e_1}t_1, \ldots, \delta^{e_N}t_N), \text{~where~} \delta^{e_j}= \delta_1^{e_j^1} \cdots \delta_\nu^{e_j^\nu}.
$$
For $\varsigma\in C_c^\infty(\mathbb{R}^N)$, and $\delta \in (0,\infty)^\nu$, we denote 
$$
\varsigma^{(\delta)}(t) = \delta^{e_1+\cdots +e_N}\varsigma(\delta t).
$$
Note that $\int \varsigma^{(\delta)}(t)\,dt = \int \varsigma(t)\,dt$.

For $k=(k_1, \ldots, k_\nu) \in \mathbb{R}^\nu$, denote $2^k:= (2^{k_1}, \ldots, 2^{k_\nu})\in (0,\infty)^\nu$. Let $t^\mu$ be the vector consisting of the coordinates $t_i$ of $t$ with $e_i^\mu \neq 0$. Note that $t^\mu$ and $t^\lambda (\mu \neq \lambda)$ might contain some of the same coordinates. For $a>0$, let $B^N(a)$ denote the $N$-dimensional ball centered at $0$ with radius $a$.

\begin{definition}
For such $N, e, a, \nu$, we say a distribution $K$ belongs to the multi-parameter singular kernel class $\mathcal{K}(N,e,a, \nu)$, if there exists a bounded sequence $\{\varsigma_k\}_{k\in \mathbb{N}^\nu}$ in $C_c^\infty(B^N(a))$ satisfying the cancellation conditions
$$
\int \varsigma_k(t)\,dt^\mu \equiv 0, \text{ for all } k\in \mathbb{N}^\nu, \mu \in \{1, \ldots, \nu\} \text{~with~} k_\mu \neq 0,
$$
where $\int \varsigma_k(t) dt^\mu$ is a function in coordinates $t_i$ not in $t^\mu$, such that we can write
$$
K(t) = \sum_{k\in \mathbb{N}^\nu} \varsigma_k^{(2^k)}(t),
$$
where the above infinite sum converges in the sense of distributions. 
\end{definition}

\begin{remark}
\begin{enumerate}
\item
By Lemma 16.1 in \cite{L2}, $\sum_{k\in \mathbb{N}^\nu} \varsigma_k^{(2^k)}$ indeed converges in the sense of distributions with compact support, for any such sequence. Thus the operator 
$$
Tf(x) =\psi(x) \int f(\gamma_t(x)) K(t) \,dt = \psi(x) \sum_{k\in \mathbb{N}^\nu} \int f(\gamma_t(x)) \varsigma_k^{(2^k)}(t) \,dt
$$ 
is well-defined for $f\in C_c^\infty(\mathbb{R}^n)$. Note the kernel $K(t)$ has a compact support contained in $\overline{B^N(a)}$.

\item
If $e$ is such that each vector $e_i$ only has one nonzero component, then $\mathcal{K}(N,e,a,\nu)$ is a class of product kernels. This definition of product kernels is equivalent to the usual definition of product kernels using differential inequalities and certain cancellation conditions (see Proposition 16.5 of \cite{L2} and Definition \ref{productkernel}).

\item
When $\nu =1$, $\mathcal{K}(N,e,a,\nu)$ is a class of standard Calder\'on-Zygmund kernels (see Proposition 16.5 of \cite{L2} and Definition \ref{productkernel}).

\end{enumerate}
\end{remark}

As in Section 1, we will use the given $C^\infty$ $\gamma_t(x)$ to produce a $C^\infty$ vector field depending on $t$: $W(t,x)$, and formulate the necessary and sufficient conditions in terms of the vector fields in the Taylor expansion of $W(t,x)$. 

For $t$ sufficiently small and fixed, $\gamma_t(\cdot)$ is a diffeomorphism from a neighborhood of $x=0$ onto its image. For a scalar $\epsilon >0$, we set $\epsilon t = (\epsilon t_1, \ldots, \epsilon t_N)$. Note this is different from the multi-parameter dilation $\delta t$, where $\delta$ is not a scalar.

\begin{definition}\label{def2.3}
Define 
$$
W(t,x)= \frac{d}{d\epsilon}\Big|_{\epsilon =1} \gamma_{\epsilon t} \circ \gamma_t^{-1}(x) \in T_x\mathbb{R}^n.
$$
Note $W(t)$ is a $C^\infty$ vector field, depending smoothly on $t$, satisfying $W(0) = 0$.
\end{definition}

\begin{lemma}\label{12.1} (Proposition 12.1 in \cite{L2})
The map $\gamma \mapsto W$ is a bijection from smooth functions $\gamma_t$, with $\gamma_0(x) \equiv x$ (thought of as germs in the $t$ variable), to smooth vector fields depending on $t$, $W(t)$, with $W(0) = 0$ (also thought of as germs in the $t$ variable). The inverse map $W \mapsto \gamma$ is obtained by setting $\gamma_t(x) = w(1,t,x)$ for the solution $w(\epsilon,t,x)$ to the ODE
$$
\frac{d}{d\epsilon} w(\epsilon,t,x) = \frac{1}{\epsilon} W(\epsilon t, w(\epsilon, t, x)), \quad w(0,t,x) =x.
$$
\end{lemma}

Write the Taylor series of $W$ in $t$ variable,
$$
W(t) \sim \sum_{|\alpha|>0} t^\alpha X_\alpha,
$$
where $\alpha \in \mathbb{N}^N$ denotes a multi-index, and the $X_\alpha$ are $C^\infty$ vector fields on a common neighborhood of $0$ in $\mathbb{R}^n$. Define the \textbf{degree}\footnote{Different $\alpha$'s may have the same degree, but each $d\in [0,\infty)^\nu \backslash\{0\}$ is the degree of only finitely many $\alpha$'s.} of each multi-index $\alpha \in \mathbb{N}^N\backslash \{0\}$ as
$$
\deg(\alpha)= \alpha_1 e_1 + \cdots + \alpha_N e_N\in [0,\infty)^\nu\backslash \{0\}.
$$
We say $\alpha$ is a \textbf{pure power} if $\deg(\alpha)$ is nonzero only in one component; otherwise we say $\alpha$ is a \textbf{nonpure power}.

\begin{definition}
Define the pure power set and the nonpure power set of $W(t,x)$ as
\begin{align*}
&\mathcal{P}=\{(X_\alpha, \deg(\alpha)): \alpha \text{~is a pure power}\},\\
&\mathcal{N}=\{(X_\alpha, \deg(\alpha)): \alpha \text{~is a nonpure power}\}.
\end{align*}
\end{definition}

\begin{definition}
For $\mathcal{S}$ a set of $C^\infty$ vector fields on a common neighborhood of $0$ in $\mathbb{R}^n$ each paired with a $\nu$-parameter degree $d\in [0,\infty)^\nu \backslash \{0\}$, we define $\mathcal{L}(S)$ to be the smallest set of vector fields on this common neighborhood paired with degrees such that
\begin{itemize}
\item 
$\mathcal{S} \subseteq \mathcal{L}(\mathcal{S})$,
\item
if $(X_1, d_1), (X_2, d_2) \in \mathcal{L}(\mathcal{S})$ then $([X_1,X_2], d_1+d_2) \in \mathcal{L}(\mathcal{S})$.
\end{itemize}
\end{definition}

Vector fields with degrees give rise to Carnot-Carath\'eodory geometry. Let $(X_1,d_1), \ldots, (X_q,d_q)$ be finitely many $C^\infty$ vector fields on some common neighborhood $\Omega$ of $0$ in $\mathbb{R}^n$, each paired with a degree in $[0, \infty)^\nu \backslash \{0\}$. Denote $(\mathbf{X}, \mathbf{d}):=\{(X_1,d_1), \ldots, (X_q,d_q)\}$.

\begin{definition}
The multi-parameter Carnot-Carath\'eodory ball in $\Omega$ with center $x\in \Omega$ and radius $\delta \in [0,\infty)^\nu$ with respect to $(\mathbf{X},\mathbf{d})$ is the set
\begin{align*}
B_{(\mathbf{X},\mathbf{d})}(x,\delta) =&\Big\{ y \in \Omega: \text{~there exists an absolutely continuous curve~}\\
& \quad \gamma: [0,1]\to \Omega \text{~with~} \gamma(0)=x, \gamma(1)=y,\text{~and satisfying}\\
&\quad \gamma'(t) =\sum_{l=1}^q a_l(t) \delta^{d_l}X_l(\gamma(t)),\text{a.e.}, ~ \Big\|\sum_{l=1}^q \big|a_l(t)\big|^2\Big\|_{L^\infty[0,1]}<1\Big\},
\end{align*}
equipped with the subspace topology.
\end{definition}

Note that $\delta$ roughly represents the size of the ``ball'' $B_{(\mathbf{X},\mathbf{d})}(x,\delta)$. But if $x\in \Omega$ is too close to the boundary of $\Omega$, $B_{(\mathbf{X},\mathbf{d})}(x,\delta)$ would be much smaller than a $\delta$-ball, since the ``frame'' $(\mathbf{X},\mathbf{d})$ is only defined on $\Omega$. The following $\mathcal{C}(x,\delta, \Omega)$ condition is to guarantee that the ``center of ball'' $x$ is not too close to the boundary of $\Omega$.

\begin{definition}
For $x\in \Omega$, and $\delta \in [0,\infty)^\nu$, we say $(\mathbf{X},\mathbf{d})$ satisfies $\mathcal{C}(x,\delta, \Omega)$, if for every $a=(a_1, \ldots, a_q):[0,1]\to \mathbb{R}^q$ with $\Big\|\sum_{l=1}^q\big|a_l(t)\big|^2\Big\|_{L^\infty[0,1]}<1$, there exists an absolutely continuous solution $\gamma: [0,1]\to \Omega$ to the ODE:
$$
\gamma'(t) = \sum_{l=1}^q a_l(t)\delta^{d_l}X_l(\gamma(t))\text{~a.e.}, \quad \gamma(0)=x.
$$
Note, by Gronwall's inequality, when this solution exists, it is unique.
\end{definition}

To state the necessary and sufficient conditions, we now start to define the notion of control, which is the uniformity condition mentioned in Section \ref{lyd}. 
\begin{itemize}
\item 
Given $(\mathbf{X}, \mathbf{d})$, denote $\delta \mathbf{X} := \{\delta^{d_1} X_1, \ldots, \delta^{d_q}X_q\}$. 
\item
An ``ordered multi-index'' $\alpha$ is a list of numbers from $\{1, \ldots, q\}$, and $|\alpha|$ denotes the length of this list.
\item
$(\delta \mathbf{X})^\alpha$ means we juxtapose vector fields from $\{\delta^{d_1} X_1, \ldots, \delta^{d_q}X_q\}$ indexed with list $\alpha$. For instance, if $\alpha = (1,2,2)$, then $|\alpha| =3$, and $(\delta \mathbf{X})^\alpha = (\delta^{d_1}X_1)(\delta^{d_2}X_2)(\delta^{d_2}X_2)$. 
\item
For every $\xi>0$, denote $\Vec{\xi} := (\xi, \ldots, \xi) \in (0, \infty)^\nu$. 
\end{itemize}
Let $(X_0, d_0)$ be another $C^\infty$ vector field on $\Omega$ with a $\nu$-parameter degree.

\begin{definition}\label{controldefn}
We say $(\boldsymbol X,\boldsymbol d)$ controls $(X_0,d_0)$ if there exists some neighborhood $U\subseteq \Omega$ of $0$ and $\xi_1>0$, such that $(\boldsymbol X,\boldsymbol d)$ satisfies $\mathcal{C}(x,\Vec{\xi_1},\Omega)$ for all $x\in U$, and that for all $x\in U, \delta\in [0,1]^\nu$,
$$
\delta^{d_0} X_0 =\sum_{l=1}^q c_l^{x,\delta}\delta^{d_l}X_l, \text{~on~} B_{(\mathbf{X},\mathbf{d})}(x,\xi_1 \delta),
$$
and\footnote{We write $\|f\|_{C(V)} := \sup_{u\in V} |f(u)|$, and if we say the norm is finite, we mean (in addition) that $f$ is continuous on $V$, where $V$ is given the subspace topology.}
$$
\sup_{\substack{\delta\in [0,1]^\nu \\ x\in U}}\sum_{|\alpha| \leq m} \Big\|(\delta \mathbf{X})^\alpha c_l^{x,\delta}\Big\|_{C(B_{(\mathbf{X},\mathbf{d})}(x,\xi_1\delta))} <\infty, \text{~for every~} m \in \mathbb{N}.
$$

When we wish to make the choice of $U, \xi_1$ explicit, we say, $(\boldsymbol X,\boldsymbol d)$ controls $(X_0,d_0)$ on $U$ with parameter $\xi_1$.
\end{definition}

As in Lemma \ref{CONTAIN}, under certain assumptions, the multi-parameter Carnot-Carath\'eodory ball generated by $(\boldsymbol Z,\boldsymbol D)$ is contained in the multi-parameter Carnot-Carath\'eodory ball generated by $(\boldsymbol X,\boldsymbol d)$ of a comparable radius, if $(\boldsymbol X,\boldsymbol d)$ controls every element in $(\boldsymbol Z,\boldsymbol D)$.

\begin{definition}\label{defn}
Given two (possibly infinite) sets $\mathcal{S}, \mathcal{T}$ of vector fields paired with $\nu$-parameter degrees, we say $\mathcal{S}$ controls $\mathcal{T}$ if there exists a finite subset $\mathcal{S}_0$ of $\mathcal{S}$ such that $\mathcal{S}_0$ controls every element of $\mathcal{S} \cup \mathcal{T}$ \textup{(}the particular neighborhood $U$ and $\xi_1$ used for control may depend on each element of $\mathcal{S} \cup \mathcal{T}$\textup{)}.
\end{definition}

\begin{remark}
The definition of control in Definition \ref{defn} differs slightly from the definition in Street and Stein \cite{ANALYTIC}. The transitivity property of control stated in Lemma 8.12 of \cite{ANALYTIC} does not seem to be fully proved there. However, the rest of that paper holds with the notion of control given in Definition \ref{defn} in place of the notion used in that paper.
\end{remark}

The main result is the following theorem (see Theorem \ref{lieri} for a more complete version).

\begin{theorem}\label{fy}
Let $\gamma_t(x)$ be a real analytic function of $(t,x)$ defined on a neighborhood of $(0,0)$ in $\mathbb{R}^N\times \mathbb{R}^n$ with $\gamma_0(x)\equiv x$. The following are equivalent.

\begin{enumerate}[label=(\alph*)]
\item\label{BOUND}
The operator
$$
Tf(x)=\psi(x) \int f(\gamma_t(x))K(t)\,dt
$$
is bounded on $L^p(\mathbb{R}^n)$, for any $1<p<\infty$, any $\psi \in C_c^\infty(\mathbb{R}^n)$ with sufficiently small support, and any $K \in \mathcal{K}(N,e,a,\nu)$ with sufficiently small $a>0$.

\item\label{CONTROL}
$\mathcal{L}(\mathcal{P})$ controls $\mathcal{N}$.
\end{enumerate}
\end{theorem}

The main goal is to show \ref{BOUND} $\Rightarrow$ \ref{CONTROL}. \ref{CONTROL} $\Rightarrow$ \ref{BOUND} has been proved by Street and Stein in Theorem 5.1 of \cite{LP} (see Section \ref{shaolinchangquan} in this paper).

\begin{remark}
When $\nu=1$, \ref{CONTROL} automatically holds, thus the single-parameter operator corresponding to real analytic $\gamma_t(x)$ is always bounded on $L^p(\mathbb{R}^n)$, for any $1<p<\infty$, any $\psi \in C_c^\infty(\mathbb{R}^n)$ with sufficiently small support, and any Calder\'on-Zygmund kernel $K$ with sufficiently small support, which is a result of Street and Stein \cite{ANALYTIC}.
\end{remark}

\begin{example}
Let $e= \{(1,0), (0,1)\}$. Then $W(t_1, t_2,x) = t_1^2 + t_2^2 + t_1 t_2$ is an example such that $\mathcal{L}(\mathcal{P})$ controls $\mathcal{N}$, where $\mathcal{L}(\mathcal{P}) = \{ (\partial_x, (2,0)), (\partial_x, (0,2))\}$, and $\mathcal{N}= \{ (\partial_x, (1,1))\}$. But $W(t_1, t_2, x) = t_1^3 +t_2^3 + t_1 t_2$ is such that $\mathcal{L}(\mathcal{P})$ does not control $\mathcal{N}$, where $\mathcal{L}(\mathcal{P}) = \{ (\partial_x, (3,0)), (\partial_x, (0,3))\}$, and $\mathcal{N}= \{ (\partial_x, (1,1))\}$. Also, $W(t_1, t_2,x_1, x_2) = t_1 \partial_{x_1} + t_2 x_1 \partial_{x_2} + t_1 t_2 \partial_{x_2}$ is another example such that $\mathcal{L}(\mathcal{P})$ controls $\mathcal{N}$, although $\mathcal{P}$ does not control $\mathcal{N}$; here $\mathcal{P} = \{ (\partial_{x_1}, (1,0)), (x_1\partial_{x_2}, (0,1))\}$, $\mathcal{N}= \{ (\partial_{x_2}, (1,1))\}$, and $\mathcal{L}(\mathcal{P}) = \{ (\partial_{x_1}, (1,0)), (x_1\partial_{x_2}, (0,1)), (\partial_{x_2}, (1,1))\}$.
\end{example}

We can also give a similar necessary and sufficient condition using a different collection of vector fields. By Theorem 8.5 in \cite{CNSW}, for any $C^\infty$ function $\gamma_t(x)$ defined on a neighborhood of $(0,0)$ in $\mathbb{R}^N \times \mathbb{R}^n$ satisfying $\gamma_0(x)\equiv x$, we can write 
\begin{equation}\label{Z}
\gamma_t(x) \sim e^{\sum_{|\alpha|>0}t^\alpha \widehat X_\alpha}x,
\end{equation}
where the $\widehat X_\alpha$ are $C^\infty$ vector fields, uniquely determined by $\gamma_t(x)$. (\ref{Z}) means
$$
\gamma_t(x) = e^{\sum_{0<|\alpha|<M}t^\alpha \widehat X_\alpha}x + O(|t|^M), \text{~for every~} M.
$$
In fact, the collection of vector fields $\{\widehat X_\alpha\}$ is equivalent to $\{X_\alpha\}$ for our purposes, which is made precise in Theorem \ref{ljj}.

\begin{definition}
The pure power set and the nonpure power set of $\gamma_t(x)$ are defined as
\begin{align*}
&\widehat{\mathcal{P}}=\{(\widehat X_\alpha, \deg(\alpha)): \alpha \text{~is a pure power}\},\\
&\widehat{\mathcal{N}}=\{(\widehat X_\alpha, \deg(\alpha)): \alpha \text{~is a nonpure power}\}.
\end{align*} 
\end{definition}

Theorem \ref{fy} is equivalent to the following theorem (see Corollary \ref{EQUIV}).

\begin{theorem}\label{newly}
Let $\gamma_t(x)$ be a real analytic function of $(t,x)$ defined on a neighborhood of $(0,0)$ in $\mathbb{R}^N\times \mathbb{R}^n$ with $\gamma_0(x)\equiv x$. The following are equivalent.

\begin{enumerate}[label=(\alph*)]
\item
The operator
$$
Tf(x)=\psi(x) \int f(\gamma_t(x))K(t)\,dt
$$
is bounded on $L^p(\mathbb{R}^n)$, for any $1<p<\infty$, any $\psi \in C_c^\infty(\mathbb{R}^n)$ with sufficiently small support, and any $K \in \mathcal{K}(N,e,a,\nu)$ with sufficiently small $a>0$.

\item
$\mathcal{L}(\widehat{\mathcal{P}})$ controls $\widehat{\mathcal{N}}$.
\end{enumerate}
\end{theorem}

\section{Real line and product kernel case}\label{real}
Consider in this section operators on the real line with product kernels. In this very simple special case, it is easy to understand our results and the main ideas of the proof.

We consider product kernels with compact support in $\mathbb{R}^2$. That is we assume in this section $K \in \mathcal{K}(2,e,a,2)$ with $a>0$, $e=\{e_1, e_2\}, e_1 = (1,0), e_2=(0,1)$. See the following for an equivalent definition.

A Calder\'on-Zygmund kernel on $\mathbb{R}$ is a distribution $K(t)$ which coincides with a smooth function away from the origin such that $\big|\frac{d^m}{dt^m}K(t)\big| \lesssim_m |t|^{-1-m}$ for all $m\in \mathbb{N}$ and such that the quantities $\int K(t)\phi(Rt)\,dt$ are bounded uniformly over all $R>0$ and all smooth $\phi$ supported in the unit ball with $\|\phi\|_{C^1} \leq 1$ (such a $\phi$ is called a normalized bump function on $\mathbb{R}$).

\begin{definition}\label{productkernel}
(Definition 2.1.1 of \cite{NRS01}) A product kernel on $\mathbb{R}^2$ is a distribution $K(s,t)$ which coincides with a $C^\infty$ function away from the coordinate axes $s=0, t=0$ and which satisfies 
\begin{enumerate}
\item 
(Differential inequalities) for each multi-index $\alpha = (\alpha_1, \alpha_2)\in \mathbb{N}^2$, there is a constant $C_\alpha$ such that $|\partial_s^{\alpha_1} \partial_t^{\alpha_2} K(s,t)| \leq C_\alpha |s|^{-1-\alpha_1} |t|^{-1-\alpha_2}$ away from the coordinate axes, and
\item
(Cancellation conditions) for any normalized bump function $\phi$ on $\mathbb{R}$, and any $R>0$, the distributions $K_{\phi, R}(s) = \int K(s,t) \phi(Rt)\,dt$ and $K^{\phi, R}(t) = \int K(s,t)\phi(Rs)\,ds$ are Calder\'on-Zygmund kernels on $\mathbb{R}$ as described above, uniformly in $\phi$,$R$.
\end{enumerate}
\end{definition}

For example, $\frac{\phi(s)\phi(t)}{s t} \in \mathcal{K}(2, e, a, 2)$ if $\phi\in C_c^\infty\big(B^1(\frac{a}{\sqrt{2}})\big)$.

We consider functions $ \gamma_{(s,t)}(x): \mathbb{R}^2\times \mathbb{R}^1 \to \mathbb{R}^1$ of the form
$$
\gamma_{(s,t)}(x) =e^{-p(s,t) \partial_x}x = x-p(s,t),
$$
where $p(s,t)$ is a polynomial without constant term.

\begin{example}\label{kitty}
For every $a>0$, there exists $K\in \mathcal{K}(2,e,a,2)$, such that the operator
$$
T_K f(x) = \int f(x-st) K(s,t)\,ds\,dt
$$
is not bounded on $L^2$. In fact, let 
\begin{equation}\label{seven}
\varsigma(s,t) = \phi(s) \phi(t), \text{ where } \phi \in C_c^\infty(B^1(\frac{a}{\sqrt{2}})) \text{ satisfies } \int \phi = 0, \int s \phi(s) \,ds \neq 0,
\end{equation}
then for the distribution $K(s,t) = \sum_{k\in \mathbb{N}} \varsigma^{(2^k, 2^{-k})}(s,t)$, by a change of coordinates, the operator
$$
T_K f(x)= \sum_{k\in \mathbb{N}} \int f(x-st) \varsigma^{(2^k, 2^{-k})}(s,t) \,ds\,dt= \sum_{k\in \mathbb{N}} \int f(x-st) \varsigma (s,t) \,ds\,dt
$$
is not bounded on $L^2$. The above $K(s,t)$ is actually not in $\mathcal{K}(2,e,a,2)$ since it is not compactly supported. But this gives the basic ideas of proving there exists a kernel in $\mathcal{K}(2,e,a,2)$ whose corresponding operator is unbounded, (see Theorem 17.8 of \cite{L2} for a proof, which involves a similar reduction as in Example \ref{know}). Note $\widehat{\mathcal{P}}$ only has zero vector fields, so does $\mathcal{L}(\widehat{\mathcal{P}})$. And since $\big(-\partial_x, (1,1) \big) \in \widehat{\mathcal{N}}$, $\mathcal{L}(\widehat{\mathcal{P}})$ does not control $\widehat{\mathcal{N}}$.
\end{example}

\begin{example}\label{know}
For every $a>0$, there exists $K\in \mathcal{K}(2,e,a,2)$, such that the operator
\begin{equation}\label{jinzhengu}
T_K f(x) = \int f(x-s^3-t^3-st) K(s,t)\,ds\,dt
\end{equation}
is not bounded on $L^2$. In fact, the existence of such a $K$ can be reduced to the existence problem in Example \ref{kitty} in the following way. Suppose by contradiction that there exists $a>0$, such that for every $K \in \mathcal{K}(2,e,a,2)$, the operator $T_K$ in (\ref{jinzhengu}) is bounded on $L^2$. For $K=\sum_{j\in \mathbb{N}^2} \varsigma_j^{(2^j)}$, we can form a Fr\'echet space where $\{\varsigma_j\}_{j\in \mathbb{N}^2}$ lies (see Definition \ref{semi}), and establish a (non-injective) closed linear map $\{\varsigma_j\} \mapsto T_K\in \mathcal{B}(L^2)$ (see Lemma \ref{yangtianxiang} for a similar situation). Then by closed graph theorem, for a fixed $\varsigma(s,t)$ satisfying (\ref{seven}), the operators
$$
T_{L, M}f(x) := \sum_{\substack{k\in \mathbb{N}\\ k\leq M}}\int f(x-s^3 -t^3 -st) \varsigma^{(2^{L+k},  2^{L-k})} (s,t)\,ds\,dt
$$
are bounded on $L^2$ uniformly in $L, M (L \geq M)$. Define a dilation $\Phi_{2^{2L}}(x) := 2^{2L}x$. Then the operators
\begin{align*}
(\Phi_{2^{2L}}^{-1})^*T_{L,M}\Phi_{2^{2L}}^* f(x) &= \sum_{\substack{k\in \mathbb{N}\\ k\leq M}}\int f\big(2^{2L}\big( 2^{-2L}x-s^3 -t^3  -st \big)\big) \varsigma^{(2^{L+k},  2^{L-k})} (s,t)\,ds\,dt\\
& = \sum_{\substack{k\in \mathbb{N}\\ k\leq M}}\int f\big(x-2^{-L} s^3 -2^{-L} t^3  -st \big) \varsigma^{(2^{k},  2^{-k})} (s,t)\,ds\,dt
\end{align*}
are bounded on $L^2$ uniformly in $L, M (L \geq M)$. By the dominated convergence theorem and Fatou's lemma, the operator
$$
\lim_{L\to \infty} (\Phi_{2^{2L}}^{-1})^*T_{L,M}\Phi_{2^{2L}}^* f(x) = \sum_{\substack{k\in \mathbb{N}\\ k\leq M}} \int f(x-st) \varsigma^{(2^{k},  2^{-k})} (s,t)\,ds\,dt  = \sum_{\substack{k\in \mathbb{N}\\ k\leq M}} \int f(x-st) \varsigma (s,t)\,ds\,dt
$$
are bounded on $L^2$ uniformly in $M$; contradiction. Note $\widehat{\mathcal{P}}$ only has two nonzero elements $\big(-\partial_x, (3,0)\big), \big(- \partial_x, (0,3)\big)$, so does $\mathcal{L}(\widehat{\mathcal{P}})$. And since $\big(-\partial_x, (1,1) \big) \in \widehat{\mathcal{N}}$, $\mathcal{L}(\widehat{\mathcal{P}})$ does not control $\widehat{\mathcal{N}}$.
\end{example}

\begin{example}\label{billy}
Consider operators of the form
\begin{equation}\label{hdmi}
T_K f(x) = \int f(x-s-st) K(s,t)\,ds\,dt.
\end{equation}
Note $\widehat{\mathcal{P}}$ only has one nonzero element $\big(-\partial_x, (1,0)\big)$, so does $\mathcal{L}(\widehat{\mathcal{P}})$. $\widehat{\mathcal{N}}$ only consists of one nonzero element $\big(-\partial_x, (1,1) \big)$. Since $\big(-\partial_x, (1,0)\big)$ controls $\big(-\partial_x, (1,1) \big)$, $\mathcal{L}(\widehat{\mathcal{P}})$ controls $\widehat{\mathcal{N}}$. Then by Theorem 5.1 of \cite{LP}, there exists $a>0$, such that for every $K\in \mathcal{K}(2,e,a,2)$, the operator $T_K$ in (\ref{hdmi}) is bounded on $L^2$. 
\end{example}

\begin{remark}
Example \ref{billy} is different from Example \ref{know} in that, when $\mathcal{L}(\widehat{\mathcal{P}})$ controls $\widehat{\mathcal{N}}$, no dilation like $\Phi_{2^{2L}}$ exists to reduce it to Example \ref{kitty}. The hardest part of this paper is to find a ``dilation\footnote{We will use local diffeomorphisms instead of dilations in the general real analytic case; see Theorem \ref{16}} map'' for every case when $\mathcal{L}(\widehat{\mathcal{P}})$ does not control $\widehat{\mathcal{N}}$ so that the problem is reduced to essentially Example \ref{kitty}, (see Section \ref{sntz} for basic ideas of how to choose a dilation for each case).
\end{remark}

In fact, the control notion in the necessary and sufficient condition of Theorem \ref{newly} is reduced in the real line polynomial case as follows. For a polynomial $p(s,t)=\sum_{|\alpha|\geq 1} c_\alpha (s,t)^\alpha$, where $\alpha\in \mathbb{N}^2$ denotes a multi-index, let $0<a,b \leq \infty$ be the smallest $a,b$ such that $c_{(a,0)} \neq 0$, and $c_{(0,b)} \neq 0$, respectively.

\begin{theorem}\label{illuminate} (\cite{17}, see also Theorem 17.8 of \cite{L2} for a proof)
The operator
\begin{equation}\label{line}
Tf(x) = \int f(x-p(s,t)) K(s,t)\,ds\,dt
\end{equation}
is bounded on $L^2$ for every $K$ with sufficiently small support if and only if whenever $c_{(e,f)} \neq 0$, we have $(e,f)$ lies on or above the line passing through $(a,0)$ and $(0,b)$.
\end{theorem}

\section{Heisenberg group and product kernel case}\label{sntz}
Consider in this section operators on the Heisenberg group with product kernels. We will explain our proof methods in this simple special case.

As in the previous section, we use product kernels on $\mathbb{R}^2$. That is we consider $\mathcal{K}(2,e,a,2)$ with $a>0$, $e=\{e_1, e_2\}, e_1 = (1,0), e_2=(0,1)$. 

Consider the Heisenberg group $\mathbb{H}^1 \cong \mathbb{R}^3$. Its Lie algebra has a basis $\{X,Y,T\}$ satisfying
$$
[X,Y] =T, [X,T]=0, [Y,T]=0.
$$
(One way is to take $X=\partial_x -  y\partial_t, Y = \partial_y + x\partial_t, T = 2\partial_t$, where $\mathbb{H}^1 \cong \mathbb{R}^3$ is given coordinates $(x,y,t)$.) The exponential map 
$$
\exp : \text{span}\,\{X,Y,T\} \to \mathbb{H}^1 \cong \mathbb{R}^3, \quad (x,y,t) \mapsto \xi = e^{x X + yY + tT}0
$$ 
is a bijection with nonzero constant Jacobian, and we let $(x,y,t)$ be the coordinates for $\xi$. Then the Lebesgue measure on $\mathbb{H}^1 \cong \mathbb{R}^3$ is equal to a constant multiple of the Lebesgue measure on $\text{span}\,\{X,Y,T\}$ using these coordinates. The group multiplication is $(x_1, y_1, t_1) \cdot (x_2, y_2, t_2) = \big(x_1 + x_2, y_1 + y_2, t_1+t_2 + \frac{1}{2}( x_1y_2 -x_2 y_1)\big)$. Thus Lebesgue measure is a two sided Haar measure on $\mathbb{H}^1$.

Consider functions $\gamma_{(s_1,s_2)}(\xi) = \gamma_s(\xi): \mathbb{R}^2\times \mathbb{H}^1 \to \mathbb{H}^1$ of the form
$$
\gamma_s(\xi) =e^{\sum_{|\alpha|>0} s^\alpha \widehat X_\alpha}\xi,
$$
where $\alpha \in \mathbb{N}^2$ denotes a multi-index, the $\widehat X_\alpha$ are left-invariant vector fields on $\mathbb{H}^1\cong \mathbb{R}^3$, and the sum $\sum_{|\alpha|>0} s^\alpha \widehat X_\alpha$ only has finitely many nonzero terms. In fact, any left-invariant vector field on $\mathbb{H}^1$ is a linear combination of $X, Y, T$ with constant coefficients. Thus we can write
$$
\gamma_s(\xi)= e^{P_1(s) X + P_2(s) Y + P_3(s) T}\xi,
$$
where $P_1(s), P_2(s), P_3(s)$ are polynomials on $\mathbb{R}^2$ without constant terms. Using the above coordinates, we can also write $\gamma_s(\xi) = (P_1(s), P_2(s), P_3(s)) \cdot (x,y,t)$.

We have $\text{deg}\,(\alpha) = \alpha \in \mathbb{N}^2\backslash \{0\}$, and the pure and nonpure power sets of $\gamma_s(\xi)$ are
\begin{align*}
& \widehat{\mathcal{P}} = \{(\widehat X_\alpha, \alpha): \alpha \text{~is nonzero in only one component}\},\\
& \widehat{\mathcal{N}} = \{(\widehat X_\alpha, \alpha): \alpha \text{~is nonzero in both components}\}.
\end{align*}
Since higher order commutators of the $\widehat X_\alpha$ vanish, $\mathcal{L}(\widehat{\mathcal{P}})$ only has finitely many nonzero vector fields, and we denote this finite subset of $\mathcal{L}(\widehat{\mathcal{P}})$ consisting of all the nonzero elements as $\mathcal{L}_0(\widehat{\mathcal{P}})$.

In the case of $\mathbb{H}^1$ and $\mathcal{K}(2,e,a,2)$, Theorem \ref{newly} becomes:

\begin{theorem}\label{clz}
For such $\gamma_s(\xi)$, the following are equivalent.
\begin{enumerate}[label=(\alph*)]
\item\label{BDD}
The operator
\begin{equation}\label{chn}
Tf(\xi) = \int f(\gamma_s(\xi)) K(s)\,ds
\end{equation}
is bounded on $L^p(\mathbb{R}^3)$, for any $1<p<\infty$ and any $K\in \mathcal{K}(2,e,a,2)$ with sufficiently small $a>0$.

\item\label{HYPERPLANE}
For any $(\widehat X_{\alpha_0}, \alpha_0) \in \widehat{\mathcal{N}}$, and any line $\pi$ in $\mathbb{R}^2$ passing through $\alpha_0=(\alpha_0^1, \alpha_0^2)$:
\begin{equation}\label{lines}
\pi =\{(\zeta_1, \zeta_2): b_1 \zeta_1 + b_2 \zeta_2 = b_1 \alpha_0^1 +b_2 \alpha_0^2 \}, \text{~satisfying~} (b_1, b_2)\in [0,\infty)^2 \backslash\{0\},
\end{equation}
we have $
\widehat X_{\alpha_0} \in \text{span}\,H_\pi$, where 
$$
H_\pi=\{\widehat X_\alpha: (\widehat X_\alpha, \alpha) \in \mathcal{L}(\widehat{\mathcal{P}}), \alpha=(\alpha^1, \alpha^2) \text{~satisfies~} b_1 \alpha^1 + b_2 \alpha^2 \leq b_1 \alpha_0^1 + b_2 \alpha_0^2\}.
$$
\end{enumerate}
\end{theorem}

\begin{remark}
The operator $T$ in (\ref{chn}) is not localized by multiplying $\psi$ as in Theorem \ref{newly}. But since $\gamma_s(\cdot)$ commutes with the right group multiplication, $T$ also commutes with the right group multiplication. Thus the boundedness of $T$ follows from the boundedness of the corresponding localized operator.
\end{remark}

We first prove \textbf{\ref{HYPERPLANE} {\boldmath $\Rightarrow$} \ref{BDD}}.
Suppose for any $(\widehat X_{\alpha_0}, \alpha_0) \in \widehat{\mathcal{N}}$, and any line $\pi$ through $\alpha_0$ satisfying (\ref{lines}), we have $\widehat X_{\alpha_0} \in \text{span}\,H_\pi$. Fix an arbitrary $(\widehat X_{\alpha_0}, \alpha_0) \in \widehat{\mathcal{N}}$. Since $\alpha_0$ is nonpure, $\delta^{\alpha_0} \widehat X_{\alpha_0} =0$ if at least one component of $\delta$ equals $0$. For any $\delta \in (0,1]^2$, there exists $\bar \delta \in (0,1]$, $(b_1, b_2) \in [0,\infty)^2 \backslash \{0\}$, such that $({\bar \delta}^{b_1}, {\bar \delta}^{b_2})= \delta$. For this $(b_1, b_2)$, let $\pi = \{(\zeta_1, \zeta_2): b_1 \zeta_1 + b_2 \zeta_2 = b_1 \alpha_0^1 +b_2 \alpha_0^2 \}$. Then there exist some constants $c_\alpha$ such that
$$
\widehat X_{\alpha_0} = \sum_{\substack{(\widehat X_{\alpha}, \alpha) \in \mathcal{L}_0(\widehat{\mathcal{P}})\\ b_1 \alpha^1 + b_2 \alpha^2 \leq b_1 \alpha_0^1 + b_2 \alpha_0^2}} c_\alpha \widehat X_{\alpha}, \quad \text{on~} \mathbb{R}^3.
$$
Thus on $\mathbb{R}^3$,
\begin{align*}
\delta^{\alpha_0} \widehat X_{\alpha_0} = {\bar \delta}^{b_1 \alpha_0^1 + b_2 \alpha_0^2} \widehat X_{\alpha_0} = \sum_{\substack{(\widehat X_{\alpha}, \alpha) \in \mathcal{L}_0(\widehat{\mathcal{P}})\\ b_1 \alpha^1 + b_2 \alpha^2 \leq b_1 \alpha_0^1 + b_2 \alpha_0^2}} \big({\bar \delta}^{b_1 \alpha_0^1 + b_2 \alpha_0^2 - b_1 \alpha^1 - b_2 \alpha^2} c_\alpha \big)  \delta^\alpha \widehat X_{\alpha}.
\end{align*}
Note $b_1 \alpha^1 + b_2 \alpha^2 \leq b_1 \alpha_0^1 + b_2 \alpha_0^2$ implies $ \big|{\bar \delta}^{b_1 \alpha_0^1 + b_2 \alpha_0^2 - b_1 \alpha^1 - b_2 \alpha^2} \big| \leq 1$. Therefore $\mathcal{L}_0(\widehat{\mathcal{P}})$ controls every element of $\widehat{\mathcal{N}}$. Hence $\mathcal{L}(\widehat{\mathcal{P}})$ controls $\widehat{\mathcal{N}}$, which is the condition \ref{CONTROL} in Theorem \ref{newly}. This implies \ref{BDD} by Theorem 5.1 in \cite{LP} (see Section \ref{shaolinchangquan}).

Before proving \ref{BDD} $\Rightarrow$ \ref{HYPERPLANE}, we introduce several notions. Associated to any left-invariant basis $\{X,Y,T\}$ satisfying $[X,Y] =T, [X,T]=0, [Y,T]=0$, and any $\tau_1, \tau_2 >0$, we can define a dilation
$$
\Phi_{\tau_1, \tau_2}: \mathbb{H}^1 \to \mathbb{H}^1, \quad e^{xX+yY +tT}0 \mapsto e^{\tau_1 x X+\tau_2 yY+ \tau_1 \tau_2 tT}0.
$$
This is a diffeomorphism with inverse $\Phi_{\tau_1^{-1}, \tau_2^{-1}}$. Using coordinates $(x,y,t)$, we have the Jacobian
$\det D\Phi_{\tau_1, \tau_2}(x,y,t) = \tau_1^2 \tau_2^2$, a fixed nonzero number. Thus by change of coordinates, for any such basis $\{X,Y,T\}$, any $\tau_1, \tau_2>0$, and any $1<p <\infty$,
\begin{equation}\label{conjugate}
\big\| (\Phi_{\tau_1, \tau_2}^{-1})^* T\Phi_{\tau_1, \tau_2}^* \big\|_{L^p\to L^p} = \big\| T\big\|_{L^p \to L^p},
\end{equation}
where $\Phi^*$ denotes the pullback via $\Phi$. By the Baker-Campbell-Hausdorff formula,
$$
(\Phi_{\tau_1, \tau_2}^{-1})^* T\Phi_{\tau_1, \tau_2}^* f(\xi) = \int f(e^{\tau_1P_1(s) X + \tau_2 P_2(s) Y+ \tau_1 \tau_2 P_3(s) T} \xi) K(s)\,ds.
$$
Recall that using $e=\{(1,0), (0,1)\}$, for $\delta= (\delta_1, \delta_2) \in (0,\infty)^2$, we define the $2$-parameter dilations on $\mathbb{R}^2$ as $\delta s = (\delta_1 s_1, \delta_2 s_2)$, and for a function $\varsigma(s)$ on $\mathbb{R}^2$, we set $\varsigma^{(\delta)}(s_1, s_2) = \delta_1 \delta_2 \varsigma(\delta_1s_1, \delta_2s_2)$. Then $(\delta s)^\alpha = \delta_1^{\alpha^1} \delta_2^{\alpha^2} s^\alpha$, for the multi-index $\alpha = (\alpha^1, \alpha^2)$.

Now we prove \textbf{\ref{BDD} {\boldmath $\Rightarrow$} \ref{HYPERPLANE}}. Assume \ref{BDD} holds. Then there exists $a>0$, such that for any bounded sequence $\{\varsigma_j\}_{j\in \mathbb{N}^2}$ in $C_c^\infty(B^2(a))$ satisfying $\int \varsigma_j(s)\,ds_1 \equiv \int \varsigma_j(s)\,ds_2 \equiv 0$ for every $j\in \mathbb{N}^2$, the operator
$$
T_{\{\varsigma_j\}}f(\xi) := \sum_{j\in \mathbb{N}^2} \int f(\gamma_s(\xi)) \varsigma_j^{(2^j)}(s)\,ds
$$
is bounded on $L^p(\mathbb{R}^3)$ for every $1<p<\infty$. We can form a Fr\'echet space where $\{\varsigma_j\}_{j\in \mathbb{N}^2}$ lies (see Definition \ref{semi}), and establish a (non-injective) closed linear map $\{\varsigma_j\} \mapsto T_{\{\varsigma_j\}} \in \mathcal{B}(L^p)$ (see Lemma \ref{yangtianxiang} for a similar situation).

Suppose by contradiction there exists $(\widehat X_{\alpha_0}, \alpha_0)\in \widehat{\mathcal{N}}$, and a line $\pi=\{(\zeta_1, \zeta_2): b_1 \zeta_1 + b_2 \zeta_2 = b_1 \alpha_0^1 +b_2 \alpha_0^2 \}$ through $\alpha_0$ satisfying (\ref{lines}), such that $\widehat X_{\alpha_0} \not\in \text{span}\,H_\pi$. Fix $b_1, b_2$. Without loss of generality, we can assume every nonpure $\beta$ below $\pi$ satisfies $\widehat X_\beta \in \text{span}\,H_{\pi_\beta}$, where $\beta = (\beta^1, \beta^2)$, and $H_{\pi_\beta} = \{\widehat X_\alpha: (\widehat X_\alpha, \alpha) \in \mathcal{L}(\widehat{\mathcal{P}}), b_1\alpha^1 + b_2\alpha^2 \leq b_1\beta^1 + b_2\beta^2\}$. Otherwise, we take instead $\alpha_0$ to be a nonpure $\beta$ satisfying $\widehat X_\beta \not\in \text{span}\,H_{\pi_\beta}$ with smallest $b_1\beta^1 + b_2\beta^2$, and take $\pi$ to be $\pi_\beta = \{(\zeta_1, \zeta_2): b_1\zeta_1 +b_2\zeta_2 = b_1 \beta^1 + b_2 \beta^2\}$.

Fix $\varsigma \in C_c^\infty(B^2(a))$ satisfying $\int \varsigma(s) \,ds_1 \equiv \int \varsigma(s)\,ds_2\equiv 0$, $\int s^{\alpha_0} \varsigma(s)\,ds \neq 0$ (see existence of such $\varsigma(s)$ in Lemma \ref{phd}), and fix a nonzero vector $\textbf{n} = (\textbf{n}_1, \textbf{n}_2)$ in $\mathbb{R}^2$ perpendicular to $\alpha_0$ and with first coordinate positive. Applying the closed graph theorem to the map $\{\varsigma_j\} \mapsto T_{\{\varsigma_j\}}$, we have for any $1<p<\infty$, the operators
$$
T_{\tau_0,M} f(\xi): = \sum_{\substack{k\in \mathbb{N}\\  1\leq k \leq M}} \int f(\gamma_s(\xi)) \varsigma^{\big((\tau_0^{b_1}, \tau_0^{b_2}) 2^{k \textbf{n}}\big)}(s)\,ds
$$
are bounded on $L^p(\mathbb{R}^3)$ uniformly in $M$ and $\tau_0(\gg  M)$, where $(\tau_0^{b_1}, \tau_0^{b_2})2^{k\textbf{n}} = (\tau_0^{b_1} 2^{k\textbf{n}_1}, \tau_0^{b_2} 2^{k\textbf{n}_2})$. To reach a contradiction, we divide into several cases.

\vspace{1em}
\noindent
\textbf{Case 1: {\boldmath $\widehat{\mathcal{P}}$} only has zero vector fields, and {\boldmath $\widehat{\mathcal{N}}$} has exactly one nonzero vector field.}

This is the simplest case, and every other case eventually reduces to this case. We can write $\gamma_s(\xi) = e^{s^{\alpha_0}\widehat X_{\alpha_0}}\xi$. This case is contained in Case $2$, so we only give a sketch of proof.

Similar to Example \ref{kitty}, for the above $\varsigma(s)$, the operator
$$
\sum_{k\in \mathbb{N}} \int f(e^{s^{\alpha_0}\widehat X_{\alpha_0}} \xi ) \varsigma^{(2^{k\textbf{n}})}(s)\,ds = \sum_{k\in \mathbb{N}} \int f(e^{s^{\alpha_0}\widehat X_{\alpha_0}} \xi ) \varsigma(s)\,ds
$$
is not bounded on $L^p(\mathbb{R}^3)$ for any $1<p<\infty$. Though the distribution $\sum_{k\in \mathbb{N}} \varsigma^{(2^{k\textbf{n}})}(s)$ is not compactly supported, this gives the basic idea for finding appropriate kernels that allow one to reach a contradiction.

\begin{figure}
\centering
\begin{subfigure}{0.3\textwidth}
\centering
\begin{tikzpicture}
\draw (0,2.5) -- (0,0) -- (1.25,0) node[anchor=north]{Case $1 (\subseteq$ Case 2)} -- (2.5,0);
\draw (2,0) node[anchor=south]{$\pi$} -- (0,2);
\draw (2.4, 0.05) -- (2.5,0) -- (2.4, -0.05);
\draw (-0.05, 2.4) -- (0,2.5) -- (0.05, 2.4);
\draw [dashed] (0,0) -- (1,2);
\filldraw[black] (0.667,1.333) circle (1pt) node[anchor=south]{$\alpha_0$};
\draw (0.667, 1.333) -- (1.667, 0.833);
\draw (1.507, 0.833) -- (1.667, 0.833) node[anchor=south]{$\mathbf{n}$} -- (1.587, 0.953);
\end{tikzpicture}
\end{subfigure}
\begin{subfigure}{0.3\textwidth}
\centering
\begin{tikzpicture}
\draw (0,2.5)  -- (0,0) -- (1.25,0) node[anchor=north]{Case 2} -- (2.5,0);
\draw (2,0) node[anchor=south]{$\pi$} -- (0,2) ;
\draw (2.4, 0.05) -- (2.5,0) -- (2.4, -0.05);
\draw (-0.05, 2.4) -- (0,2.5) -- (0.05, 2.4);
\draw [dashed] (0,0) -- (1,2);
\filldraw[black] (0.667,1.333) circle (1pt) node[anchor=south]{$\alpha_0$};
\draw (0.667, 1.333) -- (1.667, 0.833);
\draw (1.507, 0.833) -- (1.667, 0.833) node[anchor=south]{$\mathbf{n}$} -- (1.587, 0.953);
\filldraw[black] (0.6,0) circle (1pt) node[anchor=south]{$\alpha_2$};
\filldraw[black] (0,2) circle (1pt) node[anchor=east]{$\alpha_1$};
\filldraw[black] (1.5, 0.5) circle (1pt) node[anchor=east]{$\alpha_3$};
\filldraw[black] (1.3,1.7) circle (1pt) node[anchor = west]{$\alpha_4$};
\end{tikzpicture}
\end{subfigure}
\begin{subfigure}{0.3\textwidth}
\centering
\begin{tikzpicture}
\draw (0,2.5)  -- (0,0) -- (1.25,0) node[anchor=north]{Case 3} -- (2.5,0);
\draw (2,0) node[anchor=south]{$\pi$} -- (0,2) ;
\draw (2.4, 0.05) -- (2.5,0) -- (2.4, -0.05);
\draw (-0.05, 2.4) -- (0,2.5) -- (0.05, 2.4);
\draw [dashed] (0,0) -- (1,2);
\filldraw[black] (0.667,1.333) circle (1pt) node[anchor=south]{$\alpha_0$};
\draw (0.667, 1.333) -- (1.667, 0.833);
\draw (1.507, 0.833) -- (1.667, 0.833) node[anchor=south]{$\mathbf{n}$} -- (1.587, 0.953);
\filldraw[black] (0.6,0) circle (1pt) node[anchor=south]{$\alpha_2$};
\filldraw[black] (0,2) circle (1pt) node[anchor=east]{$\alpha_1$};
\filldraw[black] (1.5, 0.5) circle (1pt) node[anchor=west]{$\alpha_3$};
\filldraw[black] (1.3,0.33) circle (1pt) node[anchor = north]{$\alpha_4$};
\end{tikzpicture}
\end{subfigure}
\end{figure}

\noindent 
\textbf{Case 2: every nonzero vector field in {\boldmath $\widehat{\mathcal{N}}$} is with degree above or on the line {\boldmath $\pi$}.}

This means every nonzero $(\widehat X_\alpha, \alpha) \in \widehat{\mathcal{N}}$ satisfies $b_1 \alpha^1 + b_2 \alpha^2 \geq b_1 \alpha_0^1 + b_2 \alpha_0^2$, where $\alpha = (\alpha^1,\alpha^2), \alpha_0 = (\alpha_0^1,\alpha_0^2)$. Without loss of generality, we assume $\alpha_0\in \mathbb{N}^2$ has the smallest first coordinate and largest second coordinate among all the nonpure $\alpha$ on the line $\pi$ satisfying $\widehat X_\alpha \not\in \text{span}\,H_\pi$. Then $\textbf{n} \cdot \alpha >0$ for all other nonpure $\alpha$ on the line $\pi$ satisfying $\widehat X_\alpha \not\in \text{span}\,H_\pi$.

\noindent
\textbf{Case 2.1: {\boldmath $\widehat X_{\alpha_0} \in$} span\,{\boldmath $\{H_{\pi}, T\}$}.}

This implies $T\not\in \text{span}\,H_{\pi}$. Thus for any $(\widehat X_{\alpha_1}, \alpha_1), (\widehat X_{\alpha_2}, \alpha_2) \in \widehat{\mathcal{P}}$ satisfying $[\widehat X_{\alpha_1}, \widehat X_{\alpha_2}] \neq 0$, we have 
\begin{equation}\label{inequality}
b_1 (\alpha_1^1 + \alpha_2^1) + b_2 (\alpha_1^2+ \alpha_2^2) > b_1 \alpha_0^1 +b_2 \alpha_0^2,
\end{equation}
where $\alpha_1=(\alpha_1^1,\alpha_1^2), \alpha_2=(\alpha_2^1,\alpha_2^2)$. We now change the basis $\{X,Y,T\}$, to implement appropriate dilations.

\begin{itemize}
\item 
If $\dim \text{span}\,H_\pi =0$, let $\widetilde X=X, \widetilde Y= Y,\widetilde T=T$.

\item
If $\dim \text{span}\,H_\pi =1$, let $\widetilde X$ be a nonzero element $\widehat X_\alpha$ in $H_\pi$ with the smallest $b_1 \alpha^1 + b_2 \alpha^2$, and find left-invariant nonzero vector fields $\widetilde Y, \widetilde T$ such that $[\widetilde X,\widetilde Y]=\widetilde T, [\widetilde X,\widetilde T]=[\widetilde Y, \widetilde T]=0$.

\item
If $\dim \text{span}\,H_\pi =2$, let $\widetilde X$ be a nonzero element $\widehat X_\alpha$ in $H_\pi$ with the smallest $b_1 \alpha^1 + b_2 \alpha^2$. Since $T\not\in \text{span}\,H_\pi$, there exist elements $\widehat X_\beta$ in $H_\pi$ not commuting with $\widetilde X$. Let $\widetilde Y$ be a such $\widehat X_\beta$ with the smallest $b_1 \beta^1 + b_2 \beta^2$, where $\beta = (\beta^1, \beta^2)$. Let $\widetilde T=[\widetilde X,\widetilde Y]$. 
\end{itemize}

Changing to the new basis $\{\widetilde X, \widetilde Y, \widetilde T\}$, there exist polynomials $Q_1(s), Q_2(s), Q_3(s)$ without constant terms such that
$$
\gamma_s(\xi) = e^{Q_1(s) \widetilde X + Q_2(s) \widetilde Y + Q_3(s) \widetilde T}\xi.
$$
For $i=1,2$, let $\theta_i = \infty$ if $Q_i \equiv 0$, otherwise let 
$$
\theta_i = \min \{b_1 \alpha^1 + b_2 \alpha^2: s^\alpha \text{~is a term in~} Q_i(s) \}.
$$
Then by the way we choose the basis and by (\ref{inequality}), we have $\theta_1 + \theta_2> b_1 \alpha_0^1 + b_2 \alpha_0^2$.

Consider dilations $\Phi_{\tau_1, \tau_2}$ associated to this new basis $\{\widetilde X, \widetilde Y, \widetilde T\}$. Set $\tau_1= \tau_0^{(b_1 \alpha_0^1 + b_2 \alpha_0^2) \frac{\theta_1}{\theta_1 + \theta_2}}$, $\tau_2 =\tau_0^{(b_1 \alpha_0^1 + b_2 \alpha_0^2) \frac{\theta_2}{\theta_1 + \theta_2}}$. By (\ref{conjugate}), the operators
\begin{align*}
&\quad (\Phi_{\tau_1, \tau_2}^{-1})^* T_{\tau_0, M} \Phi_{\tau_1, \tau_2}^* f(\xi)\\
& = \sum_{\substack{k \in \mathbb{N}\\ 1\leq k \leq M}} \int f(e^{\tau_1Q_1(s) \widetilde X + \tau_2 Q_2(s) \widetilde Y+ \tau_1 \tau_2 Q_3(s) \widetilde T} \xi) \varsigma^{\big((\tau_0^{b_1}, \tau_0^{b_2}) 2^{k \textbf{n}}\big)}(s)\,ds\\
&= \sum_{\substack{k \in \mathbb{N}\\ 1\leq k \leq M}} \int f\big(e^{\tau_1Q_1\big((\tau_0^{-b_1}, \tau_0^{-b_2})s\big) \widetilde X + \tau_2 Q_2\big((\tau_0^{-b_1}, \tau_0^{-b_2})s\big) \widetilde Y+ \tau_1 \tau_2 Q_3\big((\tau_0^{-b_1}, \tau_0^{-b_2})s\big) \widetilde T} \xi\big) \varsigma^{(2^{k\textbf{n}})}(s)\,ds
\end{align*}
are bounded on $L^p(\mathbb{R}^3)$ uniformly in $\tau_0 (\gg  M), M$, for all $1<p<\infty$. Note $H_\pi \subseteq \text{span}\,\{\widetilde X,\widetilde Y\}$. And since $\tau_1 \ll \tau_0^{\theta_1}$ and $\tau_2 \ll \tau_0^{\theta_2}$ when $\tau_0\gg1$, there exist constants $c_\alpha$ such that
$$
\lim_{\tau_0 \to \infty} \tau_1 \tau_2 Q_3\big((\tau_0^{-b_1}, \tau_0^{-b_2}) s\big) = \sum_{\substack{\text{~nonpure~} \alpha \text{~on~} \pi\\ \widehat X_\alpha \not\in \text{span}\, H_\pi}} c_\alpha s^\alpha, \quad \lim_{\tau_0 \to \infty} \tau_iQ_i\big((\tau_0^{-b_1}, \tau_0^{-b_2}) s\big) = 0 \text{~for~} i=1,2.
$$

By the dominated convergence theorem and Fatou's lemma, for every $1<p<\infty$, the operators
\begin{align*}
\lim_{\tau_0 \to \infty} (\Phi_{\tau_1, \tau_2}^{-1})^* T_{\tau_0, M} \Phi_{\tau_1, \tau_2}^* f(\xi) = \sum_{\substack{k \in \mathbb{N}\\ 1\leq k \leq M}} \int f(e^{\sum c_\alpha s^\alpha \widetilde T} \xi) \varsigma^{(2^{k\textbf{n}})}(s)\,ds
\end{align*}
are bounded on $L^p(\mathbb{R}^3)$ uniformly in $M$, where the $c_\alpha$ are nonzero only for the nonpure $\alpha$'s on $\pi$ satisfying $\widehat X_\alpha \not\in \text{span}\,H_\pi$. Note $\textbf{n} \cdot \alpha_0=0$, and $\textbf{n} \cdot \alpha >0$ for all the other nonpure $\alpha$'s on $\pi$ satisfying $\widehat X_\alpha \not\in \text{span}\,H_\pi$. Thus by the dominated convergence theorem, for every $f\in C_c^\infty(\mathbb{R}^3)$,
\begin{align*}
\lim_{k \to \infty} \int f(e^{\sum c_{\alpha} s^{\alpha}\widetilde T}\xi) \varsigma^{(2^{k \textbf{n}})}(s)\,ds = \lim_{k \to \infty} \int f(e^{\sum c_{\alpha} (2^{-k \textbf{n}}s)^{\alpha}\widetilde T}\xi) \varsigma(s)\,ds  = \int f(e^{c_{\alpha_0} s^{\alpha_0} \widetilde T}\xi)\varsigma(s)\,ds.
\end{align*}
Applying the dominated convergence theorem again, the function $\int f(e^{\sum c_{\alpha} s^{\alpha}\widetilde T}\xi) \varsigma^{(2^{k \textbf{n}})}(s)\,ds$ converges to $\int f(e^{c_{\alpha_0} s^{\alpha_0} \widetilde T}\xi)\varsigma(s)\,ds$ in $L^p(\mathbb{R}^3)$, and thus $\int f(e^{c_{\alpha_0} s^{\alpha_0} \widetilde T}\xi)\varsigma(s)\,ds \equiv 0$, for every $f\in C_c^\infty(\mathbb{R}^3)$.

But since $\int s^{\alpha_0} \varsigma(s)\,ds \neq 0$, there exists $f\in C_c^\infty(\mathbb{R}^3)$ such that the function $\int f(e^{c_{\alpha_0} s^{\alpha_0} \widetilde T}\xi)\varsigma(s)\,ds$ is nonzero on some neighborhood of $0$; contradiction.

\noindent
\textbf{Case 2.2: {\boldmath $\widehat X_{\alpha_0} \not\in$} span\,{\boldmath $\{H_\pi, T\}$}.}

We first change the basis $\{X,Y,T\}$, to implement appropriate dilations. 
\begin{itemize}
\item
If $H_\pi\subseteq \text{span}\, T$, let $\widetilde Y=\widehat X_{\alpha_0}$, and find left-invariant nonzero vector fields $\widetilde X,\widetilde T$ satisfying $[\widetilde X,\widetilde Y]=\widetilde T,[\widetilde X,\widetilde T]=[\widetilde Y,\widetilde T]=0$.

\item 
If $H_\pi \not\subseteq \text{span}\,T$, let $\widetilde Y=\widehat X_{\alpha_0}$, and let $\widetilde X$ be an element $\widehat X_\alpha$ in $H_\pi$ not parallel to $T$ with the smallest $b_1 \alpha^1 +b_2 \alpha^2$, and let $\widetilde T=[\widetilde X,\widetilde Y]$.
\end{itemize}
Changing to the new basis $\{\widetilde X, \widetilde Y, \widetilde T\}$, there exist polynomials $Q_1(s), Q_2(s), Q_3(s)$ without constant terms such that
$$
\gamma_s(\xi) = e^{Q_1(s) \widetilde X + Q_2(s) \widetilde Y + Q_3(s) \widetilde T}\xi.
$$
Consider dilations $\Phi_{\tau_1, \tau_2}$ associated to the new basis $\{\widetilde X, \widetilde Y, \widetilde T\}$. Since $\widetilde T$ is parallel to $T$, $H_\pi \subseteq \text{span}\,\{\widetilde X,\widetilde T\}$. By setting $\tau_1= \tau_0^{-(b_1 \alpha_0^1 + b_2 \alpha_0^2)}, \tau_2 =\tau_0^{b_1 \alpha_0^1 + b_2 \alpha_0^2}$, we can write
\begin{align*}
&\lim_{\tau_0 \to \infty} \tau_2 Q_2\big((\tau_0^{-b_1}, \tau_0^{-b_2}) s\big) = \sum_{\substack{\text{~nonpure~} \alpha \text{~on~} \pi\\ \widehat X_\alpha \not\in \text{span}\, H_\pi}} c_\alpha s^\alpha, \text{ for some constants } c_\alpha,\\
&\lim_{\tau_0 \to \infty} \tau_1 Q_1\big((\tau_0^{-b_1}, \tau_0^{-b_2}) s\big) = \lim_{\tau_0 \to \infty} \tau_1 \tau_2 Q_3\big((\tau_0^{-b_1}, \tau_0^{-b_2}) s\big) =0.
\end{align*}

Similar to Case 2.1, for every $1<p<\infty$, the operators
\begin{align*}
\lim_{\tau_0\to \infty} (\Phi_{\tau_1, \tau_2}^{-1})^* T_{\tau_0, M} \Phi_{\tau_1, \tau_2}^* f(\xi) = \sum_{\substack{k \in \mathbb{N}\\ 1\leq k \leq M}} \int f(e^{\sum c_\alpha s^\alpha \widetilde Y} \xi) \varsigma^{(2^{k \textbf{n}})}(s)\,ds
\end{align*}
are bounded on $L^p(\mathbb{R}^3)$ uniformly in $M$, where the $c_\alpha$ are nonzero only for the nonpure $\alpha$'s on $\pi$ satisfying $\widehat X_\alpha \not\in \text{span}\,H_\pi$. Since $\textbf{n} \cdot \alpha_0=0$, and $\textbf{n} \cdot \alpha >0$ for all the other nonpure $\alpha$'s on $\pi$ satisfying $\widehat X_\alpha \not\in \text{span}\,H_\pi$, for every $f\in C_c^\infty(\mathbb{R}^3)$, by the dominated convergence theorem,
\begin{align*}
\lim_{k\to \infty} \int f(e^{\sum c_\alpha s^\alpha \widetilde Y} \xi) \varsigma^{(2^{k\textbf{n}})}(s)\,ds  = \int f(e^{c_{\alpha_0} s^{\alpha_0} \widetilde Y}\xi)\varsigma(s)\,ds.
\end{align*}
Thus we reach a contradiction by finding $f\in C_c^\infty(\mathbb{R}^3)$ such that the function $\int f(e^{c_{\alpha_0} s^{\alpha_0} \widetilde Y}\xi)\varsigma(s)\,ds$ is nonzero on some neighborhood of $0$.

\noindent
\textbf{Case 3: there exists a nonzero vector field in {\boldmath $\widehat{\mathcal{N}}$} with degree below the line {\boldmath $\pi$}. }

This means there exists a nonzero $(\widehat X_\beta, \beta)\in \widehat{\mathcal{N}}$ satisfying $b_1 \beta^1 + b_2 \beta^2 < b_1 \alpha_0^1 + b_2 \alpha_0^2$. But every nonpure $\beta$ below $\pi$ satisfies $\widehat X_{\beta} \in$ span\,$H_{\pi_\beta}$, where  $H_{\pi_\beta} = \{\widehat X_\alpha: (\widehat X_\alpha, \alpha)\in \mathcal{L}(\widehat{\mathcal{P}}), b_1 \alpha^1 + b_2 \alpha^2 \leq b_1 \beta^1 + b_2 \beta^2 \}$. Therefore this case is dealt with in a manner similar to Case 2, where the limiting process gets rid of the term $\widehat X_{\beta}$ when getting rid of every $\widehat X_\alpha \in H_{\pi_\beta}$.

\begin{remark}\label{explode}
\begin{enumerate}
\item
A key property of $\mathbb{H}^1$ used in proving \ref{BDD} $\Rightarrow$ \ref{HYPERPLANE} is that we only have first order commutators (higher order commutators vanish), which is slightly more complicated than the real line case in Section \ref{real}, where all the commutators vanish. It is easy to generalize the above proof method to the case of Nilpotent groups, where we only have bounded order commutators. For the general real analytic case, the corresponding property is the Noetherian property; see Section \ref{boruozhang}.

\item 
For the general real analytic case, we will use local diffeomorphisms instead of the global dilations $\Phi_{\tau_1, \tau_2}$; see the map $\Phi$ in Theorem \ref{16}.

\item
The above proof of $(a) \Rightarrow (b)$ gives an outline for proving the general real analytic case, where the entire limiting process above need to be applied twice for the general case in Sections \ref{A} and \ref{xj}, how we choose $\alpha_0$ among nonpure powers corresponds to first halves of Sections \ref{A} and \ref{xj}, taking limit of a sequence of uniformly bounded operators corresponds to second halves of Sections \ref{A} and \ref{xj}, and proving that the limit operator is unbounded as in Case 1 corresponds to Section \ref{as}.
\end{enumerate}
\end{remark}

\section{The role of real analyticity}\label{boruozhang}
We introduce in this section the Noetherian property of real analytic functions and several related results, which are key in proving the conditions given in Theorems \ref{fy}, \ref{newly} are necessary.

\begin{definition}
Let 
$$
\mathcal{A}_n = \{ f: \mathbb{R}^n_0 \to \mathbb{R} \text{ such that } f \text{~is real analytic}\},
$$
the set of germs of real analytic functions defined on a neighborhood of $0$ in $\mathbb{R}^n$. Note $\mathcal{A}_n$ is a ring. And let $$
\mathcal{A}_n^m = \{ f: \mathbb{R}^n_0 \to \mathbb{R}^m \text{ such that } f \text{~is real analytic}\}.
$$
Note $\mathcal{A}_n^m$ is a free $\mathcal{A}_n$-module.
\end{definition}

The next classical result is the key property we use about real analytic functions.

\begin{theorem}\label{pianzhi} 
$\mathcal{A}_n$ is a Noetherian ring, and $\mathcal{A}_n^m$ is a Noetherian $\mathcal{A}_n$-module. 
\end{theorem}

\begin{proof}[Comments on the proof]
This is a simple consequence of the Weierstrass preparation theorem. See page
148 of \cite{ZS75}. The proof in \cite{ZS75} is for the formal power series ring, however, as mentioned on page
130 of \cite{ZS75}, the proof also works for the ring of convergent power series: i.e., the ring of power series
with some positive radius of convergence. The ring of germs of real analytic functions is isomorphic to
the ring of convergent power series. Besides, for any Noetherian ring $R$, the $R$-module $R^m$ is Noetherian. 
\end{proof}

To better adapt the Noetherian property to a set of real analytic vector fields paired with degrees, we need the following notions of strong control. Suppose $(X_0, d_0), (X_1, d_1), \ldots, (X_q,d_q)$ are smooth vector fields defined on a neighborhood $\Omega$ of $0$ in $\mathbb{R}^n$, each paired with a $\nu$-parameter degree, and $W(t,x)$ is a smooth function mapping from $B^N(\rho)\times \Omega$ into $\mathbb{R}^n$ for some $\rho>0$, satisfying $W(0,x)\equiv 0$. Denote $(\mathbf{X},\mathbf{d}) = \{(X_1, d_1), \ldots, (X_q, d_q)\}$.

\begin{definition}\label{**}
We say $(\mathbf{X},\mathbf{d})$ strongly controls $(X_0,d_0)$, if on some neighborhood $U \subseteq \Omega$ of $0$,
$$
X_0 = \sum_{\substack{1\leq l \leq q \\ d_l \leq d_0}}c_{l}X_l,
$$
where the $c_l$ are $C^\infty$ on some neighborhood of the closure of $U$, and the inequalities $d_l\leq d_0$ hold coordinate-wise.

When we wish to make the choice of $U$ explicit, we say, $(\mathbf{X},\mathbf{d})$ strongly controls $(X_0,d_0)$ on $U$.
\end{definition}

\begin{definition}\label{*W}
We say $(\mathbf{X},\mathbf{d})$ strongly controls $W(t,x)$, if there exists $0<\rho' \leq \rho$ and a neighborhood $U \subseteq \Omega$ of $0$, such that
$$
W(t,x) = \sum_{\substack{1\leq l \leq q\\ \text{deg}\,(\alpha) = d_l}} c_{l,\alpha}(t,x) t^{\alpha}X_l(x), \text{~on~} B^N(\rho')\times U,
$$
where the $c_{l, \alpha}(t,x)$ are $C^\infty$ on some neighborhood of the closure of $B^N(\rho')\times U$.

When we wish to make the choice of $\rho'$ and $U$ explicit, we say $(\mathbf{X},\mathbf{d})$ strongly controls $W(t,x)$ on $B^N(\rho')\times U$.
\end{definition}

Theorem \ref{pianzhi} is related to the following ``finite-type'' conditions.

\begin{theorem}\label{K} (Theorems 9.1, 9.2 in \cite{ANALYTIC}; see also \cite{LOB70} and Theorem 1.2.5 of \cite{GAL79})
\begin{enumerate}[label=(\roman*)]
\item\label{finite-type}
For $\mathcal{S}$ a set of real analytic vector fields on a common neighborhood of $0$ in $\mathbb{R}^n$ each paired with a $\nu$-parameter degree, there exists a finite subset $\mathcal{F}$ of $\mathcal{S}$ that strongly controls\footnote{The particular neighborhood $U$ used for strong control (in Definition \ref{**}) may depend on each element of $\mathcal{S}$.} every element of $\mathcal{S}$. Moreover, the union of such $\mathcal{F}$ with any finite subset of $\mathcal{S}$ still strongly controls every element of $\mathcal{S}$. In particular, given any finite subset $\mathcal{S}_0$ of $\mathcal{S}$, there exists a finite subset of $\mathcal{S}$ that contains $\mathcal{S}_0$ and strongly controls every element of $\mathcal{S}$.

\item
If the $\gamma_t(x)$ is real analytic, then there exists a finite subset $\mathcal{F}'$ of $\mathcal{L}(\mathcal{P}\cup \mathcal{N})$ that strongly controls $W(t,x)$. Moreover, the union of such $\mathcal{F}'$ with any finite subset of $\mathcal{L}(\mathcal{P} \cup \mathcal{N})$ still strongly controls $W(t,x)$. In particular, given any finite subset $\mathcal{S}'$ of $\mathcal{L}(\mathcal{P} \cup \mathcal{N})$, there exists a finite subset of $\mathcal{L}(\mathcal{P}\cup \mathcal{N})$ that contains $\mathcal{S}'$ and strongly controls $W(t,x)$.
\end{enumerate}
\end{theorem}

\begin{remark}
\begin{enumerate}
\item 
$\gamma_t(x)$ is real analytic on a neighborhood of $(0,0)$ in $\mathbb{R}^N \times \mathbb{R}^n$, if and only if the corresponding $W(t,x)$ is real analytic on a neighborhood of $(0,0)$ in $\mathbb{R}^N \times \mathbb{R}^n$ (see the relation between $\gamma$ and $W$ in Lemma 2.4, and see Theorem 10.7.5 in \cite{18}).

\item
If $W$ is $C^\infty$ but not real analytic, then its $\mathcal{L}(\mathcal{P}\cup  \mathcal{N})$ might not satisfy the above ``finite-type'' conditions. $W(t,x,y) =t \partial_x + t^2 e^{-\frac{1}{x^2}} \partial_y$ is an example whose $\mathcal{L}(\mathcal{P}\cup \mathcal{N})$ does not satisfy the above condition \ref{finite-type}. We have
$$
\mathcal{L}(\mathcal{P}\cup  \mathcal{N}) = \Big\{ \big(\partial_x,1\big), \big(e^{-\frac{1}{x^2}}\partial_y,2\big), \Big(\frac{2}{x^3} e^{-\frac{1}{x^2}}\partial_y, 3\Big), \Big(\big(\frac{4}{x^6} - \frac{6}{x^4}\big) e^{-\frac{1}{x^2}}\partial_y,4\Big), \ldots \Big\}.
$$ 
Observe that every vector field in $\mathcal{L}(\mathcal{P} \cup \mathcal{N})$ is paired with a scalar degree, and is not in the $C^\infty$-module generated by the vector fields in $\mathcal{L}(\mathcal{P} \cup \mathcal{N})$ paired with smaller degrees, on any neighborhood of $0$. Therefore there does not exist a finite subset of $\mathcal{L}(\mathcal{P} \cup \mathcal{N})$ that strongly controls every element of $\mathcal{L}(\mathcal{P}\cup \mathcal{N})$.

\item
In the multi-parameter case, H\"ormander's condition \footnote{We say a set of $C^\infty$ vector fields $\mathcal{S}$ defined on a common neighborhood of $0\in \mathbb{R}^n$ satisfies the H\"ormander's condition at $0$, if the Lie algebra generated by the vector fields in $\mathcal{S}$ spans the tangent space of $\mathbb{R}^n$ at $0$ (and consequently spans the tangent space of $\mathbb{R}^n$ in a neighborhood of $0$).} does not imply the ``finite-type'' condition \ref{finite-type}, (although it does in the single-parameter case; see Lemma 3.2 in \cite{L2}). Let $e= \{(1,0), (0,1)\}$ and consider $W(t_1, t_2, x,y) = t_1 \partial_x + t_1^2 e^{-\frac{1}{x^2}} \partial_y + t_2 \partial_y$. 
\begin{align*}
\mathcal{L}(\mathcal{P} \cup \mathcal{N})= &\Big\{ \big(\partial_y, (0,1)\big), \big(\partial_x,(1,0)\big), \big(e^{-\frac{1}{x^2}}\partial_y,(2,0)\big), \\
&\quad \Big(\frac{2}{x^3} e^{-\frac{1}{x^2}}\partial_y, (3,0)\Big), \Big(\big(\frac{4}{x^6} - \frac{6}{x^4}\big) e^{-\frac{1}{x^2}}\partial_y,(4,0)\Big), \ldots \Big\}
\end{align*}
satisfies H\"ormander's condition, but there does not exist a finite subset that strongly controls every element of $\mathcal{L}(\mathcal{P} \cup \mathcal{N})$.
\end{enumerate}
\end{remark}

We have defined the control notion for an $(X_0,d_0)$ and the notions of strong control for both $(X_0,d_0)$ and $W(t,x)$ (Definitions \ref{controldefn}, \ref{**}, and \ref{*W}). Now we define the control notion for $W(t,x)$, which is needed in Section \ref{shaolinchangquan}. Recall that for every $\xi>0$, we let $\Vec{\xi}:= (\xi, \ldots, \xi) \in (0,\infty)^\nu$.

\begin{definition}
We say $(\mathbf{X},\mathbf{d})$ controls $W$, if there exists $0<\rho' \leq \rho$, $\xi_1>0$, and a neighborhood $U \subseteq \Omega$ of $0$, such that $(\mathbf{X},\mathbf{d})$ satisfies $\mathcal{C}(x, \Vec{\xi_1}, \Omega)$ for all $x\in U$, and that for all $x\in U, \delta\in [0,1]^\nu$, 
$$
W(\delta t,u) = \sum_{l=1}^q c_l^{x,\delta}(t,u) \delta^{d_l}X_l(u), \text{~for~}(t,u) \in B^N(\rho')\times B_{(\mathbf{X},\mathbf{d})}(x,\xi_1\delta),
$$
and
$$
\sup_{\substack{\delta\in [0,1]^\nu\\ x\in U}} \sum_{|\beta'|+|\beta| \leq m} \Big\| \partial_t^\beta (\delta \mathbf{X})^{\beta'} c_l^{x,\delta} \Big\|_{C(B^N(\rho')\times B_{(\mathbf{X},\mathbf{d})}(x,\xi_1\delta))} <\infty, \text{~for every~} m\in \mathbb{N}.
$$

When we wish to make the choice of $\rho', U, \xi_1$ explicit, we say, $(\mathbf{X},\mathbf{d})$ controls $W$ on $U$ for $t\in B^N(\rho')$ with parameter $\xi_1$.
\end{definition}

Indeed strong control implies control, see Lemma \ref{4}.

\section{Outline of the proof}\label{shaolinchangquan}
A more complete version of the main result (Theorem \ref{fy}) is:

\begin{theorem}\label{lieri}
For $N, \nu, e$ satisfying (\ref{e}), and for $n\in \mathbb{N}_{>0}$, let $\gamma_t(x)$ be a real analytic function of $(t,x)$ defined on a neighborhood of $(0,0)$ in $\mathbb{R}^N\times \mathbb{R}^n$ with $\gamma_0(x)\equiv x$. The following are equivalent.

\begin{enumerate}[label=(\alph*)]
\item\label{WTX}
$\mathcal{L}(\mathcal{P})$ controls $\mathcal{N}$.

\item\label{GAMMA}
$\mathcal{L}(\widehat{\mathcal{P}})$ controls $\widehat{\mathcal{N}}$.

\item\label{ANYP}
For every $1< p< \infty$, the operator
$$
Tf(x) = \psi(x)\int f(\gamma_t(x))K(t)\,dt
$$
is bounded on $L^p(\mathbb{R}^n)$, for every $\psi \in C_c^\infty(\mathbb{R}^n)$ with sufficiently small support, and every $K \in \mathcal{K}(N,e,a,\nu)$ with sufficiently small $a>0$.

\item\label{SOMEP}
There exists $1<p< \infty$, such that the operator
$$
Tf(x) = \psi(x)\int f(\gamma_t(x))K(t)\,dt
$$
is bounded on $L^p(\mathbb{R}^n)$, for every $\psi \in C_c^\infty(\mathbb{R}^n)$ with sufficiently small support, and every $K \in \mathcal{K}(N,e,a,\nu)$ with sufficiently small $a>0$.
\end{enumerate}
\end{theorem}

\textbf{\ref{WTX} {\boldmath $\Rightarrow$} \ref{ANYP}:} suppose $\mathcal{L}(\mathcal{P})$ controls $\mathcal{N}$. Then by Lemma \ref{yicu}, there exists a finite subset $\mathcal{F} \subseteq \mathcal{L}(\mathcal{P})$ that controls every element of $\mathcal{L}(\mathcal{P}\cup \mathcal{N})$. By Theorem \ref{K}, there exists a finite subset $\mathcal{F}' \subseteq \mathcal{L}(\mathcal{P} \cup \mathcal{N})$ that contains $\mathcal{F}$ and strongly controls $W$. Note $\mathcal{F}$ controls every element of $\mathcal{F}'$ and controls every element of $\{([X_i,X_j], d_i+d_j): X_i, X_j \in \mathcal{F} \}$. By Lemma \ref{4} and Lemma \ref{5}, $\mathcal{F}$ controls $W$. Then by Theorem 5.1 in \cite{LP}, \ref{ANYP} holds.

\textbf{\ref{ANYP} {\boldmath $\Rightarrow$} \ref{SOMEP}:} trivial.

\textbf{\ref{WTX} {\boldmath $\Leftrightarrow$} \ref{GAMMA}:} see Corollary \ref{EQUIV}.

\textbf{\ref{SOMEP} {\boldmath $\Rightarrow$} \ref{WTX}:} the main goal of this paper. Suppose \ref{SOMEP} holds but \ref{WTX} does not hold. The following three theorems will lead to a contradiction.

\vspace{1em}
Assume there exists such a real analytic $\gamma_t(x)$ that satisfies \ref{SOMEP} but does not satisfy \ref{WTX}. This implies for this $\gamma_t(x)$, we have $\nu>1$, and by Theorem \ref{K}, there exists a finite subset $\mathcal{F} \subseteq \mathcal{L}(\mathcal{P})$ that strongly controls every element of $\mathcal{L}(\mathcal{P})$, but at least one element of $\mathcal{L}(\mathcal{P} \cup \mathcal{N})$ is not controlled by $\mathcal{F}$. The following theorem reduces this to the case that all the vector fields associated with pure powers vanish.

\begin{proposition}\label{1}
Suppose there exist $N, \nu>1, e$ satisfying (\ref{e}), $p\in (1,\infty), n\in \mathbb{N}_{>0}$, and a real analytic function $\gamma_t(x)$ defined on a neighborhood of $(t,x)=(0,0)$ in $\mathbb{R}^N\times \mathbb{R}^n$ with $\gamma_0(x)\equiv x$, such that the operator
$$
Tf(x) = \psi(x)\int f(\gamma_t(x))K(t)\,dt
$$
is bounded on $L^p(\mathbb{R}^n)$ for every $\psi\in C_c^\infty(\mathbb{R}^n)$ with sufficiently small support, and every $K \in \mathcal{K}(N,e,a,\nu)$ with sufficiently small $a>0$, and such that there exists a finite subset $\mathcal{F} \subseteq \mathcal{L}(\mathcal{P})$ that strongly controls every element of $\mathcal{L}(\mathcal{P})$, but at least one element of $\mathcal{L}(\mathcal{P} \cup \mathcal{N})$ is not controlled by $\mathcal{F}$.

Then there exist $\bar N, \bar \nu>1, \bar e$ satisfying (\ref{e}), $\bar n\in \mathbb{N}_{>0}$, and a $C^\infty$ function $\bar \gamma_t(x)$ defined on a neighborhood of $(t,x)=(0,0)$ in $\mathbb{R}^{\bar N}\times \mathbb{R}^{\bar n}$  with $\bar \gamma_0(x)\equiv x$, such that the operator
$$
\overline Tf(x) = \psi(x)\int f(\bar \gamma_t(x))K(t)\,dt
$$
is bounded on $L^p(\mathbb{R}^{\bar n})$ for every $\psi\in C_c^\infty(\mathbb{R}^{\bar n})$ with sufficiently small support, and every $K \in \mathcal{K}(\bar N,\bar e,\bar a,\bar \nu)$ with sufficiently small $\bar a>0$, and such that for the corresponding $\overline W(t,x)$ and its pure and nonpure power sets $\overline{\mathcal{P}}, \overline{\mathcal{N}}$, there exists a finite subset $\overline{\mathcal{F}}'\subseteq \mathcal{L}(\overline{\mathcal{P}}\cup\overline{\mathcal{N}})$ that strongly controls $\overline W$ and every element of $\mathcal{L}(\overline{\mathcal{P}}\cup \overline{\mathcal{N}})$, and at least one vector field in $\overline{\mathcal{F}}'$ is nonzero at a sequence of points convergent to $0$ in $\mathbb{R}^{\bar n}$, but every vector field in $\overline{\mathcal{P}}$ vanishes on a common neighborhood of $0$ in $\mathbb{R}^{\bar n}$.
\end{proposition}

The next theorem further reduces the situation to the case that the only non-vanishing vector fields are associated with the same nonpure degree.

\begin{proposition}\label{2}
Suppose there exist $N, \nu>1, e$ satisfying (\ref{e}), $p\in (1,\infty), n\in \mathbb{N}_{>0}$, and a $C^\infty$ function $\gamma_t(x)$ defined on a neighborhood of $(t,x)=(0,0)$ in $\mathbb{R}^{N}\times \mathbb{R}^{n}$ with $\gamma_0(x)\equiv x$, such that the operator
$$
Tf(x) = \psi(x)\int f(\gamma_t(x))K(t)\,dt
$$
is bounded on $L^p(\mathbb{R}^{n})$ for every $\psi\in C_c^\infty(\mathbb{R}^{n})$ with sufficiently small support, and every $K \in \mathcal{K}(N, e, a, \nu)$ with sufficiently small $a>0$, and such that there exists a finite subset $\mathcal{F}'\subseteq \mathcal{L}(\mathcal{P}\cup \mathcal{N})$ that strongly controls $W$ and every element of $\mathcal{L}(\mathcal{P}\cup \mathcal{N})$, and at least one vector field in $\mathcal{F}'$ is nonzero at a sequence of points convergent to $0$ in $\mathbb{R}^n$, but every vector field in $\mathcal{P}$ vanishes on a common neighborhood of $0$ in $\mathbb{R}^n$.

Then there exist $\bar n\in \mathbb{N}_{>0}$ and a $C^\infty$ function $\bar \gamma_t(x)$ defined on a neighborhood of $(t,x)=(0,0)$ in $\mathbb{R}^{N}\times \mathbb{R}^{\bar n}$ with $\bar \gamma_0(x)\equiv x$, such that the operator
$$
\overline Tf(x) = \psi(x)\int f(\bar \gamma_t(x))K(t)\,dt
$$
is bounded on $L^p(\mathbb{R}^{\bar n})$ for every $\psi\in C_c^\infty(\mathbb{R}^{\bar n})$ with sufficiently small support, and every $K \in \mathcal{K}(N,e,\bar a,\nu)$ with sufficiently small $\bar a>0$, and such that the corresponding $\overline{W}(t,x)$ can be written as a finite sum on some neighborhood of $(0,0)$:
$$
\overline W(t,x) =\sum_{l=1}^Q t^{\alpha_l} \overline{X_l}(x),
$$ 
where the $\alpha_l\in \mathbb{N}^N$ are such that $\text{deg}\,(\alpha_1) = \cdots = \text{deg}\,(\alpha_Q)\in [0,\infty)^\nu$ is a nonpure degree, the $\overline{X}_l$ commute with each other, and there exists $1\leq l\leq Q$ such that $\overline{X}_l(0) \neq 0$.
\end{proposition}

Roughly speaking, the following theorem shows that if all the non-vanishing vector fields are associated with the same nonpure degree, then not all the corresponding operators are bounded.

\begin{proposition}\label{3}
Fix $N$, $\nu(>1), e$ satisfying (\ref{e}) and fix $n$. Suppose $\gamma_t(x)$ is a $C^\infty$ function of $(t,x)$ defined on a neighborhood of $(0,0)$ in $\mathbb{R}^N\times \mathbb{R}^n$ with $\gamma_0(x)\equiv x$, whose corresponding $W(t,x)$ can be written as a finite sum on some neighborhood of $(0,0)$:
$$
W(t,x) = \sum_{l=1}^Q t^{\alpha_l} X_l(x),
$$
where the $\alpha_l\in \mathbb{N}^N$ are such that $\text{deg}\,(\alpha_1) = \cdots = \text{deg}\,(\alpha_Q)\in [0,\infty)^\nu$ is a nonpure degree, the $X_l$ commute with each other, and there exists $1\leq l\leq Q$ such that $X_l(0) \neq 0$.

Then for every $p\in (1, \infty)$, every $a>0$, and every neighborhood $U$ of $0$ in $\mathbb{R}^n$, there exists $\psi \in C_c^\infty(U)$, $K\in \mathcal{K}(N,e,a,\nu)$, such that the operator
$$
Tf(x) = \psi(x)\int f(\gamma_t(x))K(t)\,dt
$$
is not bounded on $L^p(\mathbb{R}^n)$.
\end{proposition}

These three theorems lead to a contradiction. Hence the assumption that there exists a real analytic $\gamma_t(x)$ that satisfies \ref{SOMEP} but does not satisfy \ref{WTX} is false. Therefore \ref{SOMEP} implies \ref{WTX}.

\section{Preparation}\label{nianhuazhi}
We include in this section several definitions, lemmas, and theorems in preparation for the main proofs.

\subsection{Properties of control}
\begin{lemma}\label{4}
Suppose $(X_0, d_0), (X_1, d_1), \ldots, (X_q,d_q)$ are $C^\infty$ vector fields defined on some neighborhood of $0$ in $\mathbb{R}^n$, each paired with a $\nu$-parameter degree. And suppose $W(t,x)$ is a $C^\infty$ function mapping from a neighborhood of $(0,0)$ in $\mathbb{R}^N \times \mathbb{R}^n$ into $\mathbb{R}^n$ satisfying $W(0,x) \equiv 0$. Denote $(\mathbf{X},\mathbf{d}) :=\{(X_1,d_1), \ldots, (X_q,d_q)\}$.

\begin{enumerate}[label=(\roman*)]
\item\label{parta}
If $(\mathbf{X},\mathbf{d})$ strongly controls $(X_0, d_0)$, then $(\mathbf{X},\mathbf{d})$ controls $(X_0, d_0)$.

\item\label{partb}
If $(\mathbf{X},\mathbf{d})$ strongly controls $W$, then $(\mathbf{X},\mathbf{d})$ controls $W$.
\end{enumerate}
\end{lemma}

\begin{proof}
We only prove Part \ref{parta}; Part \ref{partb} is similar. By the assumption, on some neighborhood $\Omega$ of $0$,
$$
X_0 = \sum_{\substack{1\leq l\leq q\\ d_l \leq d_0}}c_lX_l,
$$
where the $c_l$ are $C^\infty$ on a neighborhood of the closure of $\Omega$, and the inequalities $d_l \leq d_0$ hold coordinate-wise. There exists $\xi_1>0$ and a neighborhood $U \subseteq \Omega$ of $0$, such that $(\mathbf{X},\mathbf{d})$ satisfies $\mathcal{C}(x, \Vec{\xi_1}, \Omega)$ for all $x\in U$. Thus for all $x \in U, \delta\in [0,1]^\nu$, 
$$
\delta^{d_0}X_0 = \sum_{\substack{1\leq l \leq q \\ d_l \leq d_0}} \big(\delta^{d_0-d_l} c_l\big) \delta^{d_l} X_l, \text{~on~} B_{(\mathbf{X},\mathbf{d})}(x,\xi_1\delta).
$$
For $d_l \leq d_0$ and $m \in \mathbb{N}$,
$$
\sup_{\substack{x \in U\\ \delta\in [0,1]^\nu}}\sum_{|\beta|\leq m}\Big\|(\delta \mathbf{X})^\beta \big(\delta^{d_0-d_l} c_l\big)\Big\|_{C(B_{(\mathbf{X},\mathbf{d})}(x,\xi_1\delta))} <\infty.
$$
\end{proof}

The transitivity of control holds in the following situation.

\begin{lemma}\label{5}
Suppose $(X_0,d_0), (X_1,d_1), \ldots, (X_q,d_q), (Y_1, D_1), \ldots, (Y_{Q}, D_{Q})$ are $C^\infty$ vector fields defined on a neighborhood $\Omega$ of $0$ in $\mathbb{R}^n$, each paired with a $\nu$-parameter degree. And suppose $W(t,x)$ is a $C^\infty$ function mapping from a neighborhood of $(0,0)$ in $\mathbb{R}^N \times \mathbb{R}^n$ into $\mathbb{R}^n$ satisfying $W(0,x) \equiv 0$. Denote $(\mathbf{X},\mathbf{d}) := \{(X_1,d_1), \ldots, (X_q,d_q)\}$, and $(\mathbf{Y},\mathbf{D}) := \{(Y_1, D_1), \ldots, (Y_{Q}, D_{Q})\}$. Suppose $(\mathbf{X},\mathbf{d})$ is a subset of $(\mathbf{Y},\mathbf{D})$.

\begin{enumerate}[label=(\roman*)]
\item\label{PARTA}
If $(\mathbf{Y},\mathbf{D})$ controls $(X_0,d_0)$, and $(\mathbf{X},\mathbf{d})$ controls every element of $(\mathbf{Y}, \mathbf{D})$, then $(\mathbf{X},\mathbf{d})$ controls $(X_0,d_0)$.

\item\label{PARTB}
If $(\mathbf{Y},\mathbf{D})$ controls $W$, and $(\mathbf{X},\mathbf{d})$ controls every element of $(\mathbf{Y}, \mathbf{D})$, then $(\mathbf{X},\mathbf{d})$ controls $W$.
\end{enumerate}
\end{lemma}

\begin{proof}
We only prove Part \ref{PARTA}; Part \ref{PARTB} is similar. By the assumptions, there exists a neighborhood $U\subseteq \Omega$ of $0$ and $\xi_1>0$, such that $(\mathbf{X},\mathbf{d})$ and $(\mathbf{Y}, \mathbf{D})$ both satisfy $\mathcal{C}(x, \Vec{\xi_1}, \Omega)$ for all $x\in U$, and such that for all $x\in U, \delta\in [0,1]^\nu$,
\begin{align*}
&\delta^{d_0} X_0 = \sum_{i=1}^{Q} \bar c_i^{x,\delta}\delta^{D_i} Y_i, \text{~on~} B_{(\mathbf{Y}, \mathbf{D})}(x, \xi_1\delta),\\
& \delta^{D_i} Y_i = \sum_{l=1}^q c_{i,l}^{x,\delta} \delta^{d_l} X_l, \text{~on~} B_{(\mathbf{X},\mathbf{d})}(x,\xi_1 \delta),
\end{align*}
and for all $m\in \mathbb{N}$,
\begin{align*}
&\sup_{\substack{x\in U\\ \delta \in [0,1]^\nu}} \sum_{|\beta|\leq m} \Big\| (\delta \mathbf{Y})^\beta \bar c_i^{x,\delta} \Big\|_{C(B_{(\mathbf{Y}, \mathbf{D})}(x, \xi_1 \delta)} <\infty,\\
&\sup_{\substack{x\in U\\ \delta \in [0,1]^\nu}} \sum_{|\beta| \leq m} \Big\| (\delta \mathbf{X})^\beta c_{i,l}^{x,\delta} \Big\|_{C(B_{(\mathbf{X},\mathbf{d})}(x,\xi_1 \delta))} <\infty,
\end{align*}
where $\delta \mathbf{X}= \{\delta^{d_l} X_l\}_{l=1}^q, \delta \mathbf{Y} = \{ \delta^{D_i} Y_i \}_{i=1}^Q$, and $\beta$ denotes an ordered multi-index. Since $(\mathbf{X},\mathbf{d})\subseteq (\mathbf{Y},\mathbf{D})$, for all $x\in U, \delta\in [0,1]^\nu$, $B_{(\mathbf{X},\mathbf{d})}(x,\xi_1\delta) \subseteq B_{(\mathbf{Y}, \mathbf{D})}(x,\xi_1\delta)$, and
$$
\delta^{d_0}X_0 = \sum_{i=1}^{Q} \bar c_{i}^{x,\delta} \sum_{l=1}^q c_{i,l}^{x,\delta}\delta^{d_l}X_l, \text{~on~} B_{(\mathbf{X},\mathbf{d})}(x,\xi_1\delta),
$$
and for all $m \in \mathbb{N}$,
\begin{align*}
&\quad \sup_{\substack{x\in U\\ \delta \in [0,1]^\nu}} \sum_{|\beta| \leq m}\Big\|(\delta \mathbf{X})^\beta \big( \bar c_{i}^{x,\delta} c_{i,l}^{x,\delta}\big)\Big\|_{C(B_{(\mathbf{X},\mathbf{d})}(x,\xi_1\delta))}\\
&\leq \sup_{\substack{x\in U\\ \delta \in [0,1]^\nu}} \sum_{|\beta| \leq m} \sum_{\beta' \subseteq \beta}\Big\|(\delta \mathbf{Y})^{\beta'} \bar c_{i}^{x,\delta}\Big\|_{C(B_{(\mathbf{Y},\mathbf{D})}(x,\xi_1\delta))} \Big\|(\delta \mathbf{X})^{\beta \backslash \beta'} c_{i,l}^{x,\delta}\Big\|_{C(B_{(\mathbf{X},\mathbf{d})}(x,\xi_1\delta))} <\infty,
\end{align*}
where $\beta'$ is a sublist of the ordered multi-index $\beta$, and $\beta\backslash \beta'$ is the complement sublist\footnote{For instance, if $\beta=(1,3,2)$, $\beta'$ could be $(1,2)$, then $\beta \backslash \beta' = (3)$.}. 
\end{proof}

The notion of control is preserved under Lie brackets, in the following sense.

\begin{lemma}\label{yicu}
Suppose $(X_0,d_0), (X_1, d_1), (X_2, d_2), \ldots, (X_q,d_q)$ are $C^\infty$ vector fields defined on a neighborhood $\Omega$ of $0$ in $\mathbb{R}^n$, each paired with a $\nu$-parameter degree. And suppose $(\mathbf{X},\mathbf{d}) :=\{(X_2, d_2), \ldots, (X_q,d_q)\}$ controls every element of $\{([X_i,X_j], d_i +d_j): 2\leq i,j \leq q \}$. If $(\mathbf{X},\mathbf{d})$ controls $(X_0,d_0)$ and $(X_1, d_1)$, then $(\mathbf{X},\mathbf{d})$ also controls $([X_0, X_1], d_0 + d_1)$.
\end{lemma}

\begin{proof}
By the assumptions, there exists $\xi_1>0$ and a neighborhood $U \subseteq \Omega$ of $0$, such that $(\mathbf{X},\mathbf{d})$ satisfies $\mathcal{C}(x, \Vec{\xi_1}, \Omega)$ for all $x\in U$, and such that for all $x\in U, \delta\in [0,1]^\nu$, on $B_{(\mathbf{X},\mathbf{d})}(x,\xi_1\delta)$,
\begin{align*}
&\big[ \delta^{d_l} X_l, \delta^{d_{l'}} X_{l'} \big] = \sum_{i=2}^q c_{l,l',i}^{x,\delta} \delta^{d_i} X_i \text{~~for~} 2 \leq l, l' \leq q, \quad \delta^{d_0} X_0 = \sum_{l=2}^q c_{0,l}^{x,\delta} \delta^{d_l} X_l, \quad \delta^{d_1} X_1 = \sum_{l=2}^q c_{1,l}^{x,\delta} \delta^{d_l} X_l,
\end{align*}
and for all $m\in \mathbb{N}$,
\begin{align*}
&\sup_{\substack{x\in U\\ \delta \in [0,1]^\nu}} \sum_{|\beta|\leq m} \big\| (\delta \mathbf{X})^\beta c_{l, l', i}^{x,\delta} \big\|_{C(B_{(\mathbf{X},\mathbf{d})}(x,\xi_1\delta))} <\infty,\\
&\sup_{\substack{x\in U\\ \delta \in [0,1]^\nu}} \sum_{|\beta|\leq m} \big\| (\delta \mathbf{X})^\beta c_{0,l}^{x,\delta} \big\|_{C(B_{(\mathbf{X},\mathbf{d})}(x,\xi_1\delta))} <\infty,\\
&\sup_{\substack{x\in U\\ \delta \in [0,1]^\nu}} \sum_{|\beta|\leq m} \big\| (\delta \mathbf{X})^\beta c_{1,l}^{x,\delta} \big\|_{C(B_{(\mathbf{X},\mathbf{d})}(x,\xi_1\delta))} <\infty,
\end{align*}
where $\delta \mathbf{X} = \{\delta^{d_l} X_l\}_{l=2}^q$. Therefore for all $x\in U, \delta\in [0,1]^\nu$, on $B_{(\mathbf{X},\mathbf{d})}(x,\xi_1\delta)$,
\begin{align*}
\delta^{d_0+d_1} [X_0, X_1] &= \Big[\sum_{l=2}^q c_{0,l}^{x,\delta} \delta^{d_l} X_l, \sum_{l=2}^q c_{1,l}^{x,\delta} \delta^{d_l} X_l\Big]\\
&= \sum_{l,l'=2}^q \Big( c_{0,l}^{x,\delta} \delta^{d_l}X_l \big(c_{1,l'}^{x,\delta}\big) \delta^{d_{l'}} X_{l'} - c_{1,l'}^{x,\delta} \delta^{d_{l'}} X_{l'} \big( c_{0,l}^{x,\delta} \big) \delta^{d_l} X_l + c_{0,l}^{x,\delta} c_{1,l'}^{x,\delta} \big[ \delta^{d_l}X_l, \delta^{d_{l'}} X_{l'} \big] \Big)\\
& = \sum_{l,l'=2}^q \Big( c_{0,l}^{x,\delta} \delta^{d_l}X_l \big(c_{1,l'}^{x,\delta}\big) \delta^{d_{l'}} X_{l'} - c_{1,l'}^{x,\delta} \delta^{d_{l'}} X_{l'} \big( c_{0,l}^{x,\delta} \big) \delta^{d_l} X_l + c_{0,l}^{x,\delta} c_{1,l'}^{x,\delta}  \sum_{i=2}^q c_{l,l',i}^{x,\delta} \delta^{d_i} X_i \Big)\\
&= \sum_{l=2}^q \Big( \sum_{l'=2}^q \Big( c_{0,l'}^{x,\delta} \delta^{d_{l'}} X_{l'}(c_{1,l}^{x,\delta}) - c_{1,l'}^{x,\delta} \delta^{d_{l'}} X_{l'} (c_{0,l}^{x,\delta}) + \sum_{i=2}^q c_{0,i}^{x,\delta} c_{1,l'}^{x,\delta} c_{i,l',l}^{x,\delta} \Big)  \Big) \delta^{d_l} X_l,
\end{align*}
and for all $m\in \mathbb{N}$,
\begin{align*}
&\sup_{\substack{x\in U\\ \delta \in [0,1]^\nu}} \sum_{|\beta|\leq m} \Big\| (\delta \mathbf{X})^\beta \Big( \sum_{l'=2}^q \Big( c_{0,l'}^{x,\delta} \delta^{d_{l'}} X_{l'}(c_{1,l}^{x,\delta}) - c_{1,l'}^{x,\delta} \delta^{d_{l'}} X_{l'} (c_{0,l}^{x,\delta}) + \sum_{i=2}^q c_{0,i}^{x,\delta} c_{1,l'}^{x,\delta} c_{i,l',l}^{x,\delta} \Big)  \Big) \Big\|_{C(B_{(\mathbf{X},\mathbf{d})}(x,\xi_1\delta))} <\infty.
\end{align*}
\end{proof}

The following definition includes several notations that appear in Theorem \ref{16} and in the main proofs.

\begin{definition}\label{indoor}
\begin{itemize}
\item 
Denote $m\vee n :=\max\{m,n\}, m\wedge n: = \min\{m,n\}$. For an $m\times n$ matrix $A$, and for $n_0\leq m \wedge n$, let $\det_{n_0\times n_0} A$ be the vector whose entries are the determinants of $n_0\times n_0$ submatrices of $A$. The order of the entries in this vector does not matter.

\item
Let $|\cdot|_1, |\cdot|_\infty$, and $|\cdot|$ denote the $l^1$, $l^\infty$, and $l^2$ norms of a vector, respectively.

\item
For $n_0 \leq n$, let $\mathcal{I}(n_0,n)$ be the set of all lists of integers $(i_1, \ldots, i_{n_0})$ such that $1 \leq i_1 < i_2< \cdots <i_{n_0} \leq n$.

\item 
For $\xi>0$, $\Vec{\xi}$ denotes the vector $(\xi, \ldots, \xi) \in (0,\infty)^\nu$.

\item
For $d_1, d_2 \in [0,\infty)^\nu$, $d_1 \leq d_2$ means the inequality holds coordinate-wise.

\item
For $(\mathbf{X},\mathbf{d})=\{(X_1, d_1), \ldots, (X_q, d_q)\}$ a set of $C^\infty$ vector fields defined on a neighborhood of $0$ in $\mathbb{R}^n$, each paired with a $\nu$-parameter degree, and for some subset $J \subseteq \{1, \ldots, q\}$, $\mathbf{X}$ denotes the matrix with column vectors $\{X_l\}_{l=1}^q$, $\mathbf{X}_J$ denotes the matrix with column vectors $\{X_l\}_{l\in J}$, $\delta \mathbf{X}$ denotes the matrix with column vectors $\{ \delta^{d_l} X_l\}_{l=1}^q$, and $(\delta \mathbf{X})_J$ denotes the matrix with column vectors $\{ \delta^{d_l} X_l\}_{l\in J}$.

\item
Let $|J|$ denote the number of elements in $J$. For $u\in \mathbb{R}^{|J|}$, let $u\cdot \mathbf{X}_J= \sum_{l\in J} u_l X_l$, and $u\cdot (\delta \mathbf{X})_J= \sum_{l\in J} u_l \delta^{d_l} X_l$, which are still $C^\infty$ vector fields on a neighborhood of $0$ in $\mathbb{R}^n$.

\item
An ordered multi-index $\alpha$ is a list of numbers from $\{1, \ldots, q\}$, and $|\alpha|$ denotes the length of this list. $(\delta \mathbf{X})^\alpha$ means we juxtapose vector fields from $\{\delta^{d_l} X_l\}_{l=1}^q$ indexed with list $\alpha$. For instance, if $\alpha = (1,2,2)$, then $|\alpha| =3$, and $(\delta \mathbf{X})^\alpha = (\delta^{d_1}X_1)(\delta^{d_2}X_2)(\delta^{d_2}X_2)$.

\item 
Let $(\delta \mathbf{X})_{I\times J}$ be the submatrix of $\delta \mathbf{X}$ formed with rows indexed by $I$ and columns indexed by $J$. Let $(\delta^{d_l} X_l)_I$ be the truncated vector obtained by keeping entries indexed by $I$ from $\delta^{d_l}X_l$.

\item
Note the above matrices $\mathbf{X}(x), \mathbf{X}_J(x), \delta \mathbf{X}(x)$, $(\delta \mathbf{X})_J(x)$, and $(\delta \mathbf{X})_{I\times J}(x)$ depend smoothly on $x$ in a neighborhood of $0$. 
\end{itemize}
\end{definition}

The local diffeomorphism $\Phi$ in the following theorem is the main tool of this paper (it generalizes dilations $\Phi_{2^{2L}}$ in Section \ref{real} and $\Phi_{\tau_1, \tau_2}$ in Section \ref{sntz}).

Given $(\mathbf{X},\mathbf{d}):=\{(X_1,d_1), \ldots,(X_q,d_q)\}$ a set of $C^\infty$ vector fields defined on a neighborhood $\Omega$ of $0$ in $\mathbb{R}^n$, each paired with a $\nu$-parameter degree. Fix $x_0\in \Omega, \xi_1>0$ such that $(\mathbf{X},\mathbf{d})$ satisfies $\mathcal{C}(x_0, \Vec{\xi_1}, \Omega)$. Suppose for $1\leq l,k, i\leq q$, there exist functions $c_{l,k}^i$ defined on $B_{(\mathbf{X},\mathbf{d})}(x_0, \Vec{\xi_1})$, such that
$$
[X_l, X_k]=\sum_{i=1}^q c_{l,k}^{i} X_i, \text{~~on~} B_{(\mathbf{X},\mathbf{d})}(x_0,  \Vec{\xi_1}),
$$
and\footnote{We write $\|f\|_{C^m(U)} := \sum_{|\beta|\leq m} \sup_{u\in U} |\partial_u^\beta f(u)|$, and if we say the norm is finite, we mean (in addition) that the partial derivatives up to order $m$ of $f$ exist and are continuous on $U$. If $f$ is replaced by a vector field $Y = \sum_i b_i(u) \partial_{u_i}$, then we let $\|Y\|_{C^m(U)} := \sum_i \| b_i \|_{C^m(U)}$.} for all $m\in \mathbb{N}$,
$$
\big \|X_l \big\|_{C^m(B_{(\mathbf{X},\mathbf{d})}(x_0,  \Vec{\xi_1}))} <\infty, \quad \sum_{|\beta|\leq m}\big\|\mathbf{X}^\beta c_{l,k}^{i} \big\|_{C(B_{(\mathbf{X},\mathbf{d})}(x_0,  \Vec{\xi_1}))}<\infty.
$$
Denote $n_0 := \textup{rank} \,\mathbf{X}(x_0)$. Fix a $\zeta \in (0,1]$. Choose a $J \in \mathcal{I}(n_0,q)$ such that
$$
\Big|\det_{n_0\times n_0}\mathbf{X}_J(x_0)\Big|_\infty \geq \zeta \Big|\det_{n_0\times n_0}\mathbf{X}(x_0) \Big|_\infty.
$$

\begin{definition}\label{admissible}
For $m\geq 2$, we say $C$ is an $m$-admissible constant\footnote{In different sections, $m$-admissible constants may be allowed to depend on different parameters. We will specify in each section what an $m$-admissible constant is allowed to depend on.} if $C$ is chosen to depend only on fixed upper and lower bounds for $\big\{\big|d_l \big|_1: 1\leq l \leq q\big\}, m$, $\zeta$, fixed lower bound for $\xi_1$, fixed upper bound for $n,q,\nu$, fixed upper bound for the quantities:
$$
\big\|X_l\big\|_{C^m(B_{(\mathbf{X},\mathbf{d})}(x_0,  \Vec{\xi_1}))}, \quad \sum_{|\beta|\leq m}\big\|\mathbf{X}^\beta c_{i,j}^{k}\big\|_{C(B_{(\mathbf{X},\mathbf{d})}(x_0,  \Vec{\xi_1}))}.
$$
We write $A\lesssim_m B$ to mean $A\leq CB$ with $C\geq 0$ an m-admissible constant, write $A\approx_m B$ to mean $A \lesssim_m B$ and $B\lesssim_m A$, write $A\lesssim B$ to mean $A\lesssim_2 B$, and write $A\approx B$ to mean $A\approx_2 B$.

When we wish to make the choice of $(\mathbf{X}, \mathbf{d})$ and $x_0$ clear, we say, $C$ is an $m$-admissible constant with respect to $(\mathbf{X}, \mathbf{d})$ and $x_0$.
\end{definition}

\begin{theorem}\label{16}(Theorem 4.1 in \cite{CARNOT})
There exist $2$-admissible constants $\eta_1>0, 0<\xi_2<\xi_1$, such that we can define a diffeomorphism $\Phi: B^{n_0}(\eta_1)\to \Phi(B^{n_0}(\eta_1))\subseteq \mathbb{R}^n$ as
$$
\Phi(u) = e^{u\cdot \mathbf{X}_{J}}x_0,
$$
and
$$
B_{(\mathbf{X},\mathbf{d})}(x_0, \Vec{\xi_2}) \subseteq \Phi(B^{n_0}(\eta_1)) \subseteq B_{(\mathbf{X},\mathbf{d})}(x_0,  \Vec{\xi_1}).
$$
Note $J=J(x_0)$ is fixed. For every $0<\xi_3\leq \xi_1$, and every $0<\eta_3\leq \eta_1$, there exists $\xi_4>0$ only depending on $\eta_3$ and the quantities that a 2-admissible constant is allowed to depend on, and there exists $\eta_4>0$ only depending on $\xi_3$ and the quantities that a 2-admissible constant is allowed to depend on, such that
\begin{align*}
\Phi(B^{n_0}(\eta_4)) \subseteq B_{(\mathbf{X},\mathbf{d})}(x_0, \Vec{\xi_3}), \quad B_{(\mathbf{X},\mathbf{d})}(x_0, \Vec{\xi_4})\subseteq  \Phi(B^{n_0}(\eta_3)).
\end{align*}
Let $Y_l$ be the pullback of $X_l$ under the map $\Phi$ to $B^{n_0}(\eta_1)$, and let $\mathbf{Y}$ be the matrix with column vectors $\{Y_l\}_{1\leq l \leq q}$. Then on $B^{n_0}(\eta_1)$,
$$
\Big|\det_{n_0\times n_0}\mathbf{Y} \Big| \approx \Big|\det \mathbf{Y}_{J} \Big| \approx 1, \quad \Big| \det_{n_0\times n_0} D\Phi\Big| \approx \Big| \det_{n_0\times n_0} \mathbf{X}(x_0)\Big|.
$$
For every $2$-admissible constant $0<\eta \leq \eta_1$, every $m\in \mathbb{N}$, and every $f\in C^m(B^{n_0}(\eta))$,
$$
\big\|Y_{l}\big\|_{C^m(B^{n_0}(\eta_1))} \lesssim_{m\vee 2} 1, \quad 
\big\|f\big\|_{C^m(B^{n_0}(\eta))} \approx_{(m-1)\vee 2} \sum_{|\alpha|\leq m} \big\|\mathbf{Y}^\alpha f \big\|_{C(B^{n_0}(\eta))}.
$$
\end{theorem}

The containment relation between multi-parameter Carnot-Carath\'eodory balls holds under the following assumptions.

\begin{lemma}\label{CONTAIN}
Suppose $(X_1,d_1), \ldots, (X_q,d_q), (Z_1, D_1), \ldots, (Z_{Q}, D_{Q})$ are $C^\infty$ vector fields defined on a neighborhood $\Omega$ of $0$ in $\mathbb{R}^n$, each paired with a $\nu$-parameter degree. And suppose there exists $\xi_1>0$ and a neighborhood $U\subseteq \Omega$ of $0$, such that $(\mathbf{X},\mathbf{d}):= \{(X_1, d_1), \ldots, (X_q, d_q)\}$ satisfies $\mathcal{C}(x, \Vec{\xi_1}, \Omega)$ for every $x\in U$.
If $(\mathbf{X},\mathbf{d})$ controls every element of $\{ ([X_i, X_j], d_i +d_j) : 1\leq i, j \leq q\}$ and every element of $(\mathbf{Z},\mathbf{D}) := \{(Z_1, D_1), \ldots, (Z_Q, D_Q)\}$, on $U$ with parameter $\xi_1$, then there exists $c>0$ such that for every $x\in U, \delta \in [0,1]^\nu$,
$$
B_{(\mathbf{Z},\mathbf{D})}(x, c\xi_1 \delta) \subseteq B_{(\mathbf{X},\mathbf{d})}(x,\xi_1\delta),
$$
and $(\mathbf{Z},\mathbf{D})$ satisfies $\mathcal{C}(x, c\Vec{\xi_1}, \Omega)$ for all $x\in U$.
\end{lemma}

\begin{proof} 
Since $(\mathbf{X},\mathbf{d})$ controls every element of $\{ ([X_i, X_j], d_i +d_j) : 1\leq i, j \leq q\}$ on $U$ with parameter $\xi_1$, $\{(\delta^{d_l}X_l, d_l)\}_{l=1}^q$ satisfies the conditions of Theorem \ref{16} on $B_{(\mathbf{X},\mathbf{d})} (x, \xi_1 \delta) = B_{(\delta \mathbf{X},\mathbf{d})} (x, \Vec{\xi_1})$ uniformly in $x\in U, \delta \in [0,1]^\nu$. For $x\in U, \delta \in [0,1]^\nu$, let $n_0(x,\delta) := \text{rank}\,(\delta \mathbf{X})(x)$, and choose $J(x,\delta) \in \mathcal{I}(n_0(x,\delta), q)$ such that
$$
\Big| \det_{n_0(x,\delta)\times n_0(x,\delta)} (\delta \mathbf{X})_{J(x,\delta)}(x) \Big|_\infty = \Big| \det_{n_0(x,\delta)\times n_0(x,\delta)} (\delta \mathbf{X})(x) \Big|_\infty.
$$
Then there exists $\eta_1>0$, $0< \xi_2 < \xi_1$, such that for every $x\in U, \delta\in [0,1]^\nu$, we can define a diffeomorphism 
$$
\Phi_{x, \delta}: B^{n_0(x, \delta)}(\eta_1) \to \Phi_{x, \delta}(B^{n_0(x, \delta)}(\eta_1))\subseteq \mathbb{R}^n, \quad u\mapsto e^{u\cdot (\delta \mathbf{X})_{J(x, \delta)}}x,
$$
satisfying
$$
B_{(\mathbf{X},\mathbf{d})}(x, \xi_2\delta) \subseteq \Phi_{x,\delta}(B^{n_0(x,\delta)}(\eta_1)) \subseteq B_{(\mathbf{X},\mathbf{d})}(x, \xi_1\delta).
$$
Since $(\mathbf{X},\mathbf{d})$ controls every element of $(\mathbf{Z},\mathbf{D})$ on $U$ with parameter $\xi_1$, for $1\leq l\leq Q$, $x \in U$, and $\delta \in [0,1]^\nu$,
\begin{equation}\label{111}
\delta^{D_l} Z_l = \sum_{i=1}^q c^{x, \delta}_{l,i} \delta^{d_i} X_i, \text{~on~} B_{(\mathbf{X},\mathbf{d})}(x, \xi_1 \delta),
\end{equation}
and for all $m\in \mathbb{N}$,
$$
\sup_{\substack{x\in U\\ \delta \in [0,1]^\nu}} \sum_{|\beta| \leq m} \Big\| (\delta \mathbf{X})^\beta c_{l,i}^{x, \delta} \Big\|_{C(B_{(\mathbf{X},\mathbf{d})}(x, \xi_1 \delta))} <\infty.
$$
Let $Y_i^{x, \delta}$ be the pullback of $\delta^{d_i} X_i$ under the map $\Phi_{x, \delta}$ to $B^{n_0(x,\delta)}(\eta_1)$. Pull back (\ref{111}) and obtain
$$
\Phi_{x, \delta}^*(\delta^{D_l} Z_l) = \sum_{i=1}^q ~c_{l,i}^{x, \delta} \circ \Phi_{x, \delta} ~ Y^{x, \delta}_i, \quad \text{on } B^{n_0(x, \delta)}(\eta_1).
$$
By Theorem \ref{16}, the vector fields $Y_i^{x, \delta}$ and the coefficients $c_{l,i}^{x, \delta} \circ \Phi_{x, \delta}$ are bounded in $C^m(B^{n_0(x,\delta)}(\eta_1))$ for all $m$ uniformly in $x,\delta$. Thus $\Phi_{x,\delta}^*(\delta^{D_l} Z_l)$ are bounded in $C^m(B^{n_0(x,\delta)}(\eta_1))$ for all $m$ uniformly in $x, \delta$.

By Picard-Lindel\"of theorem, for every $\{a_l(t)\}_{l=1}^Q \in \big(L^\infty([0,1])\big)^Q$ satisfying $\Big\|\sum_{l=1}^Q |a_l(t)|^2 \Big\|_{L^\infty[0,1]} \leq 1$, and for every $x\in U, \delta \in [0,1]^\nu$, there exists locally an absolutely continuous curve $\gamma(t)$ in $B^{n_0(x, \delta)}(\eta_1)$ being the solution to
\begin{equation}\label{222}
\gamma(0) = 0, \quad \gamma'(t) = \sum_{l=1}^Q a_l(t) \xi_1^{|D_l|_1} \Phi_{x,\delta}^*(\delta^{D_l} Z_l) (\gamma(t)).
\end{equation}
The above right-hand side is Liptchitz in $\gamma(t)$ uniformly in $x, \delta$, and $\{a_l(t)\}_{l=1}^Q$. Thus there exists a constant $0< c_0 \leq 1$ independent of $x,\delta, \{a_l(t)\}$, such that for $t\in [0, c_0]$, the solution $\gamma(t)$ exists and $\gamma(t) \in B^{n_0(x,\delta)}(\eta_1)$. Let $\tilde \gamma(t) := \Phi_{x, \delta} \circ \gamma(t) \in B_{(\mathbf{X},\mathbf{d})}(x, \xi_1\delta)$, for $t\in [0, c_0]$. Push forward (\ref{222}) and have
$$
\tilde \gamma(0) = x, \quad \tilde \gamma'(t) = \sum_{l=1}^Q a_l(t) \xi_1^{|D_l|_1} \delta^{D_l} Z_l(\tilde \gamma(t)).
$$
By arbitrariness of $\{a_l(t)\}$, for $c= c_0^{\frac{1}{\min_l |D_l|_1}}$,
$$
B_{(\mathbf{Z},\mathbf{D})}(x, c\xi_1\delta) \subseteq B_{(\mathbf{X},\mathbf{d})}(x, \xi_1\delta),
$$
and $(\mathbf{Z},\mathbf{D})$ satisfies $\mathcal{C}(x, c\Vec{\xi_1}, \Omega)$ for all $x\in U$.
\end{proof}

\subsection{Equivalence of \texorpdfstring{$\mathcal{P}, \mathcal{N}$}{TEXT} and \texorpdfstring{$\widehat{\mathcal{P}}, \widehat{\mathcal{N}}$}{TEXT}}
\begin{theorem}\label{ljj}
For a $C^\infty$ $\gamma_t(x)$ defined on a neighborhood of $(t,x)= (0,0)$ in $\mathbb{R}^N \times \mathbb{R}^n$ satisfying $\gamma_0(x) \equiv x$,

\begin{enumerate}[label=(\roman*)]
\item\label{partone}
every element of $\mathcal{L}(\mathcal{P} \cup \mathcal{N})$ is a linear combination of elements in $\mathcal{L}(\widehat{\mathcal{P}} \cup \widehat{\mathcal{N}})$ with the same degree with constant coefficients, and every element of $\mathcal{L}(\mathcal{P} )$ is a linear combination of elements in $\mathcal{L}(\widehat{\mathcal{P}} )$ with the same degree with constant coefficients.

\item\label{parttwo}
every element of $\mathcal{L}(\widehat{\mathcal{P}} \cup \widehat{\mathcal{N}})$ is a linear combination of elements in $\mathcal{L}(\mathcal{P} \cup \mathcal{N})$ with the same degree with constant coefficients, and every element of $\mathcal{L}(\widehat{\mathcal{P}} )$ is a linear combination of elements in $\mathcal{L}(\mathcal{P} )$ with the same degree with constant coefficients.
\end{enumerate}
\end{theorem}

\begin{proof}
By the Baker-Campbell-Hausdorff formula, 
\begin{align*}
W(t,x) = \frac{d}{d\epsilon}\Big|_{\epsilon =0} \gamma_{(1+\epsilon)t} \circ \gamma_t^{-1} (x) \sim \frac{d}{d\epsilon}\Big|_{\epsilon =0} e^{\sum_{|\alpha|>0} (1+\epsilon)^{|\alpha|} t^\alpha \widehat X_\alpha} e^{-\sum_{|\alpha|>0}t^\alpha \widehat X_\alpha}x \sim \sum_{|\alpha|>0}  t^\alpha (|\alpha| \widehat X_\alpha(x) + V_\alpha(x)) ,
\end{align*}
where $V_\alpha$ is a linear combination of elements in $\mathcal{L}(\widehat{\mathcal{P}} \cup \widehat{\mathcal{N}}) \backslash (\widehat{\mathcal{P}} \cup \widehat{\mathcal{N}})$ with degree equal to $\text{deg}\,(\alpha)$ and with constant coefficients. Thus every element of $\mathcal{P} \cup \mathcal{N}$ (respectively $\mathcal{P}$) is a linear combination of elements in $\mathcal{L}(\widehat{\mathcal{P}} \cup \widehat{\mathcal{N}})$ (respectively $\mathcal{L}(\widehat{\mathcal{P}})$) with the same degree with constant coefficients. Taking Lie brackets, Part \ref{partone} follows.

When $|\alpha|=1$, $V_\alpha =0$, and thus $\widehat X_\alpha =  X_\alpha$. We induct on $|\alpha|$. Suppose every $\widehat X_\alpha \in \widehat{\mathcal{P}} \cup \widehat{\mathcal{N}}$ with $|\alpha|\leq m$, is a linear combination of elements in $\mathcal{L}(\mathcal{P}\cup \mathcal{N})$ with degree equal to $\text{deg}\,(\alpha)$ and with constant coefficients. Then for $|\alpha| = m+1$, $V_\alpha \in \mathcal{L}(\widehat{\mathcal{P}} \cup \widehat{\mathcal{N}}) \backslash (\widehat{\mathcal{P}} \cup \widehat{\mathcal{N}})$ is a linear combination of elements in $\mathcal{L}(\mathcal{P}\cup \mathcal{N})$ with degree equal to $\text{deg}\,(\alpha)$ and with constant coefficients, and hence
$$
\widehat X_\alpha = \frac{1}{|\alpha|} ( X_\alpha - V_\alpha)
$$
is a linear combination of elements in $\mathcal{L}(\mathcal{P}\cup \mathcal{N})$ with degree equal to $\text{deg}\,(\alpha)$ and with constant coefficients. And if, in addition, $\widehat X_\alpha \in \widehat{\mathcal{P}}$, then $\text{deg}\,(\alpha)$ is pure,  and thus $\widehat X_\alpha$ is a linear combination of elements in $\mathcal{L}(\mathcal{P})$ with degree equal to $\text{deg}\,(\alpha)$ and with constant coefficients. By induction, every element of $\widehat{\mathcal{P}} \cup \widehat{\mathcal{N}}$ (respectively $\widehat{\mathcal{P}}$) is a linear combination of elements in $\mathcal{L}(\mathcal{P} \cup \mathcal{N})$ (respectively $\mathcal{L}(\mathcal{P})$) with the same degree with constant coefficients. Taking Lie brackets, Part \ref{parttwo} follows.
\end{proof}

We now can show the two conditions stated in Theorem \ref{lieri} are equivalent.

\begin{corollary}\label{EQUIV}
Let $\gamma_t(x)$ be a $C^\infty$ function of $(t,x)$ defined on a neighborhood of $(0,0)$ in $\mathbb{R}^N \times \mathbb{R}^n$ with $\gamma_0(x) \equiv x$. The following are equivalent.
\begin{enumerate}[label=(\alph*)]
\item\label{NOHAT}
$\mathcal{L}(\mathcal{P})$ controls $\mathcal{N}$.

\item\label{HAT}
$\mathcal{L}(\widehat{\mathcal{P}})$ controls $\widehat{\mathcal{N}}$.
\end{enumerate}
\end{corollary}

\begin{proof}
We only prove \ref{NOHAT} $\Rightarrow$ \ref{HAT}; the other direction is identical. Assume \ref{NOHAT} holds. Then by Definition \ref{defn} there exists a finite subset $\mathcal{F} \subseteq \mathcal{L}(\mathcal{P})$ that controls every element of $\mathcal{L}(\mathcal{P})$ and $\mathcal{N}$. By Lemma \ref{yicu}, $\mathcal{F}$ controls every element of $\mathcal{L}(\mathcal{P} \cup \mathcal{N})$. By Theorem \ref{ljj}, there exists a finite subset $\widehat{\mathcal{F}} \subseteq \mathcal{L}(\widehat{\mathcal{P}})$ such that each element of $\mathcal{F}$ is a linear combination of elements in $\widehat{\mathcal{F}}$ with the same degree with constant coefficients. Now it suffices to show $\widehat{\mathcal{F}}$ controls every element of $\mathcal{L}(\widehat{\mathcal{P}} \cup \widehat{\mathcal{N}})$. By Theorem \ref{ljj}, it suffices to show $\widehat{\mathcal{F}}$ controls every element of $\mathcal{L}(\mathcal{P} \cup \mathcal{N})$. Since $\mathcal{F}$ controls every element of $\mathcal{L}(\mathcal{P})$, by Theorem \ref{ljj}, $\mathcal{F}$ controls every element of $\mathcal{L}(\widehat{\mathcal{P}})$, and in particular, $\mathcal{F}$ controls every element of $\widehat{\mathcal{F}}$.

There exists $\xi_1>0$ and neighborhoods $U \subseteq \Omega$ of $0$ in $\mathbb{R}^n$, such that every element of $\mathcal{F}$ and $\widehat{\mathcal{F}}$ is defined on $\Omega$, $\mathcal{F}$ satisfies $\mathcal{C}(x,\Vec{\xi_1}, \Omega)$ for $x\in U$, and $\mathcal{F}$ controls every element of $\big\{\big([X_i,X_j], d_i+d_j\big): (X_i,d_i), (X_j,d_j) \in \mathcal{F}\big\}$
and every element of $\widehat{\mathcal{F}}$, on $U$ with parameter $\xi_1$.

By Lemma \ref{CONTAIN}, there exists $c>0$, such that for all $x\in U, \delta \in [0,1]^\nu$,
\begin{equation}\label{CONT}
B_{\widehat{\mathcal{F}}}(x,c\xi_1\delta) \subseteq B_{\mathcal{F}}(x,\xi_1 \delta),
\end{equation}
and $\widehat{\mathcal{F}}$ satisfies $\mathcal{C}(x, c\Vec{\xi_1}, \Omega)$ for all $x\in U$.
Fix an arbitrary $(X_0, d_0) \in \mathcal{L}(\mathcal{P}\cup \mathcal{N})$. $\mathcal{F}$ controls $(X_0, d_0)$. Thus there exists a neighborhood $U_0 \subseteq U$ of $0$ and some $0< \xi_0 \leq \xi_1$, such that for all $x\in U_0, \delta\in [0,1]^\nu$, 
$$
\delta^{d_0} X_0 = \sum_{(X_l, d_l)\in \mathcal{F}} c^{x,\delta}_l \delta^{d_l} X_l, \text{~on~} B_{\mathcal{F}} (x, \xi_0 \delta),
$$
and for all $m\in \mathbb{N}$,
$$
\sup_{\substack{x\in U_0\\ \delta\in [0,1]^\nu}} \sum_{|\beta|\leq m} \Big\| (\delta \mathbf{X})^\beta c^{x,\delta}_l \Big\|_{C(B_{\mathcal{F}} (x, \xi_0 \delta))} <\infty,
$$
where $\delta \mathbf{X}$ denotes $\{\delta^{d_l} X_l\}_{(X_l,d_l)\in \mathcal{F}}$, and $\beta$ denotes an ordered multi-index. Therefore for $x\in U_0, \delta\in [0,1]^\nu$,
$$
\delta^{d_0} X_0 = \sum_{(\widehat X_l, D_l)\in \widehat{\mathcal{F}}} \hat c^{x,\delta}_l \delta^{D_l} \widehat X_l, \text{~on~} B_{\widehat{\mathcal{F}}} (x, c\xi_0 \delta),
$$
and for all $m\in \mathbb{N}$,
$$
\sup_{\substack{x\in U_0\\ \delta\in [0,1]^\nu}} \sum_{|\beta|\leq m} \Big\| (\delta \mathbf{\widehat{X}})^\beta \hat c^{x,\delta}_l \Big\|_{C(B_{\widehat{\mathcal{F}}} (x, c\xi_0 \delta))} <\infty,
$$
where $\delta \mathbf{\widehat{X}}$ denotes $\{\delta^{D_l} \widehat X_l\}_{(\widehat X_l,D_l)\in \widehat{\mathcal{F}}}$. This is because every element of $\mathcal{F}$ is a finite linear combination of elements of $\widehat{\mathcal{F}}$ with the same degree with constant coefficients, and $\mathcal{F}$ controls every element of $\widehat{\mathcal{F}}$. Hence $\widehat{\mathcal{F}}$ controls $(X_0,d_0)$.
\end{proof}

\subsection{Properties of multi-parameter singular kernels}
\begin{lemma}\label{ff} (Proposition 16.3 in \cite{L2})
Let $\delta_0(t)$ be the Dirac function on $\mathbb{R}^N$ with point mass at $0$. For every $\nu>0$, $a>0$, and for every $e$ satisfying (\ref{e}), $\delta_0(t) \in \mathcal{K}(N,e,a,\nu)$.
\end{lemma}

For $K(t) \in \mathcal{K}(N,e,a,\nu)$, the way to decompose $K(t)=\sum \varsigma_k^{(2^k)}(t)$ is not unique, i.e. the map $\{\varsigma_k(t)\}_k \mapsto K(t):=\sum \varsigma_k^{(2^k)}(t)$ is not one-to-one. Unfortunately, $\mathcal{K}(N,e,a,\nu)$ is not a Fr\'echet space. However, in the next definition, we introduce a slight modification of this space, which is a Fr\'echet space. This will allow us to apply the closed graph theorem to this new space.

\begin{definition}\label{semi}
For every sequence $\{\varsigma_k(t)\}_{k\in \mathbb{N}^\nu}$ in $C_c^\infty(\overline{B^N(a)})$, and for every $m\in \mathbb{N}$, let $\big\| \{\varsigma_k\}_{k\in \mathbb{N}^\nu} \big\|_m = \sup_{k\in \mathbb{N}^\nu} \| \varsigma_k \|_{C^m}$. Define the space $\overline{\mathcal{K}}(N,e,a,\nu)$ as the set
\begin{align*}
&\Big\{ \{\varsigma_k\}_{k\in \mathbb{N}^\nu} \in C_c^\infty(\overline{B^N(a)}) : \big\| \{\varsigma_k\} \big\|_m < \infty \text{~for every~} m\in \mathbb{N},\\
&\qquad \int \varsigma_k(t)\,dt^\mu \equiv 0 \text{~for every~} k\in \mathbb{N}^\nu, 1\leq \mu \leq \nu \text{~satisfying~} k_\mu \neq 0\Big\}
\end{align*}
equipped with the family of semi-norms $\{\|\cdot \|_m\}_{m\in \mathbb{N}}$.
\end{definition}

\begin{lemma}\label{kk}
$\overline{\mathcal{K}}(N,e,a,\nu)$ is a Fr\'echet space.
\end{lemma}

\begin{proof}
If $\big\| \{\varsigma_k\}_{k} \big\|_m =0$ for all $m$, then $\{\varsigma_k\}_{k}=0$. It suffices to show the space is complete. Consider a Cauchy sequence $\{\{\varsigma^j_k\}_k\}_j \subseteq \overline{\mathcal{K}}(N,e,a,\nu)$. For each $m$, $k$,
\begin{align*}
\lim_{j,l\to \infty} \big\| \varsigma^j_k - \varsigma^l_k \big\|_{C^m} \leq \lim_{j,l\to \infty} \sup_k \big\| \varsigma^j_k - \varsigma^l_k \big\|_{C^m} = \lim_{j,l\to \infty} \big\| \{\varsigma^j_k-\varsigma^l_k\}_k \big\|_m =0.
\end{align*}
Fix $k$, $\varsigma_k^j$ converges in $C_c^\infty(\overline{B^N(a)})$ to some $\varsigma^0_k$ as $j\to \infty$. For every $m$, for every $\epsilon>0$, there exists $N$, for $j,l \geq N$, for every $k$, $\big\| \varsigma^j_k - \varsigma^l_k \big\|_{C^m} <\epsilon$, and thus for $j \geq N$ and for every $k$, $\big\| \varsigma^j_k - \varsigma^0_k \big\|_{C^m} = \lim_{l\to \infty} \big\| \varsigma^j_k - \varsigma^l_k \big\|_{C^m} \leq \epsilon$. Therefore for every $m$,
\begin{align*}
\lim_{j\to\infty} \big\| \{\varsigma^j_k -\varsigma^0_k \}_k\big\|_m = \lim_{j\to \infty} \sup_k \big\| \varsigma^j_k - \varsigma^0_k \big\|_{C^m} = 0.
\end{align*}
For every $k\in \mathbb{N}^\nu$, $1\leq \mu\leq \nu$ satisfying $k_\mu \neq 0$, by the dominated convergence theorem,
\begin{align*}
\int \varsigma^0_k(t) \,dt^\mu \equiv \lim_j \int \varsigma^j_k(t) \,dt^\mu \equiv 0.
\end{align*}
Hence $\{\varsigma^j_k\}$ converges to $\{\varsigma^0_k\}$ in $\overline{\mathcal{K}}(N,e,a,\nu)$.
\end{proof}

\section{Proof of Proposition \ref{1}}\label{A}
The assumptions are that there exist $N,\nu>1,e$ satisfying (\ref{e}), $p\in (1,\infty), n\in \mathbb{N}_{>0}$, and a real analytic function $\gamma_t(x)$ defined on a neighborhood of $(t,x)=(0,0)$ in $\mathbb{R}^N\times \mathbb{R}^n$ with $\gamma_0(x)\equiv x$, such that the operator
$$
Tf(x) = \psi(x)\int f(\gamma_t(x))K(t)\,dt
$$
is bounded on $L^p(\mathbb{R}^n)$ for every $\psi\in C_c^\infty(\mathbb{R}^n)$ with sufficiently small support, and every $K \in \mathcal{K}(N,e,a,\nu)$ with sufficiently small $a>0$, and such that there exists a finite subset $\mathcal{F}_0 \subseteq \mathcal{L}(\mathcal{P})$ that strongly controls every element of $\mathcal{L}(\mathcal{P})$, but at least one element of $\mathcal{L}(\mathcal{P} \cup \mathcal{N})$ is not controlled by $\mathcal{F}_0$. Roughly speaking, the goal is to construct a new $\gamma_t(x)$ whose corresponding operators are still bounded, but all the vector fields associated with pure powers vanish.

\subsection{Preliminaries}\label{natural1}
We want to pick appropriate finite subsets $\mathcal{F} \subseteq \mathcal{L}(\mathcal{P}), \mathcal{F}' \subseteq \mathcal{L}(\mathcal{P} \cup \mathcal{N})$ so that we will only need to work on $\mathcal{F}, \mathcal{F}'$.

\begin{lemma}\label{appropriate}
For any finite subset $\mathcal{F}$ of $\mathcal{L}(\mathcal{P})$ containing $\mathcal{F}_0$, at least one element of $\mathcal{L}(\mathcal{P} \cup \mathcal{N})$ is not controlled by $\mathcal{F}$.
\end{lemma}

\begin{proof}
Suppose for contradiction that there exists a finite subset $\mathcal{F}$ of $\mathcal{L}(\mathcal{P})$ containing $\mathcal{F}_0$ that controls every element of $\mathcal{L}(\mathcal{P} \cup \mathcal{N})$. By Lemma \ref{4}, $\mathcal{F}_0$ controls every element of $\mathcal{F}$. And since $\mathcal{F}_0 \subseteq \mathcal{F}$, by Lemma \ref{5}, $\mathcal{F}_0$ controls every element of $\mathcal{L}(\mathcal{P} \cup \mathcal{N})$; contradiction.
\end{proof}

\begin{proposition}\label{comb}
There exist finitely many elements $(X_1, d_1), \ldots, (X_p,d_p) \in \mathcal{P}, (X_{p+1}, d_{p+1}), \ldots, (X_q,d_q) \in \mathcal{N}$, and finite sets $\mathcal{F}, \mathcal{F}'$ satisfying
\begin{align*}
&\{(X_1, d_1), \ldots, (X_p, d_p)\} \subseteq \mathcal{F} \subseteq \mathcal{L}\big(\{(X_1, d_1), \ldots, (X_p, d_p)\}\big),\\
&\{(X_1, d_1), \ldots, (X_q, d_q)\} \subseteq \mathcal{F}' \subseteq \mathcal{L}\big(\{(X_1, d_1), \ldots, (X_q, d_q)\}\big), \quad \mathcal{F} \subseteq \mathcal{F}',
\end{align*} 
such that $\mathcal{F}$ strongly controls every element of $\mathcal{L}(\mathcal{P})$, $\mathcal{F}'$ strongly controls $W(t,x)$ and every element of $\mathcal{L}(\mathcal{P} \cup \mathcal{N})$, and at least one element of $\mathcal{F}'$ is not controlled by $\mathcal{F}$.

Denote $\mathcal{F}'= (\mathbf{X},\mathbf{d})= \{(X_1, d_1), \ldots, (X_Q,d_Q)\}$ for some $Q\geq q$. Let $\{D_1=(1,0, \ldots, 0), \ldots, D_q= (0, 0, \ldots, 1)\}$ be the standard basis of $\mathbb{R}^q$. Then there exist $D_{q+1}, \ldots, D_Q \in \mathbb{N}^q \backslash \{0\}$, such that for $1\leq i \leq Q$, $d_i = \sum_{j=1}^q D_i^j d_j$, where $D_i^j$ is the $j$-th coordinate of $D_i$, and such that $(\mathbf{X},\mathbf{D}):=\{(X_1,D_1), \ldots, (X_Q,D_Q)\}\subseteq \mathcal{L}\big(\{(X_1, D_1), \ldots, (X_q, D_q)\}\big)$. Moreover, the above $\mathcal{F}'$ can be chosen so that the corresponding $(\mathbf{X}, \mathbf{D})$ strongly controls every element of $\mathcal{L}\big(\{(X_1, D_1), \ldots, (X_q, D_q)\}\big)$. 
\end{proposition}

Recall in Section \ref{sntz}, we choose basis $\widetilde X, \widetilde Y, \widetilde T$ to implement appropriate dilations. Likewise, the choice of $\mathcal{F}$ and $\mathcal{F}'$ is important in our argument, and in particular, we will construct ``dilations'' based on $\mathcal{F}'$. The above $D_1, \ldots, D_Q$ keep track of how we generate $(X_1, d_1), \ldots, (X_Q, d_Q)$ from $(X_1, d_1), \ldots, (X_q, d_q)$ by taking commutators of the vector fields and taking sums of the degrees. For example, if $X_l (q<l\leq Q)$ equals to $[X_1, [X_1,X_2]]$ and $d_l=2d_1 +d_2$, then $D_l = (2,1,0,\ldots,0)\in \mathbb{N}^q$.

\begin{proof}
By Theorem \ref{K}, there exists a finite subset $\mathcal{F}_0'$ of $\mathcal{L}(\mathcal{P}\cup \mathcal{N})$ containing $\mathcal{F}_0$, such that $\mathcal{F}_0'$ strongly controls $W(t,x)$ and every element of $\mathcal{L}(\mathcal{P}\cup \mathcal{N})$. There exist $(X_1, d_1), \ldots, (X_p, d_p) \in \mathcal{P}$, and $(X_{p+1}, d_{p+1}), \ldots, (X_q, d_q) \in \mathcal{N}$, such that $\mathcal{F}_0' \subseteq \mathcal{L}\big(\{(X_1, d_1), \ldots, (X_q, d_q)\}\big)$. For this $q$, let $\{D_1=(1,0, \ldots, 0), \ldots, D_q= (0, 0, \ldots, 1)\}$ be the standard basis of $\mathbb{R}^q$. Then we have a surjection $\Theta$ from $\mathcal{L}\big(\{(X_1, D_1), \ldots, (X_q, D_q)\}\big)$ onto $\mathcal{L}\big(\{(X_1, d_1), \ldots, (X_q, d_q)\}\big)$ by sending each $(X_0, D_0) \in \mathcal{L}\big(\{(X_1, D_1), \ldots, (X_q, D_q)\}\big)$ to $(X_0, \sum_{j=1}^q D_0^j d_j)$, where $D_0^j$ is the $j$-th coordinate of $D_0 \in \mathbb{N}^q \backslash \{0\}$. Note the preimage under $\Theta$ of each element of $\mathcal{L}\big(\{(X_1, d_1), \ldots, (X_q, d_q)\}\big)$ is finite.

By Theorem \ref{K}, there exists a finite set $(\mathbf{X}, \mathbf{D})$ satisfying
$$
\{(X_1, D_1), \ldots, (X_q, D_q)\} \cup \Theta^{-1}(\mathcal{F}_0') \subseteq (\mathbf{X}, \mathbf{D}) \subseteq \mathcal{L}\big(\{(X_1, D_1), \ldots, (X_q, D_q)\}\big),
$$
such that $(\mathbf{X}, \mathbf{D})$ strongly controls every element of $\mathcal{L}\big(\{(X_1, D_1), \ldots, (X_q, D_q)\}\big)$. Let $\mathcal{F}'= \Theta\big((\mathbf{X}, \mathbf{D})\big)$. Then we can write $(\mathbf{X},\mathbf{D}) = \{(X_1, D_1), \ldots, (X_Q, D_Q)\}$, and\footnote{$\mathcal{F}'$ may have cardinality less than $Q$, i.e. may have repeated indices.} $\mathcal{F}' = (\mathbf{X},\mathbf{d})= \{(X_1, d_1), \ldots, (X_q, d_Q)\}$, for some $Q\geq q$. We have $\{(X_1, d_1), \ldots, (X_q, d_q)\}\cup \mathcal{F}_0' \subseteq \mathcal{F}' \subseteq \mathcal{L}\big(\{(X_1, d_1), \ldots, (X_q, d_q)\}\big)$.

Let $\mathcal{F} = \mathcal{F}' \cap \mathcal{L}(\mathcal{P})$. Then $\mathcal{F}_0 \subseteq \mathcal{F}_0' \cap \mathcal{L}(\mathcal{P}) \subseteq \mathcal{F}$, and $\{(X_1, d_1), \ldots, (X_p, d_p)\} \subseteq \mathcal{F} \subseteq \mathcal{L} \big( \{(X_1, d_1), \ldots, (X_p, d_p)\} \big)$. By Theorem \ref{K}, $\mathcal{F}'$ strongly controls $W(t,x)$ and every element of $\mathcal{L}(\mathcal{P}\cup \mathcal{N})$, and $\mathcal{F}$ strongly controls every element of $\mathcal{L}(\mathcal{P})$. By Lemma \ref{appropriate}, at least one element of $\mathcal{L}(\mathcal{P} \cup \mathcal{N})$ is not controlled by $\mathcal{F}$. Then at least one element of $\mathcal{F}'$ is not controlled by $\mathcal{F}$. Otherwise by Lemma \ref{4} and Lemma \ref{5}, $\mathcal{F}$ would control every element of $\mathcal{L}(\mathcal{P}\cup \mathcal{N})$; contradiction.
\end{proof}

We set here the parameters in some of the above statements. There exists $\rho>0$, and a neighborhood $\Omega$ of $0$, such that $(\mathbf{X},\mathbf{D})$ strongly controls every element of $\{([X_i, X_j], D_i + D_j): 1\leq i, j \leq Q\}$ on $\Omega$, $\mathcal{F}'=(\mathbf{X},\mathbf{d})$ strongly controls $W(t,x)$ on $B^N(\rho)\times \Omega$, and $X_l$ has finite $C^m(\Omega)$ norm for every $1\leq l \leq Q, m\in \mathbb{N}$. Assume $\gamma_t(x)$ is $C^\infty$ on a neighborhood of the closure of $B^N(\rho)\times \Omega$, $\gamma_t(\cdot)$ is a diffeomorphism from $\Omega$ onto $\gamma_t(\Omega)$ for $t\in B^N(\rho)$, and all the vector fields in $\mathcal{F}'$ are $C^\infty$ on a neighborhood of the closure of $\Omega$. Assume $\rho>0$ is sufficiently small so that there exists $\xi_1>0$ and a neighborhood $U \subseteq \Omega$ of $0$, such that $\mathcal{F}$ and $\mathcal{F}'$ both satisfy $\mathcal{C}(x, \Vec{\xi_1}, \Omega)$ for all $x\in U$, and such that $U\subseteq \gamma_t(\Omega)$ for $t\in B^N(\rho)$. Then $\gamma_t^{-1}(x)$ is defined and smooth on $B^N(\rho) \times U$.

For $1\leq l, k\leq Q$, we can write
$$
[X_l, X_k] = \sum_{D_i \leq D_l +D_k} c_{l,k}^{i}  X_i, \text{~on~} \Omega,
$$
where the $c_{l,k}^i$ are $C^\infty$ functions on some neighborhood of the closure of $\Omega$. $D_i \leq D_l +D_k$ implies $d_i \leq d_l + d_k$. Thus for all $x\in U, \delta\in [0,1]^\nu$,
$$
[\delta^{d_l} X_l, \delta^{d_k} X_k] = \sum_{d_i \leq d_l + d_k} c_{l,k}^{i} \delta^{d_l + d_k - d_i} \delta^{d_i} X_i, \text{~on~} B_{(\mathbf{X},\mathbf{d})}(x,\xi_1\delta),
$$
and for $d_i \leq d_l + d_k$ and $m\in \mathbb{N}$,
$$
\sup_{\substack{x\in U\\ \delta\in [0,1]^\nu}} \sum_{|\beta|\leq m} \Big\| (\delta \mathbf{X})^\beta (c_{l,k}^{i} \delta^{d_l + d_k - d_i})\Big\|_{C(B_{(\mathbf{X},\mathbf{d})}(x,\xi_1\delta))}\leq \sum_{|\beta|\leq m} \big\| \mathbf{X}^\beta c_{l,k}^i \big\|_{C(\Omega)} <\infty.
$$
Hence $(\mathbf{X},\mathbf{d})$ controls $\big\{\big([X_i, X_j], d_i +d_j\big): 1\leq i,j \leq Q\big\}$ on $U$ with parameter $\xi_1$.

Therefore $\{(\delta^{d_l}X_l, d_l)\}_{l=1}^Q$ satisfies the conditions of Theorem \ref{16} on $B_{(\mathbf{X}, \mathbf{d})}(x, \xi_1 \delta) = B_{(\delta \mathbf{X}, \mathbf{d})}(x, \Vec{\xi_1})$ uniformly in $x\in U, \delta\in [0,1]^\nu$. We will define $m$-admissible constants of this section so that every $m$-admissible constant as in Definition \ref{admissible} with respect to $\{(\delta^{d_l}X_l, d_l)\}_{l=1}^Q$ and $x$, is an $m$-admissible constant of this section, for any $\delta\in [0,1]^\nu$ and $x\in U$.

\begin{definition}\label{uniformadmissible}
Given $(\mathbf{X}, \mathbf{d})$ and $\zeta \in (0,1]$, for $m\geq 2$, we say $C$ is an $m$-admissible constant if $C$ is chosen to depend only on fixed upper and lower bounds for  $\{|d_l|_1: 1\leq l \leq Q\}$, $m$, $\zeta$, fixed lower bound for $\xi_1$, fixed upper bound for $n,Q,\nu$, fixed upper bound for the quantities:
$$
\big\|X_l\big\|_{C^m(\Omega)},\quad  \sum_{|\beta|\leq m} \Big\| \mathbf{X}^\beta c_{l,k}^{i}\Big\|_{C(\Omega)}.
$$
When we wish to make the choice of $(\mathbf{X}, \mathbf{d})$ and $\zeta$ clear, we say, $C$ is an $m$-admissible constant with respect to $(\mathbf{X}, \mathbf{d})$ and $\zeta$.
\end{definition}

In this section, $m$-admissible constants are with respect to the above $(\mathbf{X}, \mathbf{d})$ and $\zeta=1$.

\subsection{Devise a new \texorpdfstring{\boldmath $\gamma_t(x)$}{TEXT}}\label{natural2}
To devise a new $\gamma_t(x)$ whose vector fields associated with pure powers vanish, we need to ``zoom in'' on the original $\gamma_t(x)$ at certain points where its behavior is particularly extreme. Let $\mathcal{F}, \mathcal{F}'=(\mathbf{X}, \mathbf{d})$ be as in Proposition \ref{comb}. The following proposition guides us as to where to find such points.

Let $\delta \mathbf{X}$ be the matrix with column vectors $\{\delta^{d_l} X_l\}_{l=1}^Q$, and let $(\delta \mathbf{X})_{\mathcal{F}}$ be the matrix with column vectors $\{ \delta^{d_l} X_l \}_{(X_l, d_l)\in \mathcal{F}}$. Denote $n_0(x,\delta):=\text{rank}\, (\delta \mathbf{X})(x)$. $|\mathcal{F}|$ denotes the number of elements in $\mathcal{F}$. Let $\mathcal{I}(n_0, \mathcal{F})$ be the set of all lists of integers $(i_1, \ldots, i_{n_0})$ such that $1 \leq i_1 < i_2 < \cdots < i_{n_0} \leq Q$ and $(X_{i_1}, d_{i_1}), \ldots, (X_{i_{n_0}}, d_{i_{n_0}}) \in \mathcal{F}$. See other notations in Definition \ref{indoor}.

\begin{proposition}\label{R}
For every $0<\epsilon \leq 1$, $0<\xi \leq 1$, and every neighborhood $U'\subseteq U$ of $0$ in $\mathbb{R}^n$, there exists $x\in U'$ and $\delta \in [0,\xi]^\nu$ with $n_0(x,\delta)>0$, such that either $
n_0(x,\delta)> |\mathcal{F}|$, or
$$
n_0(x,\delta)\leq |\mathcal{F}|, \quad \Big| \det_{n_0(x,\delta)\times n_0(x,\delta)} (\delta \mathbf{X})_{\mathcal{F}} (x) \Big|_\infty <\epsilon \Big| \det_{n_0(x,\delta)\times n_0(x,\delta)} (\delta \mathbf{X}) (x) \Big|_\infty.
$$
\end{proposition}

To prove Proposition \ref{R}, we need the following lemma.

\begin{lemma}\label{7} (Lemmas 3.6, 3.8, 4.12, and 4.13 in \cite{CARNOT})
\begin{enumerate}[label=(\roman*)]
\item 
For every $x\in U$, $\delta\in [0,1]^\nu$ with $n_0(x,\delta)>0$, and every $I\in \mathcal{I}(n_0(x,\delta),n), J \in \mathcal{I}(n_0(x,\delta),Q)$, $m\in \mathbb{N}$, we have on $B_{(\mathbf{X},\mathbf{d})}(x,\xi_1 \delta)$, $\text{rank}\,(\delta \mathbf{X})$ is constant, and
\begin{align*}
&\sum_{|\beta|\leq m} \Big| (\delta \mathbf{X})^\beta \det (\delta \mathbf{X})_{I \times J}  \Big| \lesssim_{m\vee 2} \Big| \det_{n_0(x,\delta)\times n_0(x,\delta)} (\delta \mathbf{X}) \Big|,\\
&\Big| \det_{n_0(x,\delta)\times n_0(x,\delta)} (\delta \mathbf{X})\Big| \approx \Big| \det_{n_0(x,\delta)\times n_0(x,\delta)} (\delta \mathbf{X})(x)\Big|_\infty \approx \Big| \det_{n_0(x,\delta)\times n_0(x,\delta)} (\delta \mathbf{X})(x)\Big|.
\end{align*}
\item
Denote by $\big(I(\epsilon, x,\delta), J(\epsilon, x,\delta)\big)$ any element in $\mathcal{I}(n_0(x,\delta), n)\times \mathcal{I}(n_0(x,\delta), Q)$ such that 
$$
\Big| \det (\delta \mathbf{X})_{I(\epsilon, x,\delta) \times J(\epsilon, x,\delta)} (x) \Big| \geq \epsilon \Big| \det_{n_0(x,\delta)\times n_0(x,\delta)} (\delta \mathbf{X})(x) \Big|_\infty.
$$
For every $0<\epsilon \leq 1$, there exists a constant $\widetilde{\xi_1}(\epsilon): 0<\widetilde{\xi_1}(\epsilon) \leq \xi_1$, depending only on $\epsilon$ and the quantities that a $2$-admissible constant is allowed to depend on, such that for every $x\in U, \delta \in [0,1]^\nu$ with $n_0(x,\delta)>0$, and every $I\in \mathcal{I}(n_0(x,\delta),n), J \in \mathcal{I}(n_0(x,\delta),Q)$, $m\in \mathbb{N}$,
\begin{align*}
&\Big| \det (\delta \mathbf{X})_{I(\epsilon, x,\delta) \times J(\epsilon, x,\delta)}  \Big| \gtrsim \epsilon \Big| \det_{n_0(x,\delta)\times n_0(x,\delta)} (\delta \mathbf{X})\Big|, \text{~on~} B_{(\mathbf{X},\mathbf{d})}(x, \widetilde{\xi_1}(\epsilon) \delta),\\
&\sum_{|\beta|\leq m}\Big\| (\delta \mathbf{X})^\beta \Big(\frac{\det (\delta \mathbf{X})_{I\times J}}{\det (\delta \mathbf{X})_{I(\epsilon, x,\delta) \times J(\epsilon, x,\delta)}}\Big) \Big\|_{C(B_{(\mathbf{X},\mathbf{d})}(x, \widetilde{\xi_1}(\epsilon) \delta) )} \lesssim_{m\vee 2,\epsilon} 1,
\end{align*}
where $\lesssim_{m\vee 2,\epsilon} 1$ means less than a constant which only depends on $\epsilon$ and the quantities that an $(m\vee 2)$-admissible constant is allowed to depend on.
\end{enumerate}
\end{lemma}

Without loss of generality, we may shrink $\xi_1$ to be equal to $\widetilde{\xi_1}(1)$ so that the conclusions of Lemma \ref{7} hold with $\widetilde{\xi_1}(1)$ replaced by $\xi_1$.

\begin{proof}[\textbf{Proof of Proposition \ref{R}}]

Suppose by contradiction that there exists $0<\epsilon \leq 1, 0<\xi \leq 1$, and a neighborhood $U'\subseteq U$ of $0$, such that for all $x\in U', \delta \in [0,\xi]^\nu$ with $n_0(x,\delta)>0$, we have 
$$
n_0(x,\delta)\leq |\mathcal{F}|, \quad \Big| \det_{n_0(x,\delta)\times n_0(x,\delta)} (\delta \mathbf{X})_{\mathcal{F}} (x) \Big|_\infty \geq \epsilon \Big| \det_{n_0(x,\delta)\times n_0(x,\delta)} (\delta \mathbf{X}) (x) \Big|_\infty.
$$
Then we can pick $J(\epsilon, x,\delta) \in \mathcal{I}(n_0(x,\delta), \mathcal{F})$. By Lemma \ref{7}, 
$$
\Big| \det (\delta \mathbf{X})_{I(\epsilon, x,\delta) \times J(\epsilon, x,\delta)}  \Big| \gtrsim \epsilon \Big| \det_{n_0(x,\delta)\times n_0(x,\delta)} (\delta \mathbf{X})\Big|, \text{~on~} B_{(\mathbf{X},\mathbf{d})}(x, \widetilde{\xi_1}(\epsilon) \delta),
$$
and for every $I\in \mathcal{I}(n_0(x,\delta), n), J\in \mathcal{I}(n_0(x,\delta), Q)$, and $m\in \mathbb{N}$,
$$
\sum_{|\beta|\leq m}\Big\| (\delta \mathbf{X})^\beta \Big(\frac{\det (\delta \mathbf{X})_{I\times J}}{\det (\delta \mathbf{X})_{I(\epsilon, x,\delta) \times J(\epsilon, x,\delta)}}\Big) \Big\|_{C(B_{(\mathbf{X},\mathbf{d})}(x, \widetilde{\xi_1}(\epsilon) \delta) )} \lesssim_{m\vee 2,\epsilon} 1.
$$

By Cramer's rule, for $x\in U'$ and $\delta\in [0,\xi]^\nu$ with $n_0(x,\delta)>0$, and for $1\leq k \leq Q$, 
\begin{equation*}
\delta^{d_k} X_k =\sum_{l\in J(\epsilon, x,\delta)} g_k^{l,x,\delta} \delta^{d_l} X_l,  \text{~on~} B_{(\mathbf{X},\mathbf{d})}(x, \widetilde{\xi_1}(\epsilon) \delta),
\end{equation*}
where
$$
g_k^{l,x,\delta}= \frac{\det (\delta \mathbf{X})_{I(\epsilon, x,\delta) \times J(\epsilon, x,\delta)}^{l,k}}{\det (\delta \mathbf{X})_{I(\epsilon, x,\delta) \times J(\epsilon, x,\delta)}},
$$
and $(\delta \mathbf{X})_{I(\epsilon, x,\delta) \times J(\epsilon, x,\delta)}^{l,k}$ denotes the matrix obtained from $(\delta \mathbf{X})_{I(\epsilon, x,\delta) \times J(\epsilon, x,\delta)}$ via replacing the column vector $(\delta^{d_l}X_l)_{I(\epsilon, x,\delta)}$ by $(\delta^{d_k}X_k)_{I(\epsilon, x,\delta)}$. Note $\det (\delta \mathbf{X})_{I(\epsilon, x,\delta) \times J(\epsilon, x,\delta)}^{l,k}$ is either $0$ or $\pm \det (\delta \mathbf{X})_{I(\epsilon, x,\delta) \times J'}$ for some $J'\in \mathcal{I}(n_0( x,\delta),Q)$. Thus
\begin{align*}
\sum_{|\beta|\leq m} \Big\| \big((\delta \mathbf{X})_{\mathcal{F}}\big)^\beta g_k^{l,x,\delta}\Big\|_{C(B_{\mathcal{F}}(x, \widetilde{\xi_1}(\epsilon) \delta))} &\leq \sum_{|\beta|\leq m} \Big\| \big(\delta \mathbf{X}\big)^\beta g_k^{l,x,\delta}\Big\|_{C(B_{(\mathbf{X},\mathbf{d})}(x, \widetilde{\xi_1}(\epsilon) \delta))} \\
&\lesssim_{m\vee 2,\epsilon} 1.
\end{align*}
Write $\delta = \xi\delta'$. Then for $x\in U'$ and $\delta'\in [0,1]^\nu$ with $n_0(x,\delta') >0$ (note $n_0(x,\delta')=n_0(x,\delta))$,
$$
\delta'{}^{d_k}X_k= \sum_{l\in J(\epsilon, x,\delta)} \big(\xi^{|d_l|_1-|d_k|_1}g_k^{l,x,\delta}\big) \delta'{}^{d_l}X_l, \text{~on~} B_{\mathcal{F}}(x, \widetilde{\xi_1}(\epsilon)\xi \delta'),
$$
and
$$
\sum_{|\beta|\leq m}\Big\| \big( (\delta'\mathbf{X})_{\mathcal{F}}\big)^{\beta} \big( \xi^{|d_l|_1-|d_k|_1}g_k^{l,x,\delta}\big)\Big\|_{C(B_{\mathcal{F}}(x, \widetilde{\xi_1}(\epsilon) \xi \delta'))} \lesssim_{m\vee 2,\epsilon,\xi} 1,
$$
where $\lesssim_{m\vee 2,\epsilon, \xi}$ means less than a constant which only depends on $\epsilon, \xi$ and the quantities that an $(m\vee 2)$-admissible constant is allowed to depend on. For $x\in U'$ and $\delta'\in [0,1]^\nu$ with $n_0(x,\delta')=0$,
$$
{\delta'}^{d_k}X_k=0, \text{~on~} B_{\mathcal{F}}(x, \widetilde{\xi_1}(\epsilon) \xi \delta') =\{x\}.
$$
Hence $\mathcal{F}$ controls every $(X_k,d_k)\in \mathcal{F}'$ on $U'$ with parameter $\widetilde{\xi_1}(\epsilon) \xi$; contradiction. 
\end{proof}

By Proposition \ref{R}, either there exist sequences $x_j\to 0$ in $U, \delta_j \to 0$ coordinate-wise in $[0,1]^\nu$, such that $n_0(x_j,\delta_j) > |\mathcal{F}|$ for every $j$, or there exist sequences $\epsilon_j\to 0 (0<\epsilon_j \leq 1)$, $x_j\to 0$ in $U$, $\delta_j\to 0$ coordinate-wise in $[0,1]^\nu$, such that for every $j$,
$$
0<n_0(x_j,\delta_j)\leq |\mathcal{F}|, \quad \Big| \det_{n_0(x_j,\delta_j)\times n_0(x_j,\delta_j)} (\delta_j \mathbf{X})_{\mathcal{F}} (x_j) \Big| <\epsilon_j \Big| \det_{n_0(x_j,\delta_j)\times n_0(x_j,\delta_j)} (\delta_j \mathbf{X}) (x_j) \Big|.
$$
Such points $x_j$ are where we will ``zoom in'' on the $\gamma_t(x)$. We will get the ``zoom ratio'' by modifying numbers $\delta_j^{-d_l}$ for $1\leq l \leq Q$.

Since $\mathcal{F}'$ is finite, we can find positive numbers $b_1, \ldots, b_\nu$ such that the linear function $\theta(d) = \theta(d^1, \ldots, d^\nu) := \sum_{\mu=1}^\nu b_\mu d^\mu$ separates\footnote{This means for any $(X_j,d_j),(X_k,d_k)\in \mathcal{F}'$ such that $d_j \neq d_k$, we have $\theta(d_j) \neq \theta(d_k)$.} different degrees in $\mathcal{F}'$. For $1\leq l\leq p$, define $\theta_l = \frac{1}{2} \theta(d_l)$. For $p< l \leq q$, define $\theta_l = \theta(d_l)$. For $q< l \leq Q$, define $\theta_l = \sum_{i=1}^q D_l^i \theta_i$, where $D_l^i$ is the $i$-th coordinate of $D_l$.

For $1\leq l \leq Q$, since $d_l = D_l^1 d_1 + \cdots + D_l^q d_q$, we have $\theta(d_l) = D_l^1 \theta(d_1) + \cdots + D_l^q \theta(d_q)$, and thus $\theta_l = D_l^1 \theta_1 + \cdots + D_l^q \theta_q \leq \theta(d_l)$, where the equality holds if and only if $D_l^1 = \cdots = D_l^p = 0$. Hence for $1\leq l \leq Q$, $0< \theta_l \leq \theta(d_l)$, and $\theta_l = \theta(d_l)$ if and only if $(X_l, d_l)$ is generated by $\{(X_{p+1}, d_{p+1}), \ldots, (X_q, d_q)\}$.

\begin{remark}\label{parallel}
\begin{enumerate}
\item 
Recall in Section \ref{sntz}, we ``zoom in'' on the $\gamma_s(\xi)$ based at the origin with ``zoom ratios'' $\tau_1$ and $\tau_2$ along ``directions'' $\widetilde X$ and $\widetilde Y$, respectively, using the dilation $\Phi_{\tau_1, \tau_2}$. The situation is somewhat similar here. Let $\tau_j$ be large numbers to be chosen later. We will ``zoom in'' on the $\gamma_t(x)$ based at $x_j$ with ``zoom ratio'' $\tau_j^{\theta_l} \delta_j^{-d_l}$ along ``direction'' $X_l$ for $1\leq l \leq Q$, using the local diffeomorphism $\Phi$ of Theorem \ref{16}.
\item
In Section \ref{sntz}, the ``zoom in'' via $\Phi_{\tau_1, \tau_2}$ causes a dilation on each vector field $X_\alpha$ by $\tau_1, \tau_2$, or $\tau_1 \tau_2$; the kernel $\sum_{k\in \mathbb{N}} \varsigma^{\big((\tau_0^{b_1}, \tau_0^{b_2})2^{k\textbf{n}}\big)}(s)$ causes a dilation on each vector field $X_\alpha$ by $\tau_0^{-(b_1\alpha^1 + b_2 \alpha^2)}$, where $\alpha=(\alpha^1, \alpha^2)$ is the multi-index, and large numbers $\tau_0, \tau_1, \tau_2$ are chosen so that the net dilations on certain vector fields $X_\alpha$ are small and thus these vector fields vanish in the limiting process. Here, the ``zoom in'' via the $\Phi$ of Theorem \ref{16} will cause a dilation on each $X_l$ by $\tau_j^{\theta_l} \delta_j^{-d_l}$; we will choose a kernel that can cause a dilation on each $X_l$ by $\tau_j^{-\theta(d_l)} \delta_j^{d_l}$. Since $\theta_l \leq \theta(d_l)$ for all $l$, and $\theta_l < \theta(d_l)$ for some $l$, the net dilations $\tau_j^{\theta_l - \theta(d_l)}$ on $X_l$ for certain $l$'s are small, and thus these $X_l$'s will vanish in a limiting process elaborated later.
\end{enumerate}
\end{remark}

To make sure that for every $\tau(\geq1),j$, the $\Phi$ constructed with ``zoom ratio'' $\tau^{\theta_l} \delta_j^{-d_l}$ along ``direction'' $X_l$ for $1\leq l \leq Q$ is a valid ``zoom in'', we now verify that $\{(\tau^{-\theta_l} \delta_j^{d_l} X_l, d_l)\}_{1\leq l \leq Q}$ satisfies the conditions of Theorem \ref{16} for every $x_j, \delta_j, \tau(\geq1)$.
Note $D_i \leq D_l + D_k$ implies $\theta_i = D_i \cdot (\theta_1, \ldots, \theta_q) \leq (D_l + D_k) \cdot (\theta_1, \ldots, \theta_q) = \theta_l + \theta_k$, and $d_i = D_i^1 d_1 + \cdots + D_i^q d_q \leq (D_l^1 + D_k^1) d_1 + \cdots + (D_l^q + D_k^q) d_q = d_l + d_k$. Thus for $1\leq l, k \leq Q$,
$$
[X_l, X_k] = \sum_{\substack{d_i \leq d_l + d_k \\ \theta_i \leq \theta_l +\theta_k}} c_{l,k}^{i}  X_i, \text{~on~} \Omega,
$$
where the $c_{l,k}^i$ are $C^\infty$ functions on some neighborhood of the closure of $\Omega$. Denote $\tau^{-1} \delta_j \mathbf{X} := \{\tau^{-\theta_l} \delta_j^{d_l} X_l\}_{1\leq l \leq Q}$, and $(\tau^{-1} \delta_j \mathbf{X}, \mathbf{d}):=\{(\tau^{-\theta_l} \delta_j^{d_l} X_l, d_l)\}_{1\leq l \leq Q}$. Then for every $\tau \geq 1$ and every $j$, $(\tau^{-1} \delta_j \mathbf{X}, \mathbf{d})$ satisfies $\mathcal{C}(x_j, \Vec{\xi_1}, \Omega)$, and for $1\leq l,k \leq Q$,
\begin{equation}\label{commutator}
[\tau^{-\theta_l}\delta_j^{d_l} X_l,\tau^{-\theta_k}\delta_j^{d_k} X_k] = \sum_{\substack{d_i \leq d_l + d_k \\ \theta_i \leq \theta_l +\theta_k}} \tau^{-\theta_l-\theta_k + \theta_i} \delta_j^{d_l + d_k -d_i} c_{l,k}^{i} \tau^{-\theta_i} \delta_j^{d_i}  X_i, \text{~on~} B_{(\tau^{-1} \delta_j \mathbf{X},\mathbf{d})}(x_j,\Vec{\xi_1}), 
\end{equation}
and for $d_i \leq d_l + d_k, \theta_i \leq \theta_l +\theta_k$, and $m\in \mathbb{N}$,
\begin{equation}\label{coefficient}
\sup_{\tau\geq1,j} \sum_{|\beta|\leq m} \Big\|(\tau^{-1} \delta_j \mathbf{X})^\beta (\tau^{-\theta_l-\theta_k + \theta_i} \delta_j^{d_l + d_k -d_i} c_{l,k}^{i})\Big\|_{C(B_{(\tau^{-1} \delta_j \mathbf{X}, \mathbf{d})}(x_j,\Vec{\xi_1}))}<\infty.
\end{equation}
So $(\tau^{-1} \delta_j \mathbf{X}, \mathbf{d})$ satisfies the conditions of Theorem \ref{16} uniformly in $x_j, \delta_j, \tau(\geq1)$.

We now choose the sequence $\{\tau_j\}_j$. In the case that $n_0(x_j,\delta_j) > |\mathcal{F}|$ for each $j$, just let $1\leq \tau_j \to \infty$. In the case that $n_0(x_j,\delta_j) \leq |\mathcal{F}|$ for each $j$, let $1\leq \tau_j \to \infty$ such that
$$
\tau_j^{n\max_l \theta_l} \leq \frac{1}{\sqrt{\epsilon_j}},
$$
and then
\begin{align*}
\Big| \det_{n_0(x_j,\delta_j)\times n_0(x_j,\delta_j)} (\tau_j^{-1} \delta_j \mathbf{X})_{\mathcal{F}} (x_j) \Big| &\leq \tau_j^{-n_0(x_j,\delta_j) \min_l(\theta_l)} \Big| \det_{n_0(x_j,\delta_j)\times n_0(x_j,\delta_j)} (\delta_j \mathbf{X})_{\mathcal{F}} (x_j) \Big| \\
&\leq\tau_j^{-n_0(x_j,\delta_j) \min_l(\theta_l)} \epsilon_j \Big| \det_{n_0(x_j,\delta_j)\times n_0(x_j,\delta_j)} (\delta_j \mathbf{X}) (x_j) \Big| \\
&\leq \tau_j^{n_0(x_j,\delta_j)(\max_l \theta_l - \min_l \theta_l)} \epsilon_j  \Big| \det_{n_0(x_j,\delta_j)\times n_0(x_j,\delta_j)} (\tau_j^{-1}  \delta_j \mathbf{X}) (x_j) \Big| \\
&\leq \sqrt{\epsilon_j} \Big| \det_{n_0(x_j,\delta_j)\times n_0(x_j,\delta_j)} (\tau_j^{-1}  \delta_j \mathbf{X}) (x_j) \Big|,
\end{align*}
where $(\tau_j^{-1} \delta_j \mathbf{X})_{\mathcal{F}}$ denotes the matrix with column vectors $\{ \tau_j^{-\theta_l} \delta_j^{d_l} X_l\}_{(X_l, d_l) \in \mathcal{F}}$.

$(\tau_j^{-1} \delta_j \mathbf{X},\mathbf{d})$ satisfies the conditions of Theorem \ref{16} uniformly in $j$. Since $\tau_j^{-\theta_l}\neq 0$, $\text{rank}\, (\tau_j^{-1} \delta_j \mathbf{X}) (x_j) = n_0(x_j,\delta_j)$. Let $J_j\in \mathcal{I}(n_0(x_j,\delta_j),Q)$ be such that
$$
\Big| \det_{n_0(x_j,\delta_j)\times n_0(x_j,\delta_j)} (\tau_j^{-1}  \delta_j \mathbf{X})_{J_j} (x_j) \Big|_\infty = \Big| \det_{n_0(x_j,\delta_j)\times n_0(x_j,\delta_j)} (\tau_j^{-1}  \delta_j \mathbf{X}) (x_j) \Big|_\infty.
$$
In the case that $n_0(x_j,\delta_j) > |\mathcal{F}|$ for every $j$, it is clear that $J_j \not\subseteq \mathcal{F}$. In the case that $n_0(x_j,\delta_j) \leq |\mathcal{F}|$ for every $j$, since $\sqrt{\epsilon_j} \to 0$, $J_j \not\subseteq \mathcal{F}$ for $j\gg1$.

For $D_i \leq D_l +D_k$,
\begin{align*}
&\sup_j \Big\| \tau_j^{-\theta_l} \delta_j^{d_l} X_l \Big\|_{C^m( B_{(\tau_j^{-1} \delta_j \mathbf{X},\mathbf{d})}(x_j, \Vec{\xi_1}))}\leq \big\|X_l\big\|_{C^m(\Omega)} <\infty,\\
& \sup_j \sum_{|\beta|\leq m} \Big\|(\tau_j^{-1} \delta_j \mathbf{X})^\beta (\tau_j^{-\theta_l-\theta_k + \theta_i} \delta_j^{d_l + d_k -d_i} c_{l,k}^{i})\Big\|_{C(B_{(\tau_j^{-1} \delta_j \mathbf{X},\mathbf{d})}(x_j, \Vec{\xi_1}))} \leq  \sum_{|\beta|\leq m} \big\| \mathbf{X}^\beta c_{l,k}^{i}\big\|_{C(\Omega)}< \infty.
\end{align*}
By Theorem \ref{16}, there exist 2-admissible constants $\eta_1>0, 0<\xi_2\leq \xi_1$, such that for every $j$, we can define a diffeomorphism $\Phi_j: B^{n_0(x_j,\delta_j)}(\eta_1) \to \Phi_j(B^{n_0(x_j,\delta_j)}(\eta_1))\subseteq \mathbb{R}^n$ as
$$
\Phi_j(u)=e^{u\cdot (\tau_j^{-1} \delta_j \mathbf{X})_{J_j}}x_j,
$$
and
$$
B_{(\tau_j^{-1}\delta_j \mathbf{X},\mathbf{d})}(x_j, \Vec{\xi_2}) \subseteq \Phi_j (B^{n_0(x_j,\delta_j)}(\eta_1))\subseteq B_{(\tau_j^{-1}\delta_j \mathbf{X},\mathbf{d})}(x_j, \Vec{\xi_1}),
$$
where $u\cdot (\tau_j^{-1} \delta_j \mathbf{X})_{J_j} = \sum_{l\in J_j} u_l \tau_j^{-\theta_l} \delta_j^{d_l} X_l$.
If we let $Y^j_l$ be the pullback of $\tau_j^{-\theta_l} \delta_j^{d_l} X_l$ under the map $\Phi_j$ to $B^{n_0(x_j,\delta_j)}(\eta_1)$, then for every $j$, $1\leq l \leq Q$, and $m\in \mathbb{N}$, on $B^{n_0(x_j,\delta_j)}(\eta_1)$,
\begin{align*}
\Big| \det_{n_0(x_j,\delta_j)\times n_0(x_j,\delta_j)} \mathbf{Y}^j \Big|\approx \big| \det \mathbf{Y}^j_{J_j}\big| \approx 1, \quad \big\| Y_l^j\big\|_{C^m(B^{n_0(x_j,\delta_j)}(\eta_1))}\lesssim_{m\vee 2} 1,
\end{align*}
where $\mathbf{Y}^j$ is the matrix with column vectors $\{\mathbf{Y}^j_l\}_{1\leq l \leq Q}$, and $Y_{J_j}^j$ is the matrix with column vectors $\{Y^j_l\}_{l\in J_j}$. For $1\leq l \leq Q$, let $Z_l^j$ be the pullback of $\tau_j^{-\theta(d_l)} \delta_j^{d_l} X_l$ under the map $\Phi_j$ to $B^{n_0(x_j,\delta_j)}(\eta_1)$. Then we also have that $Z_l^j$ are bounded in $C^m(B^{n_0(x_j,\delta_j)}(\eta_1))$ for all $m$ uniformly in $j$. Since there are only finitely many choices for $n_0(x_j,\delta_j)$ and $J_j$, there exists a subsequence of $\{(x_j,\delta_j, \tau_j)\}_j$, still denoted as $\{(x_j,\delta_j, \tau_j)\}_j$, such that the corresponding $n_0(x_j,\delta_j)$, $J_j$ are fixed, denoted as $\bar n$, $J$, respectively. Since $J=J_j \not\subseteq \mathcal{F}$, there exists 
\begin{equation}\label{prairie}
p_0\in J \backslash \mathcal{F}, \text{ and thus } |Y^j_{p_0}| \approx 1, \text{~on~} B^{\bar n}(\eta_1) \text{ for all } j.
\end{equation}

Denote $\tilde \tau_j := (\tau_j^{b_1}, \ldots, \tau_j^{b_\nu})$. Then for $\delta_j = (\delta_j^1, \ldots, \delta_j^\nu)$, $\tilde \tau_j^{-1} \delta_j t = (\tau_j^{-\theta(e_1)} \delta_j^{e_1} t_1, \ldots,\tau_j^{-\theta(e_N)} \delta_j^{e_N} t_N)$. Let $(\tilde \tau_j^{-1} \delta_j \mathbf{X}, \mathbf{d})= \{ (\tau_j^{-\theta(d_l)} \delta_j^{d_l} X_l, d_l) \}_{1\leq l \leq Q}$. Note $(\tilde \tau_j^{-1} \delta_j \mathbf{X}, \mathbf{d})$ is different from $(\tau_j^{-1} \delta_j \mathbf{X},\mathbf{d})$.

\begin{proposition}\label{ppp}
There exists a subsequence of the above $\{(x_j,\delta_j, \tau_j)\}_j$, still denoted as $\{(x_j,\delta_j, \tau_j)\}_j$, such that for each $1\leq \mu \leq \nu$, $\{\delta_j^\mu\}_j$ either never equals $0$ or always equals $0$, and such that $W^j(t):=\Phi_{j}^*(W( \tilde \tau_j^{-1} \delta_j t))$ converges in $C^\infty(B^N(\rho)\times B^{\bar n}(\eta_1))$ to some $\overline W(t)$, and such that associated to each $(X_\alpha, d_\alpha)\in \mathcal{L}(\mathcal{P}\cup \mathcal{N})$, $Z^{j}_\alpha: = \Phi_j^*(\tau_j^{-\theta(d_\alpha)}\delta_j^{d_\alpha} X_\alpha)$ converges in $C^\infty(B^{\bar n}(\eta_1))$ to some $Z_\alpha$.
\end{proposition}

\begin{proof}
It is clear we can take a subsequence of the above $\{(x_j,\delta_j, \tau_j)\}_j$, still denoted as $\{(x_j,\delta_j, \tau_j)\}_j$, such that for each $1\leq \mu \leq \nu$, $\{\delta_j^\mu\}_j$ either never equals $0$ or always equals $0$.

Since $(\mathbf{X},\mathbf{d})$ strongly controls $W(t,x)$ on $B^N(\rho)\times \Omega$, we have on $B^N(\rho)\times B^{\bar n}(\eta_1)$,
\begin{align*}
W^j(t,u) &= D\Phi_j(u)^{-1}\Big(W\big(\tilde \tau_j^{-1} \delta_j t, \Phi_j(u)\big)\Big)\\ 
&= D\Phi_j(u)^{-1}\Big(\sum_{\substack{ 1\leq l \leq Q\\ \text{deg}\,(\alpha) = d_l}}  t^{\alpha}\, c_{l, \alpha}\big(\tilde \tau_j^{-1} \delta_j t, \Phi_j(u)\big) ~\tau_j^{-\theta(d_l)} \delta_j^{d_l}  X_l(\Phi_j(u)) \Big)  \\
&= \sum_{\substack{ 1\leq l \leq Q\\ \text{deg}\,(\alpha) = d_l}} t^{\alpha} \, c_{l,\alpha} \big(\tilde \tau_j^{-1} \delta_j t, \Phi_j(u)\big) \,Z_l^j(u),
\end{align*}
where the $c_{l,\alpha}$ are $C^\infty$ on a neighborhood of the closure of $B^N(\rho)\times \Omega$. By Theorem \ref{16}, for all $m\in \mathbb{N}$,
\begin{align*}
&\quad \sup_j \Big\| c_{l,\alpha}\big(\tilde \tau_j^{-1} \delta_j t, \Phi_j(u)\big) \Big\|_{C^m(B^N(\rho)\times B^{\bar n}(\eta_1))}\\
&\lesssim_{m\vee 2} \sup_j \sum_{|\beta| + |\beta'| \leq 2m} \Big\| (\mathbf{Y}^j)^\beta (\partial_t^{\beta'} c_{l,\alpha})\big(\tilde \tau_j^{-1} \delta_j t, \Phi_j(u)\big) \Big\|_{C(B^N(\rho)\times B^{\bar n}(\eta_1))}\\
&\leq \sup_j \sum_{|\beta| + |\beta'|\leq 2m} \Big\| (\tau_j^{-1} \delta_j \mathbf{X})^\beta (\partial_t^{\beta'} c_{l, \alpha})\big(\tilde \tau_j^{-1} \delta_j t,x\big) \Big\|_{C(B^N(\rho)\times \Omega)} <\infty.
\end{align*}
Hence $\{W^j\}_j$ is bounded in $C^\infty(B^N(\rho)\times B^{\bar n}(\eta_1))$, and thus there exists a subsequence of $\{(x_j,\delta_j, \tau_j)\}_j$, still denoted as $\{(x_j,\delta_j, \tau_j)\}_j$, such that $W^{j}$ converges in $C^\infty(B^N(\rho)\times B^{\bar n}(\eta_1))$ to some $\overline W$.

Fix an arbitrary $(X_\alpha, d_\alpha)\in \mathcal{L}(\mathcal{P} \cup \mathcal{N})$. $(\mathbf{X},\mathbf{d})$ strongly controls $(X_\alpha, d_\alpha)$ on some neighborhood $U_\alpha \subseteq \Omega$ of $0$. There exists $N_\alpha\in \mathbb{N}$, for $j\geq N_\alpha$, $\Phi_j(B^{\bar n}(\eta_1)) \subseteq B_{(\tau_j^{-1} \delta_j \mathbf{X}, \mathbf{d})} (x_j, \Vec{\xi_1}) \subseteq U_\alpha$. Then for $j\geq N_\alpha$, on $B^{\bar n}(\eta_1)$,
\begin{align*}
Z_\alpha^j &= \Phi_j^*( \tau_j^{-\theta(d_\alpha)} \delta_j^{d_\alpha} X_\alpha) = \Phi_j^*\Big(\sum_{d_l \leq d_\alpha} \tau_j^{-\theta(d_\alpha) + \theta(d_l)} \,\delta_j^{d_\alpha -d_l} \,c_{\alpha, l}~ \tau_j^{-\theta(d_l)} \delta_j^{d_l} X_l\Big)\\
&= \sum_{d_l \leq d_\alpha} \tau_j^{-\theta(d_\alpha) + \theta(d_l)} \,\delta_j^{d_\alpha -d_l} \,(c_{\alpha, l} \circ \Phi_j) \,Z_l^j,
\end{align*}
where the $c_{\alpha, l}$ are $C^\infty$ on a neighborhood of the closure of $U_\alpha$. By Theorem \ref{16}, for all $m\in \mathbb{N}$ and $d_l \leq d_\alpha$, 
\begin{align*}
\sup_{j\geq N_\alpha} \Big\| \tau_j^{-\theta(d_\alpha) + \theta(d_l)} \,\delta_j^{d_\alpha -d_l} \,(c_{\alpha, l} \circ \Phi_j) \Big\|_{C^m(B^{\bar n}(\eta_1))} &\lesssim_{m\vee 2} \sup_{j\geq N_\alpha} \sum_{|\beta|\leq m} \Big\| (\mathbf{Y}^j)^\beta (c_{\alpha, l} \circ \Phi_j) \Big\|_{C(B^{\bar n}(\eta_1))}\\
&\leq \sup_{j\geq N_\alpha} \sum_{|\beta|\leq m} \Big\| (\tau_j^{-1} \delta_j \mathbf{X})^\beta c_{\alpha, l}\Big\|_{C(U_\alpha)}<\infty.
\end{align*}
Hence associated with any $(X_\alpha, d_\alpha) \in \mathcal{L}(\mathcal{P}\cup \mathcal{N})$, $\{Z_\alpha^j\}_j$ is bounded in $C^\infty(B^{\bar n}(\eta_1))$.

Denote $\mathcal{L}(\mathcal{P}\cup \mathcal{N}) = \{(X_{\alpha_r}, d_{\alpha_r})\}_{r=1}^\infty$. We are going to obtain a subsequence of $\{(x_j, \delta_j, \tau_j)\}_j$ by the diagonal argument, so that the corresponding subsequence of $\{Z_{\alpha_r}^j\}_j$ converges in $C^\infty(B^{\bar n}(\eta_1))$ for each $r$. Since $\{Z_{\alpha_1}^{j}\}_{j}$ is bounded in $C^\infty(B^{\bar n}(\eta_1))$, there exists a subsequence $\{(x_{1,j},\delta_{1,j}, \tau_{1,j})\}_j$ of $\{(x_j,\delta_j, \tau_j)\}_j$  such that the corresponding subsequence $\{Z_{\alpha_1}^{1,j}\}_{j}$ of $\{Z^j_{\alpha_1}\}_{j}$ converges in $C^\infty(B^{\bar n}(\eta_1))$ to some $Z_{\alpha_1}$. Since $\{Z_{\alpha_2}^{j}\}_{j}$ is bounded in $C^\infty(B^{\bar n}(\eta_1))$, $\{Z_{\alpha_2}^{1,j}\}_{j}$ is also bounded in $C^\infty(B^{\bar n}(\eta_1))$. Thus there exists a further subsequence $\{(x_{2,j},\delta_{2,j}, \tau_{2,j})\}_j$ of $\{(x_{1,j},\delta_{1,j}, \tau_{1,j})\}_j$ such that the corresponding subsequence $\{Z_{\alpha_2}^{{2,j}}\}_{j}$ of $\{Z_{\alpha_2}^{1,j}\}_{j}$ converges in $C^\infty(B^{\bar n}(\eta_1))$ to some $Z_{\alpha_2}$. Proceeding in an inductive fashion, we obtain a subsequence $\{(x_{r+1,j},\delta_{r+1,j}, \tau_{r+1,j})\}_j$ of the previous subsequence $\{(x_{r,j},\delta_{r,j}, \tau_{r,j})\}_j$ such that the corresponding subsequence $\{Z_{\alpha_{r+1}}^{r+1, j}\}_j$ of $\{Z_{\alpha_{r+1}}^{r,j}\}_{j}$ converges in $C^\infty(B^{\bar n}(\eta_1))$ to some $Z_{\alpha_{r+1}}$. Note the superscript of $Z_{\alpha_r}^{i,j}$ means it is associated to $(x_{i,j}, \delta_{i,j}, \tau_{i,j})$, which is the $j$-th element in the $i$-th subsequence. Now we take the diagonal elements of this sequence of subsequences $\big\{ \{(x_{i,j}, \delta_{i,j}, \tau_{i,j})\}_j \big\}_i$, and obtain a new subsequence $\{(x_{j,j},\delta_{j,j}, \tau_{j,j})\}_j$ of the original sequence. Therefore for each $r$, $\{Z_{\alpha_r}^{j,j}\}_j$ converges to $Z_{\alpha_r}$ in $C^\infty(B^{\bar n}(\eta_1))$.
\end{proof}

Write the Taylor expansion of $\overline{W}(t,u)$ in $t$ as
$$
\overline W(t,u) \sim \sum_{|\alpha|>0} t^\alpha \bar Z_\alpha(u),
$$
where the $\bar Z_\alpha$ are $C^\infty$ vector fields. In fact we have

\begin{lemma}\label{AA}
Associated to each $(X_\alpha, d_\alpha)\in \mathcal{P}\cup \mathcal{N}$, $Z_\alpha =\bar Z_\alpha$ on $B^{\bar n}(\eta_1)$. 
\end{lemma}

\begin{proof}
Fix $u\in B^{\bar n}(\eta_1)$. For every $M \in \mathbb{N}_{>0}$,
\begin{align*}
W^j(t, u) &=D\Phi_j(u)^{-1}(W(\tilde \tau_j^{-1}\delta_j t, \Phi_j(u)))\\
&= D\Phi_j(u)^{-1}\big(\sum_{0<|\alpha|<M}t^\alpha \tau_j^{-\theta(d_\alpha)} \delta_j^{d_\alpha} X_\alpha (\Phi_j(u)) + O(|t|^M)\big) = \sum_{0<|\alpha|<M} t^\alpha Z_\alpha^{j}(u) + O(|t|^M),
\end{align*}
where the last equality is because $D\Phi_j(u)^{-1}$ is a linear map, which maps $O(|t|^M)$ boundedly to $O(|t|^M)$. Thus the Taylor series for $W^j$ in $t$ is
$$
W^j(t,u) \sim \sum_{|\alpha|>0} t^\alpha Z_\alpha^j(u).
$$
Since $W^j \to \overline{W}$ in $C^\infty(B^N(\rho)\times B^{\bar n}(\eta_1))$, $Z_\alpha^{j} \to \bar Z_\alpha$ in $C^\infty(B^{\bar n}(\eta_1))$. Hence $Z_\alpha = \bar Z_\alpha$ on $B^{\bar n}(\eta_1)$. 
\end{proof}

The pure and nonpure power sets of $\overline W(t,u)$ are
\begin{align*}
&\widetilde{\mathcal{P}}=\{( Z_\alpha,d_\alpha): (X_\alpha,d_\alpha)\in \mathcal{P}\},\\
&\widetilde{\mathcal{N}}=\{( Z_\alpha,d_\alpha): (X_\alpha, d_\alpha)\in \mathcal{N}\}.
\end{align*}
Let 
\begin{align*}
&\widetilde{\mathcal{F}}=\{( Z_\alpha, d_\alpha): (X_\alpha, d_\alpha)\in \mathcal{F}\},\\
&\widetilde{\mathcal{F}}'=\{( Z_\alpha,d_\alpha): (X_\alpha,d_\alpha)\in \mathcal{F}'\}.
\end{align*}
Then
\begin{align*}
& \mathcal{L}(\widetilde{\mathcal{P}}) = \{( Z_\alpha, d_\alpha): (X_\alpha, d_\alpha)\in \mathcal{L}(\mathcal{P})\},\\
& \mathcal{L}(\widetilde{\mathcal{P}}\cup \widetilde{\mathcal{N}}) = \{( Z_\alpha, d_\alpha): (X_\alpha, d_\alpha)\in \mathcal{L}(\mathcal{P}\cup \mathcal{N})\},\\
&\widetilde{\mathcal{F}} \subseteq \mathcal{L}(\widetilde{\mathcal{P}}),\quad \widetilde{\mathcal{F}} \subseteq \widetilde{\mathcal{F}}' \subseteq \mathcal{L}(\widetilde{\mathcal{P}} \cup \widetilde{\mathcal{N}}).  
\end{align*}
Note $Z^j_l = \tau_j^{-\theta(d_l) +\theta_l}\ Y^j_l$. Then $Z_l^j \to 0$ in $C^\infty(B^{\bar n}(\eta_1))$ for each $(X_l, d_l) \in \mathcal{F}$. Thus every vector field in $\widetilde{\mathcal{F}}$ vanishes on $B^{\bar n}(\eta_1)$. Denote $Z_l:= \lim_j Z_l^j$. Then we can write $\widetilde{\mathcal{F}}'=\big\{\big(Z_l, d_l\big): 1\leq l \leq Q\big\}$.

Suppose for contradiction that there exists a neighborhood $V \subseteq B^{\bar n}(\eta_1)$ of $0$, such that $Y_{p+1}^j, \ldots, Y_q^j$ all converges to $0$ in $C^\infty(V)$ as $j\to \infty$. For every $1\leq l \leq Q$, $\{Y_l^j\}_j$ is bounded in $C^\infty(B^{\bar n}(\eta_1))$. As in (\ref{prairie}), $(X_{p_0}, d_{p_0}) \in \mathcal{F}' \backslash \mathcal{F}$. So $X_{p_0}$ is an iterated commutator of elements in $\{X_1, \ldots, X_q\}$, and must involve at least one element in $\{X_{p+1}, \ldots, X_q\}$. Thus $Y_{p_0}^j = \Phi_j^*(\tau_j^{-\theta_{p_0}}\delta_j^{d_{p_0}} X_{p_0})$ is an iterated commutator of elements in $\{Y_1^j, \ldots, Y_q^j\}$, and must involve at least one element in $\{Y_{p+1}^j, \ldots, Y_q^j\}$. Therefore $Y_{p_0}^j$ converges to $0$ in $C^\infty(V)$ as $j\to \infty$; contradicting to the fact that $|Y_{p_0}^j| \approx 1$ on $B^{\bar n}(\eta_1)$. Hence for every neighborhood $V \subseteq B^{\bar n}(\eta_1)$ of $0$, there exists $p< l \leq q$, such that $Z_{l}^j = Y_{l}^j$ does not converge to $0$ in $C^\infty(V)$, and thus $Z_{l}$ is nonzero somewhere in $V$. Therefore there exists a sequence $u_k \to 0$ in $B^{\bar n}(\eta_1)$ as $k\to \infty$, such that one of $Z_{p+1}(u_k), \ldots, Z_q(u_k)$ is nonzero. By taking a subsequence of $\{u_k\}_k$, still denoted as $\{u_k\}_k$, there exists $p<l_0\leq q$, such that 
\begin{equation}\label{seminar}
Z_{l_0}(u_k) \neq 0 \text{ for all } k.
\end{equation}

\begin{remark}
How we find $l_0$ here corresponds to (but is much more complicated than) how we find $\alpha_0$ in Section \ref{sntz} using lines $\pi_\beta$ of slope $-\frac{b_1}{b_2}$. 
\end{remark}

\begin{proposition}\label{juemi}
$\widetilde{\mathcal{F}}'$ strongly controls $\overline{W}(t,u)$ on $B^N(\rho)\times B^{\bar n}(\eta_1)$ and every element of $\mathcal{L}(\widetilde{\mathcal{P}} \cup \widetilde{\mathcal{N}})$ on $B^{\bar n}(\eta_1)$. $\widetilde{\mathcal{F}}$ strongly controls every element of $\mathcal{L}(\widetilde{\mathcal{P}})$ on $B^{\bar n}(\eta_1)$, and hence every vector field in $\mathcal{L}(\widetilde{\mathcal{P}})$ vanishes on $B^{\bar n}(\eta_1)$.
\end{proposition}

\begin{proof}
In the proof of Proposition \ref{ppp}, associated to each $(Z_\alpha, d_\alpha)\in \mathcal{L}(\widetilde{\mathcal{P}} \cup \widetilde{\mathcal{N}})$, we have that for $j\geq N_\alpha$, on $B^{\bar n}(\eta_1)$,
\begin{equation}\label{jiayucun}
Z_\alpha^j = \sum_{d_l \leq d_\alpha} \tau_j^{-\theta(d_\alpha) + \theta(d_l)}\, \delta_j^{d_\alpha -d_l} \,(c_{\alpha, l} \circ \Phi_j) \,Z_l^j,
\end{equation}
and each $\{\tau_j^{-\theta(d_\alpha) + \theta(d_l)} \delta_j^{d_\alpha -d_l} c_{\alpha, l} \circ \Phi_j\}_{j\geq N_\alpha}$ is bounded in $C^\infty(B^{\bar n}(\eta_1))$. Then there exists a subsequence $\{(x_{j_k}, \delta_{j_k}, \tau_{j_k})\}_k$ such that $\tau_{j_k}^{-\theta(d_\alpha) + \theta(d_l)} \delta_{j_k}^{d_\alpha -d_l} c_{\alpha, l} \circ \Phi_{j_k}$ converges in $C^\infty(B^{\bar n}(\eta_1))$ to some $b_{\alpha, l}$. Take the limit in $C^\infty(B^{\bar n}(\eta_1))$ of this subsequence of (\ref{jiayucun}), then we have
$$
Z_\alpha = \sum_{d_l \leq d_\alpha} b_{\alpha, l} \,Z_l, \text{~on~} B^{\bar n}(\eta_1).
$$
Hence $\widetilde{\mathcal{F}}'$ strongly controls each $(Z_\alpha, d_\alpha) \in \mathcal{L}(\widetilde{\mathcal{P}}\cup \widetilde{\mathcal{N}})$ on $B^{\bar n}(\eta_1)$. Similarly, $\widetilde{\mathcal{F}}$ strongly controls each element of $\mathcal{L}(\widetilde{\mathcal{P}})$ on $B^{\bar n}(\eta_1)$. And since every vector field in $\widetilde{\mathcal{F}}$ vanishes on $B^{\bar n}(\eta_1)$, every vector field in $\mathcal{L}(\widetilde{\mathcal{P}})$ vanishes on $B^{\bar n}(\eta_1)$.

In the proof of Proposition \ref{ppp}, on $B^N(\rho)\times B^{\bar n}(\eta_1)$,
\begin{equation}\label{jiliujin}
W^j(t,u)= \sum_{\substack{ 1\leq l \leq Q\\ \text{deg}\,(\alpha) = d_l}} t^{\alpha} \,c_{l,\alpha} \big(\tilde \tau_j^{-1} \delta_j t, \Phi_j(u)\big) \,Z_l^j(u), 
\end{equation}
and each $\{c_{l,\alpha}(\tilde \tau_j^{-1} \delta_j t, \Phi_j(u))\}_{j}$ is bounded in $C^\infty(B^N(\rho)\times B^{\bar n}(\eta_1))$; thus there exists a subsequence $\{(x_{j_i}, \delta_{j_i}, \tau_{j_i})\}_i$ such that $c_{l,\alpha} (\tilde \tau_{j_i}^{-1} \delta_{j_i} t, \Phi_{j_i}(u))$ converges in $C^\infty(B^N(\rho)\times B^{\bar n}(\eta_1))$ to some $b_{l, \alpha}$. Take the limit in $C^\infty(B^N(\rho)\times B^{\bar n}(\eta_1))$ of this subsequence of (\ref{jiliujin}), then we have
\begin{equation}\label{INDEP}
\overline{W}(t,u) = \sum_{\substack{ 1\leq l \leq Q\\ \text{deg}\,(\alpha) = d_l}} t^{\alpha} b_{l, \alpha}(t,u)\, Z_l(u), \text{~on~} B^N(\rho)\times B^{\bar n}(\eta_1).
\end{equation}
Hence $\widetilde{\mathcal{F}}'$ strongly controls $\overline{W}(t,u)$ on $B^N(\rho)\times B^{\bar n}(\eta_1)$.
\end{proof}

By Proposition \ref{ppp}, for each $i$, either $\delta_j^{e_i} \neq 0$ for every $j$, or $\delta_j^{e_i} =0$ for every $j$. Let
\begin{align*}
&B= \{1\leq i \leq N: \delta_j^{e_i}\neq 0 \text{~for all~} j\}, \quad B^C= \{1\leq i \leq N: \delta_j^{e_i} =0 \text{~for all~} j\},\quad |B|=\bar N,\\
& A = \{1\leq \mu \leq \nu: e_i^\mu \neq 0 \text{~for some~} i \in B\}, \quad A^C = \{1\leq \mu \leq \nu: e_i^\mu = 0 \text{~for every~} i \in B\}, \quad |A|=\bar \nu.
\end{align*}
Then for $\mu \in A$, $\delta_j^\mu \neq 0$ for all $j$.
Let $t_B$ be the vector consisting of the coordinates $t_i$ of $t$ with $i\in B$, and let $t_{B^C}$ be the vector consisting of the coordinates $t_i$ of $t$ with $i\in B^C$. Without loss of generality, we can assume $B= \{1, \ldots, \bar N\}$, and $A=\{1, 2, \ldots, \bar \nu\}$. Let 
$$
\bar e= \{e_i^\mu: 1\leq i \leq \bar N, 1\leq \mu \leq \bar \nu\}, \quad \bar e_i = (e_i^1, \ldots, e_i^{\bar \nu}) \text{ for } 1\leq i \leq \bar N.
$$
Note $\bar N, \bar \nu, \bar e$ satisfy (3). The multi-parameter dilations on $\mathbb{R}^{\bar N}$ using $\bar e$ are defined by $\delta t_B = (\delta^{\bar e_1} t_1, \ldots, \delta^{\bar e_{\bar N}} t_{\bar N})$ for $\delta \in [0,\infty)^{\bar \nu}$. For a function $\varsigma(t_B)$ on $\mathbb{R}^{\bar N}$, and for $\delta\in (0,\infty)^{\bar \nu}$, set $\varsigma^{(\delta)}(t_B) := \delta^{\bar e_1 + \cdots + \bar e_{\bar N}} \varsigma(\delta t_B)$.

For $1\leq l \leq Q$, if there exists $\alpha$ with $\text{deg}\,\alpha = d_l$ satisfying $\alpha_i \neq 0$ for some $i\in B^C$, then $\delta_j^{d_l} =0$ for all $j$, and thus the corresponding $Z_l$ and $Z_l^j=\Phi_j^*(\tau_j^{-\theta(d_l)} \delta_j^{d_l} X_l)$ vanish on $B^{\bar n}(\eta_1)$. Similarly, for $(Z, d) \in \mathcal{L}(\widetilde{\mathcal{P}} \cup \widetilde{\mathcal{N}})$ with $d^\mu \neq 0$ for some $\mu \in A^C$, since $\delta_j^{d} =0$ for all $j$, $Z$ vanishes on $B^{\bar n}(\eta_1)$. Note $\tilde \tau_j^{-1} \delta_j t= ( \tau_j^{-\theta(e_1)} \delta_j^{e_1}t_1, \ldots, \tau_j^{-\theta(e_N)} \delta_j^{e_N}t_N)$, and $\delta_j^{e_i} =0$ for $i\in B^C$. So the functions $c_{l,\alpha} (\tilde \tau_j^{-1} \delta_j t, \Phi_j(u))$ in (\ref{jiliujin}) and $b_{l, \alpha}(t,u)$ in (\ref{INDEP}) are independent of $t_{B^C}$. Therefore $W^j(t,u)$ and $\overline{W}(t,u)$ are independent of $t_{B^C}$. So we can define $C^\infty$ functions on $B^{\bar N}(\rho) \times B^{\bar n}(\eta_1)$ as $W^j(t_B, u):= W^j(t,u)$, $\overline{W}(t_B,u): = \overline{W}(t,u)$. As in (\ref{seminar}), $Z_{l_0}(u_k) \neq 0$ for all $k$. So $\delta_j^{d_{l_0}} \neq 0$ for all $j$, and thus $\bar N>0$, and $d_{l_0}^\mu = 0$ for all $\mu \in A^C$. And since $d_{l_0}$ is nonpure, $\bar \nu >1$.

For $d = (d^1, \ldots, d^{\nu}) \in \mathbb{R}^\nu$, denote $\bar d:= (d^1, \ldots, d^{\bar \nu})$. Using $\bar e$, the pure and nonpure power sets of $\overline W(t_B,u)$ are
\begin{align*}
\overline{\mathcal{P}}=\{(Z,\bar d): (Z,d)\in \widetilde{\mathcal{P}}, d^\mu = 0 \text{~for~} \mu \in A^C\},\quad \overline{\mathcal{N}}=\{(Z,\bar d): (Z,d)\in \widetilde{\mathcal{N}}, d^\mu = 0 \text{~for~} \mu \in A^C\}.
\end{align*}
Let 
\begin{align*}
&\overline{\mathcal{F}}=\{(Z, \bar d): (Z, d) \in \widetilde{\mathcal{F}}, d^\mu = 0 \text{~for~} \mu \in A^C\},\quad \overline{\mathcal{F}}'=\{(Z, \bar d): (Z, d) \in \widetilde{\mathcal{F}}', d^\mu = 0 \text{~for~} \mu \in A^C\}.
\end{align*}
Then
\begin{align*}
& \mathcal{L}(\overline{\mathcal{P}}) = \{(Z, \bar d): (Z, d)\in \mathcal{L}(\widetilde{\mathcal{P}}), d^\mu =0 \text{~for~} \mu \in A^C\},\\
& \mathcal{L}(\overline{\mathcal{P}}\cup \overline{\mathcal{N}}) = \{(Z, \bar d): (Z, d)\in \mathcal{L}(\widetilde{\mathcal{P}} \cup \widetilde{\mathcal{N}}), d^\mu =0 \text{~for~} \mu \in A^C\},\\
&\overline{\mathcal{F}} \subseteq \mathcal{L}(\overline{\mathcal{P}}),\quad \overline{\mathcal{F}} \subseteq \overline{\mathcal{F}}' \subseteq \mathcal{L}(\overline{\mathcal{P}} \cup \overline{\mathcal{N}}).
\end{align*}

Now we have $(Z_{l_0}, \bar d_{l_0}) \in \overline{\mathcal{F}}'$, and $Z_{l_0}(u_k) \neq 0$ for a sequence $u_k \to 0$. Every vector field in $\mathcal{L}(\overline{\mathcal{P}})$ vanishes on $B^{\bar n}(\eta_1)$. $\overline{\mathcal{F}}'$ strongly controls $\overline{W}(t_B,u)$ on $B^{\bar N}(\rho)\times B^{\bar n}(\eta_1)$, and strongly controls every element of $\mathcal{L}(\overline{\mathcal{P}} \cup \overline{\mathcal{N}})$ on $B^{\bar n}(\eta_1)$.

By Lemma \ref{12.1}, $\overline W(t_B,u)$ corresponds to a $C^\infty$ function $\bar \gamma_{t_B}(u)$ defined on a neighborhood of $(t_B,u)=(0,0)$ in $\mathbb{R}^{\bar N}\times \mathbb{R}^{\bar n}$ with $\bar \gamma_0(u)\equiv u$, and for each $j$, $W^j(t_B,u)$ corresponds to a $C^\infty$ function $\gamma^j_{t_B}(u)$ defined on a neighborhood of $(t_B,u)=(0,0)$ in $\mathbb{R}^{\bar N}\times \mathbb{R}^{\bar n}$ with $\gamma^j_0(u)\equiv u$. $\bar \gamma_{t_B}(u)$ is actually the desired function; the pure and nonpure power sets of $\overline W(t_B, u)$ have satisfied all the relevant conditions in Proposition \ref{1}, and it remains to show the operators corresponding to $\bar \gamma_{t_B}(u)$ are bounded.

\subsection{Boundedness of the corresponding operators}\label{natural3}
Let $\bar \gamma_{t_B}(u)$ and $\gamma^j_{t_B}(u)$ be as above, which correspond to $\overline{W}(t_B,u)$ and $W^j(t_B, u)$, respectively. Note $W^j(t_B,u) \to \overline{W}(t_B,u)$ in $C^\infty(B^{\bar N}(\rho) \times B^{\bar n}(\eta_1))$. $\{(x_j, \delta_j, \tau_j)\}$ is the sequence obtained in Proposition \ref{ppp}. $\Phi_j(u) = e^{u\cdot (\tau_j^{-1} \delta_j \mathbf{X})_J} x_j$ is given in Subsection \ref{natural2}, and satisfies $\Phi_j(B^{\bar n}(\eta_1)) \subseteq B_{(\tau_j^{-1} \delta_j \mathbf{X}, \mathbf{d})}(x_j, \Vec{\xi_1})$. Recall $\tilde \tau_j^{-1} \delta_j t = (\tau_j^{-\theta(e_1)} \delta_j^{e_1} t_1, \ldots, \tau_j^{-\theta(e_N)} \delta_j^{e_N} t_N)$, and recall from Subsection \ref{natural1} that $W(t,x)$ is strongly controlled by $(\mathbf{X}, \mathbf{d})$ on $B^N(\rho)\times \Omega$, and that $\gamma_t^{-1}(x)$ is defined and smooth on $B^N(\rho)\times U \subseteq B^N(\rho)\times \Omega$.
It remains to show
\begin{proposition}\label{UU}
The operator
\begin{equation}\label{ops}
\overline Tf(u) = \bar \psi(u)\int f(\bar{\gamma}_{t_B}(u))K(t_B)\,dt_B
\end{equation}
is bounded on $L^p(\mathbb{R}^{\bar n})$ for every $\bar \psi \in C_c^\infty(\mathbb{R}^{\bar n})$ with sufficiently small support, and every $K\in \mathcal{K}(\bar N, \bar e, \bar a, \bar \nu)$ with sufficiently small $\bar a>0$.
\end{proposition}

To prove this proposition, we need to realize the operator in (\ref{ops}) as the limit, in some sense, of the operators corresponding to $\gamma^j_{t_B}(u)$ as $j \to \infty$, and show the uniform boundedness of these operators. We first show $\Phi_j^{-1} \circ \gamma_{\tilde \tau_j^{-1} \delta_j {t}} \circ \Phi_j(u)\big|_{t_{B^C}=0}$ and $\gamma_{t_B}^j(u)$ are well-defined and equal to each other on a common neighborhood of $(t_B,u) =(0,0)$ in $\mathbb{R}^{\bar N}\times \mathbb{R}^{\bar n}$ for every $j\gg1$.

\begin{lemma}\label{VV'}
For every $2$-admissible constant $0<\eta_2 \leq \eta_1$, there exists $0<\rho_0 \leq \rho$, $N_0\in \mathbb{N}_{>0}$, and a $2$-admissible constant $0<\eta_0 < \eta_2$, such that for $j \geq N_0$, $u\in B^{\bar n}(\eta_0), t\in B^N(\rho_0)$,
$$
\gamma_{\tilde \tau_j^{-1} \delta_j t} \circ \Phi_j(u) \in \Phi_j (B^{\bar n}(\eta_2)), \quad \gamma_{\tilde \tau_j^{-1} \delta_j t}^{-1} \circ \Phi_j(u) \in \Phi_j (B^{\bar n}(\eta_2)).
$$
\end{lemma}

\begin{proof}
By Theorem \ref{16}, there exists a $2$-admissible constant $0<\xi_3<\xi_1$ such that $B_{(\tau_j^{-1} \delta_j \mathbf{X}, \mathbf{d})}(x_j, \Vec{\xi_3}) \subseteq \Phi_j(B^{\bar n}(\eta_2))$ for each $j$. Thus it suffices to show there exists $0<\rho_0 \leq \rho$, $N_0 \in \mathbb{N}_{>0}$, and a $2$-admissible constant $0<\xi_0<\xi_3$ such that for $j \geq N_0, y\in B_{(\tau_j^{-1} \delta_j \mathbf{X},\mathbf{d})}(x_j, \Vec{\xi_0}), t\in B^N(\rho_0)$,
$$
\gamma_{{\tilde \tau_j}^{-1} \delta_j t} (y) \in B_{(\tau_j^{-1} \delta_j \mathbf{X},\mathbf{d})}(x_j, \Vec{\xi_3}), \quad \gamma_{{\tilde \tau_j}^{-1} \delta_j t}^{-1} (y) \in B_{(\tau_j^{-1} \delta_j \mathbf{X},\mathbf{d})}(x_j, \Vec{\xi_3}).
$$
If the above holds, then by Theorem \ref{16}, there exists a $2$-admissible constant $0<\eta_0 <\eta_2$ such that $\Phi_j(B^{\bar n}(\eta_0)) \subseteq B_{(\tau_j^{-1} \delta_j \mathbf{X},\mathbf{d})}(x_j, \Vec{\xi_0})$ for each $j$, and thus for $j\geq N_0$, $u\in B^{\bar n}(\eta_0)$, $t\in B^N(\rho_0)$,
$$
\gamma_{\tilde \tau_j^{-1} \delta_j t} \circ \Phi_j(u), \gamma_{\tilde \tau_j^{-1} \delta_j t}^{-1} \circ \Phi_j(u) \in B_{(\tau_j^{-1} \delta_j \mathbf{X},\mathbf{d})}(x_j, \Vec{\xi_3}) \subseteq  \Phi_j(B^{\bar n}(\eta_2)).
$$

Therefore it suffices to show for any $\xi_4>0$, there exists $0<\rho_0\leq \rho$ and a neighborhood $U' \subseteq U$ of $0$ in $\mathbb{R}^n$ such that for $\delta \in [0,1]^\nu$, $t\in B^N(\rho_0), y\in U'$,
$$
\gamma_{\delta t}(y) \in B_{(\mathbf{X},\mathbf{d})}(y,\xi_4 \delta), \quad \gamma_{\delta t}^{-1}(y) \in B_{(\mathbf{X},\mathbf{d})}(y, \xi_4 \delta).
$$
Indeed for $\xi_4 = \frac{\xi_3}{2^{\frac{1}{\min_l |d_l|_1}}}$ and for the corresponding $U'$, by taking $\xi_0 = \frac{\xi_3}{2^{\frac{1}{\min_l |d_l|_1}}}$, there exists $N_0 \in \mathbb{N}_{>0}$ such that for $j\geq N_0$, $B_{(\tau_j^{-1} \delta_j \mathbf{X},\mathbf{d})}(x_j, \Vec{\xi_0}) \subseteq U'$, and thus for $j\geq N_0$, by taking $\delta = \tilde \tau_j^{-1} \delta_j$, for $y\in B_{(\tau_j^{-1} \delta_j \mathbf{X},\mathbf{d})}(x_j, \Vec{\xi_0}) \subseteq U'$, $t\in B^N(\rho_0)$,
$$
\gamma_{\tilde \tau_j^{-1} \delta_j t}(y), \gamma_{\tilde \tau_j^{-1} \delta_j t}^{-1}(y) \in B_{(\mathbf{X},\mathbf{d})}(y, \xi_4 \tilde \tau_j^{-1} \delta_j) \subseteq B_{(\tau_j^{-1} \delta_j \mathbf{X},\mathbf{d})}(y, \Vec{\xi_4}) \subseteq B_{(\tau_j^{-1} \delta_j \mathbf{X},\mathbf{d})}(x_j, \Vec{\xi_3}).
$$

Let $0<\rho_0'\leq \min\{1,\rho\}$ and a neighborhood $U''\subseteq U$ of $0$ be such that for every $\delta\in [0,1]^\nu, t\in B^N(\rho_0')$, $\gamma_{\delta t}(\cdot)$ is a diffeomorphism from $U''$ onto $\gamma_{\delta t}(U'') \subseteq U$. Fix arbitrary $y\in U'', \delta \in [0,1]^\nu, t \in B^N(\rho_0')$. Consider the smooth curve $g: [0,1] \to \mathbb{R}^n, \sigma \mapsto \gamma_{\sigma \delta t}(y)$, where $\sigma \delta t = (\sigma \delta^{e_1}t_1, \cdots, \sigma \delta^{e_N} t_N)$. Then $g(0)= y, g(1) = \gamma_{\delta t}(y)$, and for every $\sigma \in [0,1]$, $\gamma_{\sigma \delta t}(\cdot)$ is a diffeomorphism from $U''$ onto $\gamma_{\sigma \delta t}(U'') \subseteq U$. We have
\begin{equation}\label{curve}
\begin{aligned}
g'(\sigma) & = \frac{1}{\sigma} \frac{d}{d\epsilon}\Big|_{\epsilon=1}  g(\epsilon \sigma ) = \frac{1}{\sigma} \frac{d}{d\epsilon}\Big|_{\epsilon =1} \gamma_{\epsilon \sigma \delta t} (y)\\
& = \frac{1}{\sigma} \frac{d}{d\epsilon}\Big|_{\epsilon =1} \gamma_{\epsilon \sigma \delta t} \circ \gamma_{\sigma \delta t}^{-1} \big( \gamma_{\sigma \delta t}(y)\big) = \frac{1}{\sigma} W(\sigma \delta t, \gamma_{\sigma \delta t}(y))\\
& = \sum_{\substack{ 1\leq l \leq Q\\ \text{deg}\,(\alpha) = d_l}} t^{\alpha} \big(\frac{1}{\sigma} c_{l, \alpha}(\sigma \delta t, g(\sigma)) \sigma^{|\alpha|}\big) \delta^{d_l}  X_l(g(\sigma)).
\end{aligned}
\end{equation}
Note $c_{l, \alpha}(\sigma \delta t, g(\sigma)) \sigma^{|\alpha|-1}$ is bounded uniformly in $y\in U'', \delta \in [0,1]^\nu, t \in B^N(\rho_0'), \sigma \in [0,1]$. Let $\rho_0'$ in addition satisfy $\rho_0' < \frac{(\min\{1, \xi_4\})^{\max_l |d_l|_1}}{\sum_{l,\alpha} \|c_{l,\alpha}\|_{C(B^N(\rho)\times \Omega)}}$. Then $\gamma_{\delta t}(y) = g(1) \in B_{(\mathbf{X},\mathbf{d})} (y, \xi_4 \delta)$.

We now show there exists a neighborhood $U' \subseteq U''$ of $0$ and $0<\rho_0 \leq \rho_0'$, such that for every $\delta \in [0,1]^\nu, t\in B^N(\rho_0), y\in U'$, $\gamma_{\delta t}^{-1}(y) \in B_{(\mathbf{X},\mathbf{d})}(y, \xi_4 \delta)$. There exists $0\in U'\subseteq U''$ and $0<\rho_0 \leq \rho_0'$, such that for every $\delta \in [0,1]^\nu, t\in B^N(\rho_0)$, $\gamma_{\delta t}^{-1} (\cdot)$ is a diffeomorphism from $U'$ onto $\gamma_{\delta t}^{-1}(U') \subseteq U''$. Fix arbitrary $y\in U', \delta \in [0,1]^\nu, t \in B^N(\rho_0)$. Consider the smooth curve $f: [0,1] \to \mathbb{R}^n, \sigma \mapsto \gamma_{(1-\sigma) \delta t} \circ \gamma_{ \delta t}^{-1}(y) \in \gamma_{(1-\sigma)\delta t}(U'') \subseteq U$. Then $f(0)= y, f(1) = \gamma_{\delta t}^{-1}(y)$. Then by (\ref{curve}),
$$
-f'(1-\sigma)=  \sum_{\substack{ 1\leq l \leq Q\\ \text{deg}\,(\alpha) = d_l}} t^{\alpha} \big(\frac{1}{\sigma} c_{l, \alpha}(\sigma \delta t, f(1-\sigma)) \sigma^{|\alpha|}\big) \delta^{d_l}  X_l(f(1-\sigma)),
$$
where $c_{l, \alpha}(\sigma \delta t, f(1-\sigma)) \sigma^{|\alpha|-1}$ is bounded uniformly in $y\in U', \delta \in [0,1]^\nu, t\in B^N(\rho_0), \sigma \in [0,1]$. Let $\rho_0$ in addition satisfy $\rho_0 < \frac{(\min\{1, \xi_4\})^{\max_l |d_l|_1}}{\sum_{l,\alpha} \|c_{l,\alpha}\|_{C(B^N(\rho)\times \Omega)}}$. Then $\gamma_{\delta t}^{-1}(y) = f(1) \in B_{(\mathbf{X},\mathbf{d})} (y, \xi_4 \delta)$.
\end{proof}

So there exists $0<\rho_0 \leq \rho$, $N_0\in \mathbb{N}_{>0}$, and a $2$-admissible constant $0<\eta_0 < \frac{1}{2}\eta_1$, such that for $j \geq N_0, u\in B^{\bar n}(\eta_0), t\in B^N(\rho_0)$, we have $\tilde \tau_j^{-1} \delta_j \in [0,\frac{1}{2}]^\nu$,
$\Phi_j^{-1} \circ \gamma_{\tilde \tau_j^{-1} \delta_j t} \circ \Phi_j(u) \in B^{\bar n}(\frac{1}{2}\eta_1)$, and the functions in the following equation are well-defined:
\begin{align*}
&\quad \frac{d}{d\epsilon}\Big|_{\epsilon=1} \Phi_j^{-1} \circ \gamma_{\epsilon \tilde \tau_j^{-1} \delta_j  t} \circ \Phi_j \circ \Phi_j^{-1} \circ \gamma_{\tilde \tau_j^{-1} \delta_j t}^{-1} \circ \Phi_j (u)\Big|_{t_{B^C}=0}\\
&= D\Phi_j(u)^{-1} \frac{d}{d\epsilon}\Big|_{\epsilon =1}  \gamma_{\epsilon \tilde \tau_j^{-1} \delta_j  t} \circ \gamma_{\tilde \tau_j^{-1} \delta_j t}^{-1}\big( \Phi_j(u)\big)\Big|_{t_{B^C}=0}\\
&= W^j(t,u)\Big|_{t_{B^C}=0} = W^j(t_B,u).
\end{align*}
Thus by Lemma \ref{12.1}, we have:

\begin{corollary}\label{VV}
For $j \geq N_0, u \in B^{\bar n}(\eta_0), t \in B^{N}(\rho_0)$, 
$$
\gamma_{\tilde \tau_j^{-1} \delta_j {t}} \circ \Phi_j(u) \in \Phi_j(B^{\bar n}(\frac{1}{2}\eta_1)), \quad \gamma_{t_B}^j(u)= \Phi_j^{-1} \circ \gamma_{\tilde \tau_j^{-1} \delta_j t} \circ \Phi_j(u) \big|_{t_{B^C}=0}.
$$
\end{corollary}

By Lemma \ref{12.1}, by shrinking $\rho_0$ if necessary \footnote{It is easy to see, via the contraction mapping theorem, using the fact that $W^j(0, \cdot) \equiv 0$ and $\{W^j: j\geq N_0\}$ is bounded in $C^\infty(B^{\bar N}(\rho) \times B^{\bar n}(\eta_1))$, that a common $\rho_0>0$ can be taken for all $j\geq N_0$.}, for $j\geq N_0$, we have $\gamma_{t_B}^j(u) = w^j(1,t_B,u)$ on $B^{\bar N}(\rho_0) \times B^{\bar n}(\eta_0)$, where $w^j(\epsilon, t_B,u)$ is the unique solution on $[0,1] \times B^{\bar N}(\rho_0) \times B^{\bar n}(\eta_0)$ to the ODE
$$
\frac{d}{d\epsilon} w^j(\epsilon, t_B, u) = \frac{1}{\epsilon} W^j(\epsilon t_B, w^j(\epsilon, t_B, u)), \quad w^j(0,t_B,u) = u.
$$
By Lemma \ref{12.1}, by shrinking $\rho_0$ if necessary, we have $\bar \gamma_{t_B}(u) = \bar w(1,t_B,u)$ on $B^{\bar N}(\rho_0) \times B^{\bar n}(\eta_0)$, where $\bar w(\epsilon, t_B, u)$ is the unique solution on $[0,1]\times B^{\bar N}(\rho_0) \times B^{\bar n}(\eta_0)$ to the ODE
$$
\frac{d}{d\epsilon} \bar w(\epsilon, t_B, u) = \frac{1}{\epsilon} \overline{W}(\epsilon t_B, \bar w(\epsilon, t_B, u)), \quad \bar w(0,t_B,u) = u.
$$

The proofs of the following two lemmas are in Appendix \ref{jinzhangzhi}.

\begin{lemma}\label{contraction}
There exists $0<\rho''\leq \rho_0$ such that $\bar w$ and every $w^j (j\geq N_0)$ map from $[0,1]\times B^{\bar N}(\rho'') \times B^{\bar n}(\eta_0)$ to $B^{\bar n}(\frac{1}{2}\eta_1)$.
\end{lemma}

\begin{lemma}\label{wujinzang}
There exists $0<\rho'\leq \rho''$ such that $w^j \to \bar w$ uniformly on $[0,1] \times B^{\bar N}(\rho')\times B^{\bar n}(\eta_0)$.
\end{lemma}

Thus $\bar \gamma_{t_B}(u)$ and every $\gamma_{t_B}^j(u) (j\geq N_0)$ map from $B^{\bar N}(\rho') \times B^{\bar n}(\eta_0)$ to $B^{\bar n}(\frac{1}{2}\eta_1)$, and $\gamma_{t_B}^j(u) \to \bar \gamma_{t_B}(u)$ uniformly on $B^{\bar N}(\rho') \times B^{\bar n}(\eta_0)$.

Next we show the uniform boundedness of the operators corresponding to $\gamma^j_{t_B}(u)$. By the assumptions, the operator $Tf(x) = \psi(x)\int f(\gamma_t(x))K(t)\,dt$ is bounded on $L^p(\mathbb{R}^{n})$ for all $\psi\in C_c^\infty(\mathbb{R}^n)$ with sufficiently small support and all $K\in \mathcal{K}(N,e,a,\nu)$ with sufficiently small $a>0$. Fix a sufficiently small $a>0$ and fix such a $\psi$ with $\psi(0)\neq 0$. Set 
$$
0< \bar a \leq \min \{ \frac{a}{2}, \frac{\rho'}{2} \}.
$$
Then we have the following lemma, whose proof is in Appendix \ref{jinzhangzhi}.

\begin{lemma}\label{yangtianxiang}
If we fix an arbitrary $\{\varsigma_k\}_{k\in \mathbb{N}^{\bar \nu}}\in  \overline{\mathcal{K}}(\bar N, \bar e, \bar a, \bar \nu)$, the operator
$$
T_{j,M} f(x) := \psi(x) \sum_{\substack{k\in \mathbb{N}^{ \bar \nu} \\ |k|\leq M}} \int f(\gamma_{\tilde \tau_j^{-1} \delta_j t}(x)) \delta_0(t_{B^C}) \varsigma_k^{(2^k)}(t_B)\,dt
$$
is bounded on $L^p(\mathbb{R}^n)$ uniformly in $j (\geq N_0), M$, where $\delta_0(t_{B^C})$ is the Dirac function on $\mathbb{R}^{N-\bar N}$ with point mass at $0$.
\end{lemma}

For any $\bar \psi \in C_c^\infty(B^{\bar n}(\eta_0))$, and for any $\{\varsigma_k\}_{k\in \mathbb{N}^{\bar \nu}}\in  \overline{\mathcal{K}}(\bar N, \bar e, \bar a, \bar \nu)$, we set the operators corresponding to $\gamma^j_{t_B}(u)$ as
\begin{equation}\label{coco}
\begin{aligned}
&\quad \bar \psi(u)  \psi(\Phi_j(u)) \sum_{\substack{k\in \mathbb{N}^{\bar \nu}\\ |k|\leq M}} \int f(\gamma^j_{t_B}(u)) \varsigma_k^{(2^k)}(t_B)\,dt_B  \\
&= \bar \psi(u)  \psi(\Phi_j(u)) \sum_{\substack{k\in \mathbb{N}^{\bar \nu}\\ |k|\leq M}} \int f(\Phi_j^{-1} \circ \gamma_{\tilde \tau_j^{-1} \delta_j t} \circ \Phi_j(u)) \delta_0(t_{B^C}) \varsigma_k^{(2^k)}(t_B)\,dt\\
&= \bar \psi \Phi_j^* T_{j,M} (\Phi_j^{-1})^*f(u).
\end{aligned}
\end{equation}
Then we need the following proposition, to obtain uniform boundedness of $\bar \psi \Phi_j^* T_{j,M} (\Phi_j^{-1})^*$.

\begin{proposition}\label{W}
For any $\bar \psi \in C_c^\infty(B^{\bar n}(\eta_0))$, there exists $C>0$ and $N_2 \geq N_0$, such that for $j \geq N_2, M\in \mathbb{N}$, 
$$
\Big\| \bar \psi \Phi_j^* T_{j,M} (\Phi_j^{-1})^*\Big\|_{L^p(\mathbb{R}^{\bar n}) \to L^p(\mathbb{R}^{\bar n})} \leq C \Big\| T_{j,M} \Big\|_{L^p(\mathbb{R}^n) \to L^p(\mathbb{R}^n)}.
$$
\end{proposition}

This inequality is actually trivial when $\bar n = n$ because in (\ref{coco}), $\gamma_{\tilde \tau_j^{-1} \delta_j t}$ maps $\Phi_j(u)$ into the image of $\Phi_j$, and $\Phi_j$ is a local diffeomorphism with ``essentially constant'' Jacobian. But when $\bar n <n$, we need to show that, roughly speaking, the image of a perturbed $\Phi_j(u)$ under the map $\gamma_{\tilde \tau_j^{-1} \delta_j t}$ does not deviate far away from the image of $\Phi_j$. So we need the following two lemmas whose proofs are in Appendix \ref{jinzhangzhi}.

For every $j$, $\Phi_j(B^{\bar n}(\eta_1)) \subseteq \mathbb{R}^n$ is a submanifold. The normal space to $\Phi_j(B^{\bar n}(\eta_1))$ at $x_j$ has an orthonormal basis $\{V_{\bar n+1}, \ldots, V_n\}$. For $v\in \mathbb{R}^{n-\bar n}$, $v \cdot \mathbf{V}$ denotes $\sum_{\bar n +1 \leq i \leq n} v_i V_i$. By shrinking $\eta_1$ if necessary, there exists $\bar \sigma>0$, $N_1\geq N_0$, and $C'>0$, such that for every $j \geq N_1$, the map
\begin{align*}
H_j: B^{\bar n}(\eta_1) \times B^{n-\bar n}(\bar \sigma) \to \Omega \subseteq \mathbb{R}^n, \quad (u,v) \mapsto e^{u\cdot \mathbf{X}_{J}}(x_j + v \cdot \mathbf{V})
\end{align*}
is defined and smooth, and $\big\| H_j \big\|_{C^2(B^{\bar n}(\eta_1) \times B^{n-\bar n}(\bar \sigma))} \leq C'$. Define
\begin{align*}
E_j: B^{\bar n}(\eta_1) \times B^{n-\bar n}(\bar \sigma) \to \Omega \subseteq \mathbb{R}^n, \quad (u,v) \mapsto e^{u\cdot (\tau_j^{-1} \delta_j \mathbf{X})_{J}}(x_j + v \cdot \mathbf{V}).
\end{align*}
Note $E_j(u,0) = \Phi_j(u)$.

\begin{lemma}\label{modu}
There exists $N_2 \geq N_1$, such that for $j \geq N_2$, there exists $0< \sigma_j\leq \bar \sigma$, such that $E_j$ is a diffeomorphism from $B^{\bar n}(\frac{3}{4}\eta_1) \times B^{n-\bar n}(\sigma_j)$ onto its image, and that for $(u,v)\in B^{\bar n}(\frac{3}{4}\eta_1) \times B^{n-\bar n}(\sigma_j)$,
$$
|\det DE_j(u,v)| \approx |\det DE_j(0,0)|>0.
$$
\end{lemma}

\begin{lemma}\label{strada}
There exists $C_0\geq 1$, such that for $j \geq N_2$, there exists $0<\epsilon_j \leq \frac{1}{2 C_0} \sigma_j$, such that for $0<\epsilon \leq \epsilon_j$ and $t\in B^N(\rho')$, $E_j^{-1} \circ \gamma_{\tilde \tau_j^{-1}\delta_j t} \circ E_j$ maps from $B^{\bar n}(\eta_0)$ $\times B^{n-\bar n}(\epsilon)$ into $B^{\bar n}(\frac{3}{4}\eta_1)\times B^{n-\bar n}(C_0\epsilon)$.
\end{lemma}

\begin{proof}[\textbf{Proof of Proposition \ref{W}}]
Fix arbitrary $\bar \psi \in C_c^\infty(B^{\bar n}(\eta_0))$, $j \geq N_2$, and $M\in \mathbb{N}$. Let $0\leq \lambda_j \leq 1$ be a function in $C_c^\infty(B^{n-\bar n}(\sigma_j))$ satisfying $\lambda_j \equiv 1$ on $B^{n-\bar n}(\frac{\sigma_j}{2})$, and let $0\leq \lambda_{\epsilon} \leq 1$ be a function in $C_c^\infty(B^{n-\bar n}(2C_0\epsilon))$ satisfying $\lambda_\epsilon \equiv 1$ on $B^{n-\bar n}(C_0\epsilon)$. Consider an arbitrary $f\in C_c^\infty(\mathbb{R}^{\bar n})$. For $0< \epsilon \leq \epsilon_j$, we can define $F_j, F_{j,\epsilon} \in C_c^\infty(\mathbb{R}^n)$ such that on $E_j(B^{\bar n} (\frac{3}{4}\eta_1) \times B^{n-\bar n}(\sigma_j))$,
$$
F_j = (f \otimes \lambda_j)\circ  (E_j)^{-1}, \quad F_{j, \epsilon}= (f \otimes \lambda_\epsilon)\circ  (E_j)^{-1},
$$
and 
$$
\|F_{j,\epsilon} \|_{L^p(\mathbb{R}^n)} \approx \|F_{j,\epsilon} \|_{L^p(E_j(B^{\bar n}(\frac{3}{4}\eta_1)\times B^{n-\bar n}(2C_0\epsilon)))}.
$$
Then $F_j =  f \circ \Phi_j^{-1}$ on $\Phi_j(B^{\bar n}(\frac{3}{4}\eta_1))$, and 
\begin{equation}\label{strada2}
F_{j, \epsilon}  = F_j, \text{ on } E_j(B^{\bar n}(\frac{3}{4}\eta_1) \times B^{n-\bar n}(C_0\epsilon)).
\end{equation}
By the dominated convergence theorem, $T_{j,M}F_j$ is continuous. Then
\begin{align*}
&\quad \Big\| \bar \psi \Phi_j^* T_{j,M} (\Phi_j^{-1})^* f\Big\|_{L^p(\mathbb{R}^{\bar n})}\\
&= \Big\| \bar \psi(u) \psi(\Phi_j(u)) \sum_{\substack{k\in \mathbb{N}^{\bar \nu}\\ |k|\leq M}} \int f(\Phi_j^{-1}\circ \gamma_{\tilde \tau_j^{-1} \delta_j t} \circ \Phi_j(u)) \delta_0(t_{B^C})  \varsigma_k^{(2^k)}(t_B)\,dt \Big\|_{L^p(B^{\bar n}(\eta_0))}\\
&\leq \|\bar \psi\|_{C(B^{\bar n}(\eta_0))} \Big\| \psi(\Phi_j(u)) \sum_{\substack{k\in \mathbb{N}^{\bar \nu}\\ |k|\leq M}} \int F_j( \gamma_{\tilde \tau_j^{-1} \delta_j t} \circ \Phi_j(u)) \delta_0(t_{B^C}) \varsigma_k^{(2^k)}(t_B)\,dt \Big\|_{L^p(B^{\bar n}(\eta_0))}\\
&= \|\bar \psi\|_{C(B^{\bar n}(\eta_0))} \Big( \int_{B^{\bar n}(\eta_0)} \big| (T_{j,M}F_j)\circ E_j(u,0)\big|^p \,du \Big)^{\frac{1}{p}}\\
&\approx \|\bar \psi \|_{C(B^{\bar n}(\eta_0))} \lim_{\epsilon\to 0}  \Big( \frac{1}{\epsilon^{n-\bar n}} \int_{B^{\bar n}(\eta_0)\times B^{n-\bar n}(\epsilon)} |(T_{j,M} F_j)\circ E_j|^p \,du dv \Big)^{\frac{1}{p}}\\
&= \|\bar \psi \|_{C(B^{\bar n}(\eta_0))} \lim_{\epsilon\to 0}  \Big( \frac{1}{\epsilon^{n-\bar n}} \int_{E_j(B^{\bar n}(\eta_0)\times B^{n-\bar n}(\epsilon))} \frac{|T_{j,M} F_j|^p \,dx}{|\det DE_j(u,v)|} \Big)^{\frac{1}{p}}\\
&\leq \|\bar \psi \|_{C(B^{\bar n}(\eta_0))} \frac{1}{\inf_{B^{\bar n}(\eta_0)\times B^{n-\bar n}(\sigma_j)} |\det DE_j|^{\frac{1}{p}}} \lim_{\epsilon\to 0}  \Big( \frac{1}{\epsilon^{n-\bar n}} \int_{E_j(B^{\bar n}(\eta_0)\times B^{n-\bar n}(\epsilon))} |T_{j,M} F_j|^p \,dx \Big)^{\frac{1}{p}}.
\end{align*}
By Lemma \ref{strada} and (\ref{strada2}), for $j\geq N_2$,
\begin{align*}
&\quad \|\bar \psi \|_{C(B^{\bar n}(\eta_0))} \frac{1}{\inf_{B^{\bar n}(\eta_0)\times B^{n-\bar n}(\sigma_j)} |\det DE_j|^{\frac{1}{p}}} \lim_{\epsilon\to 0}  \Big( \frac{1}{\epsilon^{n-\bar n}} \int_{E_j(B^{\bar n}(\eta_0)\times B^{n-\bar n}(\epsilon))} |T_{j,M} F_j|^p \,dx \Big)^{\frac{1}{p}}\\
&= \|\bar \psi \|_{C(B^{\bar n}(\eta_0))} \frac{1}{\inf_{B^{\bar n}(\eta_0)\times B^{n-\bar n}(\sigma_j)} |\det DE_j|^{\frac{1}{p}}} \lim_{\epsilon\to 0}  \Big( \frac{1}{\epsilon^{n-\bar n}} \int_{E_j(B^{\bar n}(\eta_0)\times B^{n-\bar n}(\epsilon))} |T_{j,M} F_{j, \epsilon}|^p \,dx \Big)^{\frac{1}{p}}\\
&\leq \|\bar \psi\|_{C(B^{\bar n}(\eta_0))}  \frac{1}{\inf_{B^{\bar n}(\eta_0)\times B^{n-\bar n}(\sigma_j)} |\det DE_j |^{\frac{1}{p}}} \|T_{j,M}\|_{L^p(\mathbb{R}^n)\to L^p(\mathbb{R}^n)} \lim_{\epsilon\to 0}  \frac{1}{\epsilon^{(n-\bar n)/p}}  \|F_{j, \epsilon}\|_{L^p(\mathbb{R}^n)}.
\end{align*}
Hence
\begin{align*}
&\quad \Big\| \bar \psi \Phi_j^* T_{j,M} (\Phi_j^{-1})^* f\Big\|_{L^p(\mathbb{R}^{\bar n})}\\
&\lesssim \|\bar \psi\|_{C(B^{\bar n}(\eta_0))}  \frac{1}{\inf_{B^{\bar n}(\eta_0)\times B^{n-\bar n}(\sigma_j)} |\det DE_j|^{\frac{1}{p}}}  \|T_{j,M}\|_{L^p(\mathbb{R}^n)\to L^p(\mathbb{R}^n)} \cdot\\
&  \qquad \qquad \qquad \qquad \qquad \lim_{\epsilon\to 0} \Big( \frac{1}{\epsilon^{n-\bar n}} \int_{B^{\bar n}(\frac{3}{4}\eta_1)\times B^{n-\bar n}(2C_0\epsilon)} |F_{j, \epsilon} \circ E_j|^p |\det DE_j| \,dudv\Big)^{\frac{1}{p}}\\
&\leq \|\bar \psi\|_{C(B^{\bar n}(\eta_0))}  \frac{\sup_{B^{\bar n}(\frac{3}{4}\eta_1)\times B^{n-\bar n}(\sigma_j)} |\det DE_j|^{\frac{1}{p}}}{\inf_{B^{\bar n}(\eta_0)\times B^{n-\bar n}(\sigma_j)} |\det DE_j|^{\frac{1}{p}}}  \|T_{j,M}\|_{L^p(\mathbb{R}^n)\to L^p(\mathbb{R}^n)} \cdot\\
&\quad \qquad \qquad \qquad \qquad \qquad \qquad \lim_{\epsilon\to 0} \Big( \frac{1}{\epsilon^{n-\bar n}} \int_{B^{\bar n}(\frac{3}{4}\eta_1)\times B^{n-\bar n}(2C_0\epsilon)} |F_j\circ E_j|^p  \,dudv\Big)^{\frac{1}{p}}\\
&\approx (2C_0)^{\frac{n-\bar n}{p}} \|\bar \psi\|_{C(B^{\bar n}(\eta_0))} \|T_{j,M}\|_{L^p(\mathbb{R}^n)\to L^p(\mathbb{R}^n)} \Big( \int_{B^{\bar n}(\frac{3}{4}\eta_1)} |F_j\circ \Phi_j|^p\,du \Big)^{\frac{1}{p}}\\
&\leq (2C_0)^{\frac{n-\bar n}{p}} \|\bar \psi\|_{C(B^{\bar n}(\eta_0))} \|T_{j,M}\|_{L^p(\mathbb{R}^n)\to L^p(\mathbb{R}^n)} \|f\|_{L^p(\mathbb{R}^{\bar n})}.
\end{align*}
Therefore $\bar \psi \Phi_j^* T_{j,M} (\Phi_j^{-1})^*$ is bounded on $L^p(\mathbb{R}^{\bar n})$ uniformly in $j(\geq N_2),M$.
\end{proof}

\begin{proof}[\textbf{Proof of Proposition \ref{UU}}]
For every $\bar \psi \in C_c^\infty(B^{\bar n}(\eta_0))$, every $\{\varsigma_k\}_{k\in \mathbb{N}^{\bar \nu}}\in  \overline{\mathcal{K}}(\bar N, \bar e, \bar a, \bar \nu)$, every $M$, and every $f\in C_c^\infty(\mathbb{R}^{\bar n})$, by the dominated convergence theorem,
\begin{align*}
T_Mf(u)&:= \bar \psi(u) \sum_{\substack{k\in \mathbb{N}^{\bar \nu}\\ |k|\leq M}} \int f(\bar \gamma_{t_B}(u)) \varsigma_k^{(2^k)}(t_B)\,dt_B\\
&= \lim_{j\to \infty} \frac{1}{\psi(0)} \bar \psi(u)\psi(\Phi_j(u)) \sum_{\substack{k\in \mathbb{N}^{\bar \nu}\\ |k|\leq M}} \int f(\gamma_{t_B}^j(u)) \varsigma_k^{(2^k)}(t_B)\,dt_B\\
&= \lim_{j\to \infty} \frac{1}{\psi(0)} \bar \psi \Phi_j^* T_{j,M} (\Phi_j^{-1})^* f(u).
\end{align*}
By Fatou's Lemma, 
\begin{align*}
\big\|T_Mf\Big\|_{L^p(\mathbb{R}^{\bar n})} &\leq \frac{1}{\psi(0)} \liminf_{j\to \infty} \big\| \bar \psi \Phi_j^* T_{j,M} (\Phi_j^{-1})^* f\Big\|_{L^p(\mathbb{R}^{\bar n})}\\
&\leq \Big(\frac{1}{\psi(0)} \sup_{j \geq N_2,M} \Big\| \bar \psi \Phi_j^* T_{j,M} (\Phi_j^{-1})^* \Big\|_{L^p\to L^p}\Big) \big\| f\big\|_{L^p(\mathbb{R}^{\bar n})}.
\end{align*}
Thus $T_M$ is bounded on $L^p(\mathbb{R}^{\bar n})$ uniformly in $M$, for every fixed $\bar \psi \in C_c^\infty(B^{\bar n}(\eta_0))$, $\{\varsigma_k\}_{k\in \mathbb{N}^{\bar \nu}}\in  \overline{\mathcal{K}}(\bar N, \bar e, \bar a, \bar \nu)$. By Lemma 16.1 in \cite{L2}, $\sum_{k\in \mathbb{N}^{\bar \nu}} \varsigma_k^{(2^k)}(t_B)$ converges in the sense of distributions with compact support in $\mathbb{R}^{\bar N}$. Thus for $f\in C_c^\infty(\mathbb{R}^{\bar n})$,
$$
\overline T f(u):= \bar \psi(u)\sum_{k\in \mathbb{N}^{\bar \nu}}\int f(\bar \gamma_{t_B}(u))\varsigma_k^{(2^k)}(t_B)\,dt_B 
$$
converges and is equal to $\lim_{M\to \infty} T_M f(u)$.
Then by Fatou's lemma, $\overline T$ is bounded on $L^p(\mathbb{R}^{\bar n})$ for every $\bar \psi \in C_c^\infty(B^{\bar n}(\eta_0))$, $\{\varsigma_k\}_{k\in \mathbb{N}^{\bar \nu}}\in \overline{\mathcal{K}}(\bar N,\bar e,\bar a,\bar \nu)$. Hence $\overline T$ is bounded on $L^p(\mathbb{R}^{\bar n})$ for every $\bar \psi \in C_c^\infty(B^{\bar n}(\eta_0))$, $K\in \mathcal{K}(\bar N,\bar e, \bar a, \bar \nu)$.
\end{proof}

\section{Proof of Proposition \ref{2}}\label{xj}
The assumptions are that there exist $N, \nu>1, e$ satisfying (\ref{e}), $n\in \mathbb{N}_{>0}, p\in (1,\infty)$, and a $C^\infty$ function $\gamma_t(x)$ defined on a neighborhood of $(t,x)=(0,0)$ in $\mathbb{R}^{N}\times \mathbb{R}^{n}$ with $\gamma_0(x)\equiv x$, such that the operator
$$
Tf(x) = \psi(x)\int f(\gamma_t(x))K(t)\,dt
$$
is bounded on $L^p(\mathbb{R}^{n})$ for every $\psi\in C_c^\infty(\mathbb{R}^{n})$ with sufficiently small support, and every $K \in \mathcal{K}(N, e, a, \nu)$ with sufficiently small $a>0$, and such that there exists a finite subset $\mathcal{F}'\subseteq \mathcal{L}(\mathcal{P}\cup \mathcal{N})$ that strongly controls $W$ and every element of $\mathcal{L}(\mathcal{P}\cup \mathcal{N})$, and at least one vector field in $\mathcal{F}'$ is nonzero at a sequence of points convergent to $0$ in $\mathbb{R}^n$, but every vector field in $\mathcal{P}$ vanishes on a common neighborhood of $0$ in $\mathbb{R}^n$. Roughly speaking, the goal is to construct a new $\gamma_t(x)$ whose corresponding operators are still bounded, but the only non-vanishing vector fields are associated with the same nonpure degree.

Without loss of generality, by removing from $\mathcal{F}'$ the elements vanishing on some neighborhood of $0$ in $\mathbb{R}^n$, we can assume every vector field in $\mathcal{F}'$ is nonzero at a sequence of points convergent to $0$. Note every vector field in $\mathcal{P}$ vanishes on a common neighborhood of $0$, and so does every vector field in $\mathcal{L}(\mathcal{P})$. Thus every element of $\mathcal{F}'$ has nonpure degree. Denote $\mathcal{F}'=(\mathbf{X},\mathbf{d})=\{(X_l,d_l)\}_{l=1}^q$. There exists $\rho>0$ and a neighborhood $\Omega$ of $0$ in $\mathbb{R}^n$, such that every vector field in $\mathcal{P}$ vanishes on $\Omega$, $\mathcal{F}'$ strongly controls every element of $\big\{ \big( [X_i, X_j], d_i + d_j \big): 1\leq i, j \leq q \big\}$ on $\Omega$, $\mathcal{F}'$ strongly controls $W$ on $B^N(\rho)\times \Omega$, and every $X_l$ has finite $C^m(\Omega)$ norm for all $m\in \mathbb{N}$.

Since $\mathcal{F}'$ is finite, there exist positive numbers $b_1, \ldots, b_\nu$ such that the linear function $\theta(d)=\theta(d^1, \ldots, d^\nu) := \sum_{\mu=1}^\nu {b_\mu}d^\mu$ separates different degrees $d_l$ in $\mathcal{F}'$. By reordering $(\mathbf{X},\mathbf{d})$, assume $\theta(d_1)\leq \cdots \leq \theta(d_q)$. Let $s= \max\{1\leq l \leq q: \theta(d_l) = \theta(d_1)\}$. Then $d_1=d_2=\cdots =d_s$. (If $s=q$, we already achieve that every element of $\mathcal{F}'$ has the same nonpure degree; if $s<q$, we proceed as follows.) Let $\theta_l=\theta(d_l)$ for $1\leq l \leq s$, and $\theta_{l}=\min\{ 2\theta_1, \frac{\theta_1 + \theta(d_{s+1})}{2}\}$ for $s+1 \leq l \leq q$. Then for $1\leq l \leq s< l' \leq q$, we have $0<\theta_l<\theta_{l'}<\theta(d_{l'})$, $\theta_{l'}\leq 2\theta_l$. For $\tau \geq 1$, denote $\tilde \tau := (\tau^{{b_1}}, \ldots, \tau^{{b_\nu}})$. Then $\tilde \tau^{-1} t = (\tau^{-\theta(e_1)} t_1, \ldots, \tau^{-\theta(e_N)} t_N)$. Let $(\tau^{-1} \mathbf{X},\mathbf{d}) = \{(\tau^{-\theta_l} X_l, d_l)\}_{l=1}^q$, and $(\tilde \tau^{-1} \mathbf{X},\mathbf{d}) = \{ (\tau^{-\theta(d_l)} X_l,d_l)\}_{l=1}^q$.

For $1\leq l, k\leq q$,
$$
[X_l,X_k]=\sum_{\substack{1\leq i\leq q\\ d_i \leq d_l+d_k}} c_{l,k}^i X_i, \text{~on~} \Omega,
$$
where the $c_{l,k}^i$ are $C^\infty$ on a neighborhood of the closure of $\Omega$. Let $\{x_j\}$ be a sequence convergent to $0$ in $\Omega$ such that $X_1(x_j) \neq 0$ for every $j$ (the choice of $X_1$ is unimportant; we can replace $X_1$ with $X_l$ for any $1\leq l \leq s$).
Assume the sequence $\{x_j\}$ is contained in a compact subset of $\Omega$. Then there exists $\xi_1>0$ such that $(\mathbf{X},\mathbf{d})$ satisfies $\mathcal{C}(x_j, \Vec{\xi_1}, \Omega)$ for every $j$. For $\tau\geq 1$, $(\tau^{-1}\mathbf{X},\mathbf{d})$ also satisfies $\mathcal{C}(x_j, \Vec{\xi_1}, \Omega)$ for every $j$. Note $\theta_i \leq \theta_l + \theta_k$ for every $1\leq l,k,i \leq q$. Then on $B_{(\tau^{-1} \mathbf{X},\mathbf{d})}(x_j, \Vec{\xi_1})$,
$$
[\tau^{-\theta_l}X_l, \tau^{-\theta_k} X_k] = \sum_{ d_i\leq d_l+d_k} \tau^{-\theta_l-\theta_k+\theta_i} c_{l,k}^{i} \tau^{-\theta_i}X_i,
$$
and for all $m\in \mathbb{N}$,
$$
\sup_{\tau \geq 1, j} \sum_{|\beta|\leq m} \Big\| (\tau^{-1} \mathbf{X})^\beta \big( \tau^{-\theta_l-\theta_k+\theta_i} c_{l,k}^{i} \big) \Big\|_{C(B_{(\tau^{-1} \mathbf{X},\mathbf{d})}(x_j, \Vec{\xi_1}))} \leq \sum_{|\beta|\leq m} \big\|\mathbf{X}^\beta c_{l,k}^{i}\big\|_{C(\Omega)} <\infty.
$$
Thus $(\tau^{-1}\mathbf{X},\mathbf{d})$ satisfies the conditions of Theorem \ref{16} uniformly in $x_j$ and $\tau (\geq 1)$. In this section, $m$-admissible constants are with respect to the above $(\mathbf{X}, \mathbf{d})$ and $\zeta =1$; see Definition \ref{uniformadmissible}.

We use the notations from Definition \ref{indoor}. Let $n_j= \text{rank}\, \mathbf{X}(x_j)$. By taking a subsequence of $\{x_j\}$, still denoted as $\{x_j\}$, we can assume $n_j$ is constant, denoted as $\bar n$. We have $\text{rank}\,(\tau^{-1} \mathbf{X})(x_j) = \text{rank}\, \mathbf{X}(x_j) =\bar n>0$. Note $X_1(x_j) \neq 0$ for every $j$, and $\tau^{-\theta_1} = \cdots = \tau^{-\theta_s}\gg \tau^{-\theta_{s+1}} = \cdots = \tau^{-\theta_q}$ for $\tau\gg1$. Then for every $j$, we can take a sufficiently large $\tau_j \geq 1$, such that there exists $J_j \in \mathcal{I}(\bar n, q)$ satisfying 
$$
J_j \cap \{1, \ldots, s\} \neq \emptyset, \quad \Big| \det_{\bar n\times \bar n} (\tau_j^{-1} \mathbf{X})_{J_j}(x_j) \Big|_\infty = \Big| \det_{\bar n\times \bar n} (\tau_j^{-1} \mathbf{X})(x_j) \Big|_\infty.
$$
We can let $\tau_j \to \infty$ as $j\to \infty$. By taking a subsequence of $\{(x_j, \tau_j)\}_j$, still denoted as $\{(x_j, \tau_j)\}_j$, we can assume $J_j$ is fixed, denoted as $J$. We have $J \cap \{1, \ldots, s\} \neq \emptyset$.

By Theorem \ref{16}, there exist $2$-admissible constants $\eta_1>0, 0<\xi_2<\xi_1$, such that for every $j$, we can define a diffeomorphism $\Phi_j: B^{\bar n}(\eta_1) \to \Phi_j(B^{\bar n}(\eta_1)) \subseteq \mathbb{R}^n$ as
$$
\Phi_j(u) = e^{u\cdot (\tau_j^{-1} \mathbf{X})_{J}}x_j,
$$
and 
$$
B_{(\tau_j^{-1} \mathbf{X}, \mathbf{d})}(x_j, \Vec{\xi_2}) \subseteq \Phi_j(B^{\bar n}(\eta_1)) \subseteq B_{(\tau_j^{-1} \mathbf{X}, \mathbf{d})}(x_j, \Vec{\xi_1}).
$$
Let $Y_l^j$ be the pullback of $\tau_j^{-\theta_l} X_l$ under the map $\Phi_j$ to $B^{\bar n}(\eta_1)$. $\mathbf{Y}^{j}$ denotes the matrix with column vectors $\{Y_l^j\}_{l=1}^q$, and $\mathbf{Y}_J^{j}$ denotes the matrix with column vectors $\{Y_l^j\}_{l\in J}$. Then for every $j$,
\begin{align*}
\big| \det_{\bar n \times \bar n} \mathbf{Y}^{j} \big| \approx \big| \det \mathbf{Y}_{J}^{j} \big| \approx 1 \text{~on~} B^{\bar n}(\eta_1), \quad \big\|Y_l^j \big\|_{C^m(B^{\bar n}(\eta_1))} \lesssim_{m\vee 2} 1 \text{~for~} 1\leq l \leq q, m\in \mathbb{N}.
\end{align*}

Pick a $p_0\in J\cap \{1, \ldots, s\}$. Then for every $j$,
$$
|Y^j_{p_0}|\approx 1,  \text{~on~} B^{\bar n}(\eta_1).
$$
Let $Z^{j}_l$ be the pullback of $\tau_j^{-\theta(d_l)} X_l$ under the map $\Phi_j$ to $B^{\bar n}(\eta_1)$. Then 
$$
Z_l^{j}= \left\{
\begin{aligned}
& Y_l^j, ~1\leq l\leq s,\\
& \tau_j^{-\theta(d_l)+\theta_{l}} Y_l^j, ~s+1 \leq l \leq q.
\end{aligned}
\right.
$$
Thus $\big|Z_{p_0}^{j} \big| \approx 1 \text{~on~} B^{\bar n}(\eta_1)$ for every $j$, and $\big|Z_l^{j}\big|\to 0 \text{~in~} C^\infty(B^{\bar n}(\eta_1))$ as $j\to \infty$ for $s+1\leq l \leq q$.

Similarly to what was done for Proposition \ref{ppp}, one can prove that there exists a subsequence of the above $\{(x_j,\tau_j)\}_j$, still denoted as $\{ (x_j,\tau_j)\}_j$, such that $W^j(t):=\Phi_j^*(W( \tilde \tau_j^{-1} t))$ converges in $C^\infty(B^N(\rho)\times B^{\bar n}(\eta_1))$ to some $\overline W(t)$ as $j\to \infty$, and associated to each $(X_\alpha,d_\alpha)\in \mathcal{L}(\mathcal{P}\cup \mathcal{N})$, $Z^{j}_\alpha := \Phi_j^*(\tau_j^{-\theta(d_\alpha)} X_\alpha)$ converges in $C^\infty(B^{\bar n}(\eta_1))$ to some $Z_\alpha$ as $j\to \infty$. Denote $Z_l = \lim_j Z_l^j$. Similarly to what was done for Lemma \ref{AA}, one can prove that the Taylor series of $\overline{W}(t,u)$ in $t$ is
\begin{align*}
\overline{W}(t,u) \sim \sum_{|\alpha|>0} t^\alpha Z_\alpha(u).
\end{align*}
Then the pure and nonpure power sets of $\overline{W}(t,u)$ are
\begin{align*}
&\overline{\mathcal{P}}=\{(Z_\alpha ,d_\alpha): (X_\alpha,d_\alpha)\in \mathcal{P}\},\\
&\overline{\mathcal{N}}=\{(Z_\alpha ,d_\alpha): (X_\alpha, d_\alpha)\in \mathcal{N}\}.
\end{align*}
Let $\overline{\mathcal{F}}' =\big\{ \big( Z_l, d_l\big): 1\leq l \leq q\big\}$. Then 
\begin{align*}
\mathcal{L}(\overline{\mathcal{P}} \cup \overline{\mathcal{N}}) = \{(Z_\alpha,d_\alpha): (X_\alpha, d_\alpha)\in \mathcal{L}(\mathcal{P} \cup \mathcal{N})\}, \quad \overline{\mathcal{F}}'\subseteq \mathcal{L}(\overline{\mathcal{P}}\cup \overline{\mathcal{N}}).
\end{align*}

On $B^{\bar n}(\eta_1)$, $\big|Z_{p_0}\big| \approx 1$, and $Z_l =0$ for $s+1\leq l \leq q$. Thus at least one vector field in $\overline{\mathcal{F}}'$ is nonzero at $0$, and the only non-vanishing elements in $\overline{\mathcal{F}}'$ are $(Z_l, d_l)$ for $1\leq l \leq s$, which are of the same nonpure degree. Similarly to what was done for Proposition \ref{juemi}, one can prove that $\overline{\mathcal{F}}'$ strongly controls every element of $\mathcal{L}(\overline{\mathcal{P}} \cup \overline{\mathcal{N}})$ on $B^{\bar n}(\eta_1)$. This still holds if we remove from $\overline{\mathcal{F}}'$ the elements vanishing on $B^{\bar n}(\eta_1)$. Every element of $\overline{\mathcal{F}}'$ has the same nonpure degree $d_1= \cdots = d_s$. Then the only non-vanishing vector fields in $\mathcal{L}(\overline{\mathcal{P}}\cup \overline{\mathcal{N}})$ are with degrees greater than or equal to $d_1$ coordinate-wise. Thus $\overline{\mathcal{F}}' \subseteq \overline{\mathcal{P}} \cup \overline{\mathcal{N}}$. For every $(Z_\alpha, d_\alpha)\in \mathcal{L}(\overline{\mathcal{P}}\cup \overline{\mathcal{N}})$ with $d_\alpha \geq d_1$ coordinate-wise but $d_\alpha\neq d_1$, we have $\theta(d_\alpha)>\theta_1$, and since $\mathcal{F}'$ strongly controls the corresponding $(X_\alpha, d_\alpha)\in \mathcal{L}(\mathcal{P}\cup \mathcal{N})$, we have on $B^{\bar n}(\eta_1)$,
\begin{align*}
Z_\alpha  =\lim_{j} \Phi_{j}^*\big(\tau_j^{-\theta(d_\alpha)}X_\alpha\big) =\lim_j  \sum_{1\leq l \leq s} \tau_j^{-\theta(d_\alpha)+\theta_1} \,(c_{\alpha,l} \circ \Phi_j)\, Y_l^j + \lim_j \sum_{\substack{s+1\leq l \leq q\\ d_l \leq d_\alpha}} \tau_j^{-\theta(d_\alpha)+\theta_{s+1}} \,(c_{\alpha,l}\circ \Phi_j)\, Y_l^j,
\end{align*}
where the $c_{\alpha,l}$ are $C^\infty$ on some neighborhood of $0$. The above first term vanishes since $\theta(d_\alpha)>\theta_1$, and the above second term vanishes since $\theta(d_\alpha) \geq \theta(d_l) >\theta_{s+1}$ for $s+1 \leq l \leq q$. Therefore the only non-vanishing vector fields in $\mathcal{L}(\overline{\mathcal{P}}\cup \overline{\mathcal{N}})$ are with the same nonpure degree $d_1$.

Since $\mathcal{F}'$ strongly controls $W(t,x)$ on $B^N(\rho)\times \Omega$, we have on $B^N(\rho)\times B^{\bar n}(\eta_1)$,
\begin{equation}\label{hello}
\begin{aligned}
\overline{W}(t,u) &= \lim_j D\Phi_j(u)^{-1} \Big(W\big(\tilde \tau_j^{-1} t, \Phi_j(u)\big)\Big) \\
&= \lim_j D\Phi_j(u)^{-1} \Big(\sum_{\substack{1\leq l \leq q\\ \text{deg}\,(\alpha) = d_l}} c_{l, \alpha}\big(\tilde \tau_j^{-1} t, \Phi_j(u)\big)\, t^{\alpha}\, \tau_j^{-\theta(d_l)} X_l(\Phi_j(u)) \Big)\\
& = \lim_j \sum_{\substack{1\leq l \leq q\\ \text{deg}\,(\alpha) = d_l}} c_{l, \alpha}\big(\tilde \tau_j^{-1} t, \Phi_j(u)\big) \,t^{\alpha} \,Z_l^j(u),
\end{aligned}
\end{equation}
where the $c_{l,\alpha}$ are $C^\infty$ on a neighborhood of $B^N(\rho) \times \Omega$. Then for all $m\in \mathbb{N}$, 
\begin{align*}
&\quad \sup_{j} \Big\| c_{l, \alpha}\big(\tilde \tau_j^{-1} t, \Phi_j(u)\big) \Big\|_{C^m(B^N(\rho)\times B^{\bar n}(\eta_1))}\\
&\lesssim_{m\vee 2} \sup_{j} \sum_{|\beta| + |\beta'| \leq 2m} \Big\| (\mathbf{Y}^j)^\beta (\partial_t^{\beta'} c_{l, \alpha})\big(\tilde \tau_j^{-1} t, \Phi_j(u)\big)  \Big\|_{C(B^N(\rho)\times B^{\bar n}(\eta_1))}\\
& \leq \sup_{j} \sum_{|\beta| + |\beta'|\leq 2m} \Big\| (\tau_j^{-1} \mathbf{X})^\beta \partial_t^{\beta'} c_{l, \alpha}( t,x) \Big\|_{C(B^N(\rho)\times \Omega)} <\infty.
\end{align*}
Similarly for $1\leq i \leq N$,
\begin{align*}
\lim_{j\to \infty} \Big\|\partial_{t_i}\big( c_{l, \alpha}(\tilde \tau_j^{-1} t, \Phi_j(u))\big) \Big\|_{C(B^N(\rho)\times B^{\bar n}(\eta_1))} = \lim_{j\to \infty} \tau_j^{-\theta(e_i)} \Big\| (\partial_{t_i} c_{l, \alpha})(\tilde \tau_j^{-1} t, \Phi_j(u))  \Big\|_{C(B^N(\rho)\times B^{\bar n}(\eta_1))} = 0.
\end{align*}
Hence there exists a subsequence of $\{c_{l, \alpha}(\tilde \tau_j^{-1} t, \Phi_j(u))\}_j$ convergent in $C^\infty(B^N(\rho)\times B^{\bar n}(\eta_1))$ to some $\bar c_{l,\alpha}(t,u)$, which is independent of $t$ and is denoted as $\bar c_{l,\alpha}(u)$. Thus by taking the limit of this subsequence of (\ref{hello}), we have on $B^N(\rho)\times B^{\bar n}(\eta_1)$,
$$
\overline{W}(t,u) = \sum_{\substack{1\leq l \leq q\\ \text{deg}\,(\alpha) = d_l}} t^{\alpha} \big(\bar c_{l,\alpha}(u) Z_l(u)\big)  = \sum_{\substack{1\leq l \leq s\\ \text{deg}\,(\alpha) = d_1}} t^{\alpha} \big(\bar c_{l,\alpha}(u) Z_l(u)\big).
$$
Therefore $\overline{W}(t,u)$ can be written as a finite sum on $B^N(\rho)\times B^{\bar n}(\eta_1)$:
$$
\overline W(t,u)  = \sum_{k=1}^Q t^{\alpha_k} \overline{X}_k(u),
$$
where the $\alpha_k\in \mathbb{N}^N$ are such that $\text{deg}\,(\alpha_1) = \cdots = \text{deg}\,(\alpha_Q)\in [0,\infty)^\nu$ is a nonpure degree, the $\overline{X}_k$ commute on $B^{\bar n}(\eta_1)$, and there exists $1\leq k \leq Q$ such that $\overline{X}_k(0)\neq 0$. This is because all the non-vanishing vector fields in $\mathcal{L}(\overline{\mathcal{P}} \cup \overline{\mathcal{N}})$ are with the same nonpure degree $d_1$, and there exists $(Z_{p_0},d_{p_0})\in \overline{\mathcal{F}}' \subseteq \overline{\mathcal{P}} \cup \overline{\mathcal{N}}$ satisfying $\big| Z_{p_0} \big| \approx 1$.

By Lemma \ref{12.1}, $\overline W(t,u)$ corresponds to a $C^\infty$ function $\bar \gamma_t(u)$ defined on a neighborhood of $(t,u)=(0,0)$ in $\mathbb{R}^{N}\times \mathbb{R}^{\bar n}$ with $\bar \gamma_0(u)\equiv u$. Similarly to what was done for Proposition \ref{UU}, one can prove that the operator
$$
\overline Tf(u) = \bar \psi(u)\int f(\bar{\gamma}_t(u))K(t)\,dt
$$
is bounded on $L^p(\mathbb{R}^{\bar n})$ for every $\bar \psi \in C_c^\infty(\mathbb{R}^{\bar n})$ with sufficiently small support, and every $K\in \mathcal{K}( N,  e, \bar a,  \nu)$ with sufficiently small $\bar a>0$.

\section{Proof of Proposition \ref{3}}\label{as}
Roughly speaking, if all the non-vanishing vector fields are associated with the same nonpure degree, then after ``zoom in'' on such $\gamma_t(x)$, we are able to construct explicit kernels with arbitrarily small supports such that the corresponding operators are unbounded.

Fix $N, \nu>1, e$ satisfying (\ref{e}) and fix $n\in \mathbb{N}_{>0}, p\in (1, \infty)$. Suppose $\gamma_t(x)$ is a $C^\infty$ function of $(t,x)$ defined on $B^N(\rho)\times \Omega$ with $\gamma_0(x)\equiv x$, for some $\rho>0$ and some neighborhood $\Omega$ of $0$ in $\mathbb{R}^n$. And suppose the corresponding $W(t,x)$ can be written as a finite sum on $B^N(\rho)\times \Omega$:
$$
W(t,x) = \sum_{l=1}^Q t^{\alpha_l} X_l(x),
$$
where the $\alpha_l \in \mathbb{N}^N$ are such that the degrees $\text{deg}\,(\alpha_1) = \cdots = \text{deg}\,(\alpha_Q)=d$ for some nonpure $d\in [0,\infty)^\nu$, the $X_l$ commute on $\Omega$, and there exists $1\leq l\leq Q$ such that $X_l(0) \neq 0$. Without loss of generality, by reordering $\alpha_1, \ldots, \alpha_Q$, we assume $X_1(0) \neq 0$. By shrinking $\Omega$, we can assume every $X_l$ has finite $C^m(\Omega)$ norm for all $m\in \mathbb{N}$. Suppose by contradiction that there exists $0<a\leq \rho$, and a neighborhood $U \subseteq \Omega$ of $0$ in $\mathbb{R}^n$, such that for every $\psi\in C_c^\infty(U)$, every $K\in \mathcal{K}(N,e,a,\nu)$, the operator
$$
Tf(x)=\psi(x)\int f(\gamma_t(x)) K(t)\,dt
$$
is bounded on $L^p(\mathbb{R}^n)$. Fix such a $\psi$ with $\psi(0) \neq 0$.

Write $\alpha_1= (\alpha_1^1, \ldots, \alpha^{N}_1)$. Let 
\begin{align*}
&B=\{ 1\leq i \leq N: \alpha_1^i\neq 0\}, \quad B^C=\{1\leq i \leq  N: \alpha_1^i=0\}, \quad |B| = \bar N,\\
&A = \{1\leq \mu \leq \nu: e_i^\mu \neq 0 \text{~for some~} i\in B\}, \quad |A| = \bar \nu.
\end{align*}
Note $\bar N, \bar \nu, \bar e$ satisfy (\ref{e}). Let $t_B$ be the vector consisting of the coordinates $t_i$ of $t$ with $i\in B$, and let $t_{B^C}$ be the vector consisting of the coordinates $t_i$ of $t$ with $i\in B^C$. Without loss of generality, we can assume $B=\{1, 2, \ldots, \bar N\}$, and $A = \{1, 2, \ldots, \bar \nu\}$. Let 
$$
\bar e =\{e_i^\mu: 1\leq i \leq \bar N, 1\leq \mu \leq \bar \nu\}, \quad \bar e_i = (e_i^1, \ldots, e_i^{\bar \nu}) \text{ for } 1\leq i \leq \bar N.
$$
Then $d =\text{deg}\,(\alpha_l) \in (0,\infty)^{\bar \nu} \times \{0\} \times \cdots \times \{0\}$. For $d =(d^1, \ldots, d^\nu)$, denote $\bar d = (d^1, \ldots, d^{\bar \nu})$. Since $d$ is nonpure, we have $\bar \nu>1$, and there exists $0\neq \textbf{n} \in \mathbb{R}^{\bar \nu}$ perpendicular to $\bar d$.

Denote $\gamma_{t_B}(x) := \gamma_{t}(x)\big|_{t_{B^C}=0}, W(t_B,x) := W(t,x)\big|_{t_{B^C}=0} = \sum_{l=1}^Q t^{\alpha_l} \big|_{t_{B^C}=0} X_l(x)$. If $\alpha_l^i \neq 0$ for some $i\in B^C$, then the term $t^{\alpha_l} \big|_{t_{B^C}=0} X_l(x)$ vanishes. From now on we use $\gamma_{t_B}(x), W(t_B,x)$ instead of $\gamma_t(x), W(t,x)$, respectively. So without loss of generality, we assume every $\alpha_l$ satisfies $\alpha_l^i =0$ for $i\in B^C$, and 
$$
W(t_B, x) = \sum_{l=1}^Q t_B^{\alpha_l} X_l(x).
$$
Set $\theta = \sum_{\mu=1}^{\bar \nu} d^\mu$, and set $\theta(e_i)= \sum_{\mu=1}^{\bar \nu} e_i^\mu$ for $1\leq i \leq \bar N$. The multi-parameter dilations on $\mathbb{R}^{\bar N}$ using $\bar e$ are defined by $\delta t_B = (\delta^{\bar e_1} t_1, \ldots, \delta^{\bar e_{\bar N}} t_{\bar N})$, for $\delta \in [0,\infty)^{\bar \nu}$. For a function $\varsigma(t_B)$ on $\mathbb{R}^{\bar N}$, and for $\delta \in (0,\infty)^{\bar \nu}$, denote $\varsigma^{(\delta)}(t_B):= \delta^{\bar e_1 + \cdots + \bar e_{\bar N}} \varsigma(\delta t_B)$. For $\tau \geq 1$, denote $\tilde \tau := ( \tau^{b_1},\ldots, \tau^{b_{\bar \nu}})$. For $2^k = (2^{k_1}, \ldots, 2^{k_{\bar \nu}})$, denote $\tilde \tau 2^k := (\tau^{b_1} 2^{k_1}, \ldots, \tau^{b_{\bar \nu}} 2^{k_{\bar \nu}})$. Then $\tilde \tau^{-1} t_B = (\tau^{-\theta(e_1)} t_1, \ldots, \tau^{-\theta(e_{\bar N})} t_{\bar N})$, $2^{-k} \tilde \tau^{-1} t_B = (2^{-k\cdot \bar e_1} \tau^{-\theta(e_1)} t_1, \ldots, 2^{-k\cdot \bar e_{\bar N}} \tau^{-\theta(e_{\bar N})} t_{\bar N})$. Let $(\mathbf{X}, \mathbf{\bar{d}}) = \{ (X_l, \bar d)\}_{l=1}^Q$, and $(\tilde \tau^{-1} \mathbf{X}, \mathbf{\bar d}) = \{ (\tau^{-\theta} X_{l}, \bar d)\}_{l=1}^Q$. For $\xi>0$, denote $\Vec{\xi} = (\xi, \ldots, \xi) \in (0,\infty)^{\bar \nu}$. See other notations in Definition \ref{indoor}.

There exists $\xi_1>0$ such that $(\mathbf{X},\mathbf{\bar d})$ satisfies $\mathcal{C}(0, \Vec{\xi_1}, \Omega)$. Since the $X_l$ commute on $\Omega$, we have $(\tilde \tau^{-1} \mathbf{X}, \mathbf{\bar d})$ satisfies the conditions of Theorem \ref{16} for $x_0=0$ uniformly in $\tau \geq 1$. Let $\bar n:= \text{rank}\, \mathbf{X}(0) = \text{rank}\,(\tilde \tau^{-1} \mathbf{X})(0)$. Since $X_{1}(0)\neq 0$, we have $\bar n>0$, and we can expand $X_{1}(0)$ to a basis $\mathbf{X}_{J}(0)$ of $\text{span}\,\mathbf{X}(0)$ for some $1\in J\in \mathcal{I}(\bar n,Q)$. We have
$$
\frac{|\det_{\bar n\times \bar n} (\tilde \tau^{-1} \mathbf{X})_J(0)|_\infty}{|\det_{\bar n\times \bar n} (\tilde \tau^{-1}\mathbf{X})(0)|_\infty}= \frac{|\det_{\bar n\times \bar n} \mathbf{X}_J(0)|_\infty}{|\det_{\bar n\times \bar n} \mathbf{X}(0)|_\infty},
$$
which is a fixed positive number, denoted as $\zeta$. Without loss of generality, we can assume $J=\{1, \ldots, \bar n\}$. In this section, $m$-admissible constants are with respect to $(\mathbf{X}, \mathbf{\bar d})$ and $\zeta$; see Definition \ref{uniformadmissible}.

By Theorem \ref{16}, there exist $2$-admissible constants $\eta_1>0, 0<\xi_2\leq \xi_1$, such that for every $\tau \geq 1$, we can define a diffeomorphism $\Phi_{\tau}: B^{\bar n}(\eta_1)\to \Phi_{\tau}(B^{\bar n}(\eta_1)) \subseteq \mathbb{R}^n$ as
$$
\Phi_{\tau} (u)= e^{u\cdot(\tilde \tau^{-1} \mathbf{X})_{J}}0,
$$
such that
$$
B_{(\mathbf{X},\mathbf{\bar d})}(0,\xi_2 \tilde \tau^{-1})\subseteq \Phi_{\tau}(B^{\bar n}(\eta_1)) \subseteq B_{(\mathbf{X},\mathbf{\bar d})}(0,\xi_1 \tilde \tau^{-1}).
$$
For $1\leq l \leq Q$, if we let $Z^\tau_l$ be the pullback of $\tau^{-\theta} X_l$ under the map $\Phi_\tau$ to $B^{\bar n}(\eta_1)$, then for $\tau \geq 1$, $\|Z_l^\tau\|_{C^m(B^{\bar n}(\eta_1))} \lesssim_{m\vee 2} 1$ for all $m\in \mathbb{N}$, and $|\det \mathbf{Z}^\tau_J|\approx 1$ on $B^{\bar n}(\eta_1)$, where $\mathbf{Z}_J^\tau$ denotes the matrix with column vectors $\{Z_l^\tau\}_{l\in J}$.

\begin{proposition}
There exists a sequence $\tau_j\to \infty$, such that $W^{\tau^j}(t_B):= \Phi_{\tau_j}^*(W(\tilde \tau_j^{-1} t_B))$ converges in $C^\infty(B^{\bar N}(\rho)\times B^{\bar n}(\eta_1))$ to some $\overline{W}(t_B)$, and $Z_l^{\tau_j}$ converges in $C^\infty(B^{\bar n}(\eta_1))$ to some $Z_l$. Moreover, on $B^{\bar N}(\rho)\times B^{\bar n}(\eta_1)$,
$$
\overline{W}(t_B,u)  =  \sum_{l=1}^Q t_B^{\alpha_l} Z_{l}(u).
$$
\end{proposition}

\begin{proof}
On $B^{\bar N}(\rho)\times B^{\bar n}(\eta_1)$,
\begin{equation}\label{tau}
W^{\tau}(t_B,u)  := D\Phi_\tau(u)^{-1} \big(W(\tilde \tau^{-1} t_B, \Phi_\tau(u))\big) =  D\Phi_{\tau}(u)^{-1}\Big(\sum_{l=1}^Q t_B^{\alpha_l} \tau^{-\theta} X_l\big(\Phi_\tau(u)\big)\Big) = \sum_{l=1}^Q t_B^{\alpha_l} Z_{l}^{\tau}(u).
\end{equation}
And since $\{Z_l^{\tau}: \tau \geq 1\}$ is bounded in $C^\infty(B^{\bar n}(\eta_1))$, $\{W^\tau(t_B,u): \tau \geq 1\}$ is bounded in $C^\infty(B^{\bar N}(\rho)\times B^{\bar n}(\eta_1))$. Therefore there exists a sequence $\tau_j \to \infty$, such that $W^{\tau^j}$ converges in $C^\infty(B^{\bar N}(\rho)\times B^{\bar n}(\eta_1))$ to some $\overline W$, $Z_l^{\tau_j}$ converges in $C^\infty(B^{\bar n}(\eta_1))$ to some $Z_l$, and by taking the limit of this sequence of (\ref{tau}),
$$
\overline{W}(t_B,u)  =  \sum_{l=1}^Q t^{\alpha_l}_B Z_{l}(u), \text{~on~} B^{\bar N}(\rho)\times B^{\bar n}(\eta_1).
$$ 
\end{proof}

Thus $\big| \det \mathbf{Z}_J \big| \approx 1$ on $B^{\bar n}(\eta_1)$, $\{Z_1, \ldots, Z_{\bar n}\}$ is a basis on $B^{\bar n}(\eta_1)$, and the $Z_l$ commmute on $B^{\bar n}(\eta_1)$. Therefore for $\bar n < s \leq Q$,  $Z_s = \sum_{l=1}^{\bar n} c_{s,l} Z_l$ on $B^{\bar n}(\eta_1)$, for some constants $c_{s,l}$. By the commutativity of the $Z_l$,
\begin{align*}
\frac{d}{d\epsilon}\Big|_{\epsilon =1} e^{\sum_{l=1}^Q \frac{1}{|\alpha_l|} \epsilon^{|\alpha_l|} t_B^{\alpha_l} Z_{l}} e^{- \sum_{l=1}^Q \frac{1}{|\alpha_l|} t_B^{\alpha_l} Z_{l}}u = \sum_{l=1}^Q t_B^{\alpha_l} Z_{l}(u) = \overline W(t_B,u).
\end{align*}
So $\overline{W}(t_B,u)$ corresponds to $\bar \gamma_{t_B}(u) := e^{\sum_{l=1}^Q \frac{1}{|\alpha_l|} t_B^{\alpha_l} Z_{l}}u$. Similarly, by the commutativity of the $X_l$, $W(t_B,x)$ corresponds to $\gamma_{t_B}(x) = e^{\sum_{l=1}^Q \frac{1}{|\alpha_l|} t_B^{\alpha_l} X_{l}}x$. For $k\in \mathbb{N}$, we have $(2^{-k\textbf{n}} \tilde \tau^{-1} t_B)^{\alpha_l} = 2^{-k\textbf{n}\cdot \bar d} \tau^{-\theta} t_B^{\alpha_l} = \tau^{-\theta} t_B^{\alpha_l}$. Thus $\gamma_{2^{-k\textbf{n}} \tilde \tau^{-1} t_B}(x) = \gamma_{\tilde \tau^{-1} t_B}(x)$.

Similar to Corollary \ref{VV} and Lemma \ref{wujinzang}, there exists $0<\rho'\leq \rho$, $N_0 \in \mathbb{N}_{>0}$, and a $2$-admissible constant $0<\eta_0 < \frac{1}{2}\eta_1$, such that for $j \geq N_0$, $\Phi_{\tau_j}^{-1} \circ \gamma_{2^{-k\textbf{n}} \tilde \tau_j^{-1} t_B} \circ \Phi_{\tau_j}(u) = \Phi_{\tau_j}^{-1} \circ \gamma_{\tilde \tau_j^{-1} t_B} \circ \Phi_{\tau_j}(u)$ is well-defined and convergent to $\bar \gamma_{t_B}(u)$ as $j \to \infty$, uniformly on $B^{\bar N}(\rho')\times B^{\bar n}(\eta_0)$. Fix $0< \bar a \leq \min \{\frac{a}{2}, \frac{\rho'}{2}\}$. Fix an arbitrary sequence $\{\varsigma_{k}(t_B)\}_{k\in \mathbb{N}}$ bounded in $C_c^\infty(\overline{B^{\bar N}(\bar a)})$ satisfying $\int \varsigma_k(t_B) \,dt_B^\mu \equiv 0$, for $1\leq \mu \leq \bar \nu$, $k\in \mathbb{N}$, where $t_B^\mu$ denotes the vector consisting of the coordinates $t_i$ of $t$ with $1\leq i \leq \bar N$ and $e_i^\mu \neq 0$.

\begin{lemma}
For every $M \in \mathbb{N}_{>0}$ and $\tau\gg M$, there exists $\{\tilde \varsigma_i(t_B)\}_{i\in \mathbb{N}^{\bar \nu}} \in \overline{\mathcal{K}}(\bar N, \bar e, \bar a, \bar \nu)$ such that
$$
\sum_{\substack{k\in \mathbb{N}\\ k \leq M}} \varsigma_k^{(\tilde \tau 2^{k\textbf{n}})}(t_B) = \sum_{i\in \mathbb{N}^{\bar \nu}} \tilde \varsigma_i^{(2^i)}(t_B).
$$
Moreover, the above $\{\tilde \varsigma_i(t_B)\}_{i\in \mathbb{N}^{\bar \nu}}$ can be chosen to be bounded in $\overline{\mathcal{K}}(\bar N, \bar e, \bar a, \bar \nu)$ for $M\in \mathbb{N}_{>0}$ and $\tau\gg  M$. In particular, $\sum_{\substack{k\in \mathbb{N}\\ k \leq M}} \varsigma_k^{(\tilde \tau 2^{k\textbf{n}})}(t_B) \in \mathcal{K}(\bar N, \bar e, a, \bar \nu)$ for $M\in \mathbb{N}_{>0}, \tau\gg  M$.
\end{lemma}

\begin{proof}
Denote $\log_2(\tilde \tau) := (\log_2(\tau^{b_1}), \ldots, \log_2(\tau^{b_{\bar \nu}}))$. For $\tau\gg  M$, $\log_2(\tilde \tau) +k\textbf{n} \in (0,\infty)^{\bar \nu}$ for $k \leq M$, then we can define $\{\tilde \varsigma_i(t_B)\}_{i\in \mathbb{N}^{\bar \nu}}$ as
$$
\tilde \varsigma_i(t_B) = \sum_{\substack{k:\,k\textbf{n} +\log_2(\tilde \tau) -i \in [0,1)^{\bar \nu} \\ k\leq M}} \varsigma_k^{(\tilde \tau 2^{k\textbf{n} -i})}(t_B).
$$
The above $\tilde \tau 2^{k\textbf{n} -i} \in [1,2)^{\bar \nu}$, and the number of terms in the above sum is bounded by a constant only depending on $\textbf{n}$, ($\tilde \varsigma_i(t_B) =0$ when there exists no $k\leq M$ such that $k\textbf{n} + \log_2(\tilde \tau) -i \in [0,1)^{\bar \nu}$). Then $\{\tilde \varsigma_i(t_B)\}_{i\in \mathbb{N}^{\bar \nu}}$ is bounded in $\overline{\mathcal{K}}(\bar N, \bar e, \bar a, \bar \nu)$ for $M\in \mathbb{N}_{>0}$, $\tau\gg  M$. Note for each $k$, there exists a unique $i\in [0,1)^{\bar \nu}$ such that $k\textbf{n}+ \log_2(\tilde \tau)-i \in [0,1)^{\bar \nu}$. Therefore 
$$
\sum_{k\leq M} \varsigma_k^{(\tilde \tau 2^{k\textbf{n}})}(t_B) = \sum_{i\in \mathbb{N}^{\bar \nu}} \sum_{\substack{k:\,k\textbf{n} +\log_2(\tilde \tau) -i \in [0,1)^{\bar \nu} \\ k\leq M}} \varsigma_k^{(\tilde \tau 2^{k\textbf{n}})}(t_B) = \sum_{i\in \mathbb{N}^{\bar \nu}} \tilde \varsigma_i^{(2^i)}(t_B).
$$
\end{proof}

Thus similarly to what was done for Lemma \ref{yangtianxiang}, one can prove that the operator
$$
\psi(x) \sum_{\substack{k\in \mathbb{N} \\ k \leq M}} \int f(\gamma_{\tilde \tau^{-1} t_B}(x)) \varsigma_k^{(2^{k\textbf{n}})} (t_B)\,dt_B
$$
is bounded on $L^p(\mathbb{R}^n)$ uniformly in $M\in \mathbb{N}_{>0}, \tau\gg  M$. Fix arbitrary $\varsigma_1, \ldots, \varsigma_{\bar N}\in C_c^\infty(\overline{B^1(\frac{\bar a}{2\sqrt{\bar N}})})$ satisfying
$$
\int \varsigma_i(t_i)\,dt_i =0, \text{~for~} 1\leq i \leq \bar N.
$$
We can let $\varsigma_k(t_B) = \prod_{i=1}^{\bar N} \varsigma_i(t_i)$ for every $k$. Then the operator
\begin{align*}
T_{\tau,M}f(x) := \psi(x) \sum_{\substack{k\in \mathbb{N}\\ k \leq M}} \int f(\gamma_{2^{-k\textbf{n}} \tilde \tau^{-1}t_B} (x)) \prod_{i\in B} \varsigma_i(t_i)\,dt_B = \psi(x) (M+1) \int f(\gamma_{ \tilde \tau^{-1}t_B} (x)) \prod_{i\in B} \varsigma_i(t_i)\,dt_B
\end{align*}
is bounded on $L^p(\mathbb{R}^n)$ uniformly in $M\in \mathbb{N}_{>0}, \tau \gg   M$.

Fix a $\bar \psi \in C_c^\infty(B^{\bar n}(\eta_0))$ with $\bar \psi(0) \neq 0$. Similar to Proposition \ref{W}, there exists $C>0$ and $N_2 \geq N_0$, such that for $j\geq N_2$, $M\in \mathbb{N}_{>0}$,
$$
\Big\| \bar \psi \Phi_{\tau_j}^* T_{\tau_j ,M} (\Phi_{\tau_j}^{-1})^* \Big\|_{L^p (\mathbb{R}^{\bar n}) \to L^p(\mathbb{R}^{\bar n})} \leq C \big\| T_{\tau_j ,M} \big\|_{L^p(\mathbb{R}^n) \to L^p(\mathbb{R}^n)}.
$$
Thus $\bar \psi \Phi_{\tau_j}^* T_{\tau_j ,M} (\Phi_{\tau_j}^{-1})^*$ is bounded on $L^p(\mathbb{R}^{\bar n})$ uniformly in $M\in \mathbb{N}_{>0}, j\gg  M$.

For any $M\in \mathbb{N}_{>0}$ and $f\in C_c^\infty(\mathbb{R}^{\bar n})$, by the dominated convergence theorem,
\begin{align*}
T_Mf(u) :&= (M+1) \bar \psi(u) \int f(\bar \gamma_{t_B}(u)) \prod_{i\in B} \varsigma_i(t_i)\,dt_B  \\
&= \lim_{j\to \infty} \bar \psi(u) \frac{\psi(\Phi_{\tau_j}(u))}{\psi(0)} (M+1) \int f(\Phi_{\tau_j}^{-1} \circ \gamma_{\tilde \tau_j^{-1} t_B} \circ \Phi_{\tau_j}(u)) \prod_{i\in B} \varsigma_i(t_i)\,dt_B\\
&= \frac{1}{\psi(0)}\lim_{j \to \infty}\bar \psi \Phi_{\tau_j}^* T_{\tau_j,M} (\Phi_{\tau_j}^{-1})^*f(u).
\end{align*}
By Fatou's lemma,
\begin{align*}
\big\|T_Mf\big\|_{L^p(\mathbb{R}^{\bar n})} &\leq \frac{1}{\psi(0)} \liminf_{j\to\infty} \Big\| \bar \psi \Phi_{\tau_j}^* T_{\tau_j, M} (\Phi_{\tau_j}^{-1})^*f\Big\|_{L^p(\mathbb{R}^{\bar n})} \\
&\leq \Big( \frac{1}{\psi(0)} \sup_{\substack{M\in \mathbb{N}_{>0}\\j\gg  M}} \Big\|\bar \psi \Phi_{\tau_j}^* T_{\tau_j, M} (\Phi_{\tau_j}^{-1})^*\Big\|_{L^p(\mathbb{R}^{\bar n})\to L^p(\mathbb{R}^{\bar n})}\Big) \big\|f\big\|_{L^p(\mathbb{R}^{\bar n})}.
\end{align*}
Hence $T_M$ are bounded on $L^p(\mathbb{R}^{\bar n})$ uniformly in $M$.

Now we show $T_M$ are in fact not bounded uniformly in $M$, to reach a contradiction. It suffices to find $f\in C_c^\infty(\mathbb{R}^{\bar n})$ and $\varsigma_1, \ldots, \varsigma_{\bar N} \in C_c^\infty(\overline{B^1(\frac{\bar a}{2\sqrt{\bar N}})})$ with $\int \varsigma_i(t_i) =0$ for every $1\leq i \leq \bar N$, such that 
$$
\Big\|\bar \psi(u) \int f(\bar \gamma_{t_B}(u)) \prod_{i\in B} \varsigma_i(t_i)\,dt_B \Big\|_{L^p(\mathbb{R}^{\bar n})} \neq 0.
$$
By the dominated convergence theorem, $\bar \psi(u) \int f(\bar \gamma_{t_B}(u)) \prod_{i\in B} \varsigma_i(t_i)\,dt_B$ is continuous. Thus it suffices to find $f\in C_c^\infty(\mathbb{R}^{\bar n})$ and $\varsigma_1, \ldots, \varsigma_{\bar N} \in C_c^\infty(\overline{B^1(\frac{\bar a}{2\sqrt{\bar N}})})$ with $\int \varsigma_i(t_i) =0$ for $1\leq i \leq \bar N$, such that 
$$
\int f(\bar \gamma_{t_B}(0)) \prod_{i\in B} \varsigma_i(t_i)\,dt_B  = \int f(e^{\sum_{l=1}^Q \frac{1}{|\alpha_l|} t_B^{\alpha_l} Z_{l}}0) \prod_{i\in B} \varsigma_i(t_i)\,dt_B \neq 0.
$$

Define $h(\tilde u_1, \ldots, \tilde u_{\bar n}) = e^{\sum_{l=1}^{\bar n} \tilde u_l Z_l}0$. Note $\big| \det Dh(0) \big| = \big| \det Z_J(0) \big| \approx 1$. Then there exist neighborhoods $\tilde U, \tilde U' \subseteq B^{\bar n}(\eta_1)$ of $0$ such that $h$ is a diffeomorphism from $\tilde U$ onto $\tilde U'$. Let $f\in C_c^\infty(\mathbb{R}^{\bar n})$ be such that $f\circ h(\tilde u)= \tilde u_1$ for $\tilde u\in \tilde U$. In the following equation, $\{ V_l \}_{l=1}^{\bar  n}$ denotes the vector in $\mathbb{R}^{\bar n}$ with the $l$-th entry being $V_l$. By shrinking $\bar a>0$, we can assume $\big\{ \frac{1}{|\alpha_l|}t_B^{\alpha_l}+ \sum_{s>\bar n} \frac{1}{|\alpha_s|}t_B^{\alpha_s} c_{s,l} \big\}_{l=1}^{\bar  n} \in \tilde U$ for $t_B \in B^{\bar N}(\bar a)$. Then
\begin{align*}
\int f\big(e^{\sum_{l=1}^Q \frac{1}{|\alpha_l|} t_B^{\alpha_l} Z_{l}}0\big) \prod_{i\in B} \varsigma_i(t_i)\,dt_B &= \int f\Big( e^{\sum_{l=1}^{\bar n} \big(\frac{1}{|\alpha_l|}t_B^{\alpha_l}+ \sum_{s>\bar n} \frac{1}{|\alpha_s|}t_B^{\alpha_s} c_{s,l} \big) Z_l} 0\Big) \prod_{i\in B} \varsigma_i(t_i)\,dt_B\\
&= \int f\circ h \Big(\Big\{ \frac{1}{|\alpha_l|}t_B^{\alpha_l}+ \sum_{s>\bar n} \frac{1}{|\alpha_s|}t_B^{\alpha_s} c_{s,l} \Big\}_{l=1}^{\bar  n}\Big) \prod_{i\in B} \varsigma_i(t_i)\,dt_B\\
&= \int \Big( \frac{1}{|\alpha_1|}t_B^{\alpha_1}+ \sum_{s>\bar n} \frac{1}{|\alpha_s|}t_B^{\alpha_s} c_{s,1} \Big) \prod_{i\in B} \varsigma_i(t_i)\,dt_B.
\end{align*}
Then it remains to find $\varsigma_1, \ldots, \varsigma_{\bar N} \in C_c^\infty(\overline{B^1(\frac{\bar a}{2\sqrt{\bar N}})})$ such that
\begin{align*}
\int t_B^{\alpha_1}\prod_{i=1}^{\bar N}\varsigma_i(t_i)\,dt_B \neq 0,\quad \int t_B^{\alpha_l}\prod_{i=1}^{\bar N}\varsigma_i(t_i)\,dt_B =0 \text{~for~} l\neq 1,\quad \int \varsigma_i(t_i)\,dt_i=0 \text{~for~}1\leq i \leq \bar N.
\end{align*}

For $1\leq l \leq Q$,
$$
\int t_B^{\alpha_l}\prod_{i=1}^{\bar N}\varsigma_i(t_i)\,dt_B= \prod_{i=1}^{\bar N} \int t_i^{\alpha_l^i}\varsigma_i(t_i)dt_i.
$$
For $l\neq 1$, we have $\alpha_l \neq \alpha_1$, and thus there exists $1\leq i_l\leq \bar N$ such that $\alpha_1^{i_l} \neq \alpha_l^{i_l}$. So it suffices to find $\varsigma_1, \ldots, \varsigma_{\bar N} \in C_c^\infty(\overline{B^1(\frac{\bar a}{2\sqrt{\bar N}})})$ such that for $1\leq i \leq \bar N, 2 \leq l \leq Q$,
\begin{align*}
\int \varsigma_i(t_i)\,dt_i =0, \quad \int t_i^{\alpha_1^i}\varsigma_i(t_i)\,dt_i \neq 0, \quad \int t_{i_l}^{\alpha_l^{i_l}}\varsigma_{i_l}(t_{i_l})\,dt_{i_l}=0,
\end{align*}
or equivalently, for $1\leq i \leq \bar N$,
$$
\int \varsigma_i(t_i)\,dt_i =0, \quad \int t_i^{\alpha_1^i}\varsigma_i(t_i)\,dt_i \neq 0, \quad \int t_i^{\alpha_l^i} \varsigma_i(t_i)\,dt_i =0 \text{ for every } l \text{ satisfying } i_l = i.
$$
This can be achieved by the following lemma.

\begin{lemma}\label{phd}
For every $a>0$, there exists $\varsigma \in C_c^\infty(0,a) \subseteq C_c^\infty(\mathbb{R})$, such that $\int \varsigma(t) =0, \int t^{a_1}\varsigma(t) \neq 0, \int t^{a_2}\varsigma(t) = \cdots = \int t^{a_k}\varsigma(t)=0$, where $a_1, \ldots, a_k$ are distinct positive integers.
\end{lemma}

\begin{proof}
Fix a $\psi \in C_c^\infty(0,1), \psi \geq 0, \int \psi =1$. Then for each $m\in \mathbb{N}$, $b_m:= \int t^m \psi(t)\,dt >0$, and $1=b_0>b_1>b_2>\cdots>b_m >\cdots$. For every $x, r>0$, let 
$$
\psi_{x,r}(t) = \frac{1}{r}\psi(\frac{t-x}{r}).
$$
Then $\int \psi_{x,r} =1$. Let $\varsigma = c_1 \psi_{x_1,r_1}+\cdots +c_{k+1} \psi_{x_{k+1},r_{k+1}}$, and we need to find $c_1, \ldots, c_{k+1}, x_1, \ldots, x_{k+1}, r_1, \ldots, r_{k+1}$ such that $\varsigma$ satisfies all the required conditions.

$\int \varsigma =0$ if and only if $c_1 \int \psi_{x_1,r_1} +\cdots + c_{k+1}\int \psi_{x_{k+1},r_{k+1}} =0$, i.e., $c_1+\cdots + c_{k+1} =0$. For $1\leq l \leq k$,
\begin{align*}
\int t^{a_l} \varsigma = &c_1\int t^{a_l} \psi_{x_1,r_1} +\cdots + c_{k+1}\int t^{a_l}\psi_{x_{k+1},r_{k+1}}\\
=& c_1\int t^{a_l} \frac{1}{r_1}\psi(\frac{t-x_1}{r_1})\,dt + \cdots +c_{k+1} \int t^{a_l} \frac{1}{r_{k+1}}\psi(\frac{t-x_{k+1}}{r_{k+1}})\,dt\\
=&c_1\int r_1^{a_l}\Big(\frac{t-x_1+x_1}{r_1}\Big)^{a_l}\psi\Big(\frac{t-x_1}{r_1}\Big)\frac{dt}{r_1} + \cdots +\\
&c_{k+1}\int r_{k+1}^{a_l}\Big(\frac{t-x_{k+1}+x_{k+1}}{r_{k+1}}\Big)^{a_l}\psi\Big(\frac{t-x_{k+1}}{r_{k+1}}\Big)\frac{dt}{r_{k+1}}\\
=& c_1 \int r_1^{a_l} \Big( s+ \frac{x_1}{r_1}\Big)^{a_l} \psi(s)\,ds + \cdots + c_{k+1} \int r_{k+1}^{a_l} \Big(s+ \frac{x_{k+1}}{r_{k+1}}\Big)^{a_l} \psi(s)\,ds\\
=&c_1 r_1^{a_l}\sum_m \binom{a_l}{m}\big(\frac{x_1}{r_1}\big)^{a_l-m}b_m + \cdots + c_{k+1} r_{k+1}^{a_l}\sum_m \binom{a_l}{m} \big(\frac{x_{k+1}}{r_{k+1}}\big)^{a_l-m}b_m\\
=& c_1\sum_{m}\binom{a_l}{m}b_m r_1^mx_1^{a_l-m}+ \cdots + c_{k+1}\sum_{m}\binom{a_l}{m}b_m r_{k+1}^mx_{k+1}^{a_l-m}.
\end{align*}
Thus we need 
\begin{align*}
&\begin{bmatrix}
0\\
\int t^{a_1} \zeta \\
\cdots\\
0
\end{bmatrix}
=
\begin{bmatrix}
1&1 & \cdots & 1\\
\sum \binom{a_1}{m}b_mr_1^mx_1^{a_1-m} & \sum \binom{a_1}{m}b_mr_2^mx_2^{a_1-m} & \cdots & \sum\binom{a_1}{m}b_mr_{k+1}^mx_{k+1}^{a_1-m}\\
\cdots & \cdots & \cdots & \cdots\\
\sum \binom{a_k}{m}b_mr_1^mx_1^{a_k-m} & \sum \binom{a_k}{m}b_mr_2^mx_2^{a_k-m} &\cdots & \sum \binom{a_k}{m}b_mr_{k+1}^mx_{k+1}^{a_k-m}
\end{bmatrix}
\begin{bmatrix}
c_1\\
c_2\\
\cdots\\
c_{k+1}
\end{bmatrix} .
\end{align*}
To have a solution $(c_1,\ldots, c_{k+1})$ to the above equation, it suffices to have
\begin{align*}
&\det 
\begin{bmatrix}
1&1 & \cdots & 1\\
\sum \binom{a_1}{m}b_mr_1^mx_1^{a_1-m} & \sum \binom{a_1}{m}b_mr_2^mx_2^{a_1-m} & \cdots & \sum\binom{a_1}{m}b_mr_{k+1}^mx_{k+1}^{a_1-m}\\
\cdots & \cdots & \cdots & \cdots\\
\sum \binom{a_k}{m}b_mr_1^mx_1^{a_k-m} & \sum \binom{a_k}{m}b_mr_2^mx_2^{a_k-m} &\cdots & \sum \binom{a_k}{m}b_mr_{k+1}^mx_{k+1}^{a_k-m}
\end{bmatrix} \neq 0.    
\end{align*}

Let $c$ be a fixed number in $(0,1)$ to be determined. For $1\leq j \leq k+1$, let 
$$
r_j=r_1 c^{j-1},\quad x_j = x_1 c^{j-1}.
$$    
Denote
$$
y_0=1, \quad y_l=\sum_m \binom{a_l}{m}b_mr_1^mx_1^{a_l-m} \text{~for~}1\leq l\leq k.
$$
Let $a_0=0$, then $c^{a_0}=1$. Then we have
\begin{align*}
&\quad
\det \begin{bmatrix}
1&1 & \cdots & 1\\
\sum \binom{a_1}{m}b_mr_1^mx_1^{a_1-m} & \sum \binom{a_1}{m}b_mr_2^mx_2^{a_1-m} & \cdots & \sum\binom{a_1}{m}b_mr_{k+1}^mx_{k+1}^{a_1-m}\\
\cdots & \cdots & \cdots & \cdots\\
\sum \binom{a_k}{m}b_mr_1^mx_1^{a_k-m} & \sum \binom{a_k}{m}b_mr_2^mx_2^{a_k-m} &\cdots & \sum \binom{a_k}{m}b_mr_{k+1}^mx_{k+1}^{a_k-m}
\end{bmatrix}\\
&=\det
\begin{bmatrix}
y_0 & y_0 c^{a_0} & \cdots & y_0 (c^{a_0})^k\\
y_1 & y_1 c^{a_1} & \cdots & y_1 (c^{a_1})^k\\
\cdots & \cdots & \cdots & \cdots\\
y_k & y_k c^{a_k} & \cdots & y_k (c^{a_k})^k
\end{bmatrix} 
=y_0\cdot y_1 \cdot \cdots \cdot y_k \cdot \prod_{0\leq l<l'\leq k} (c^{a_{l'}}- c^{a_l}) \neq 0,
\end{align*}
since $y_0, y_1,\cdots, y_k>0$, and $a_0, a_1, \ldots, a_k$ are distinct nonnegative integers.

Now let $x_1=r_1 = \frac{a}{2}$, $c\in (0,1)$, $r_j=r_1 c^{j-1}$, and $x_j = x_1 c^{j-1}$ for $1\leq j \leq k+1$. We then have a solution $(c_1, \ldots, c_{k+1})$ for the above matrix equation, and $\varsigma = c_1\psi_{x_1,r_1} + \cdots + c_{k+1}\psi_{x_{k+1}, r_{k+1}}$ satisfies all the required conditions.
\end{proof}

\appendix
\section{Some lemmas and proofs}\label{jinzhangzhi}

\begin{proof}[\textbf{Proof of Lemma \ref{contraction}}]
We need to show there exists $0<\rho''\leq \rho_0$ such that $\bar w$ and every $w^j (j \geq N_0)$ map from $[0,1]\times B^{\bar N}(\rho'') \times B^{\bar n}(\eta_0)$ to $B^{\bar n}(\frac{1}{2}\eta_1)$.

We have on $[0,1] \times B^{\bar N}(\rho_0) \times B^{\bar n}(\eta_0)$,
\begin{align*}
&w^j(\epsilon, t_B, u) = u + \int_0^\epsilon \frac{1}{\sigma} W^j(\sigma t_B, w^j(\sigma, t_B,u))\,d\sigma,\\
&\bar w(\epsilon, t_B, u) = u + \int_0^\epsilon \frac{1}{\sigma} \overline W(\sigma t_B, \bar w(\sigma, t_B,u))\,d\sigma.
\end{align*}
Let $\rho'' \in (0, \rho_0]$ be a small number, to be chosen later. We have
\begin{align*}
\Big\| w^j(\epsilon,t_B, u)\Big\|_{C([0,1]\times B^{\bar N}(\rho'')\times B^{\bar n}(\eta_0))} & \leq \eta_0 + \Big\| \frac{1}{\epsilon} W^j(\epsilon t_B,  w^j(\epsilon, t_B,u))\Big\|_{C([0,1]\times B^{\bar N}(\rho'')\times B^{\bar n}(\eta_0))}\\
&\leq \eta_0 + \rho'' \Big\| DW^j\Big\|_{C(B^{\bar N}(\rho)\times B^{\bar n}(\eta_1))},\\
\Big\| \bar w(\epsilon,t_B, u)\Big\|_{C([0,1]\times B^{\bar N}(\rho'')\times B^{\bar n}(\eta_0))} &\leq \eta_0 + \rho'' \Big\| D\overline W\Big\|_{C(B^{\bar N}(\rho)\times B^{\bar n}(\eta_1))}.
\end{align*}
Since $\eta_0< \frac{1}{2}\eta_1$, we can choose $0<\rho'' \leq \rho_0$ such that 
$$
\eta_0 + \rho'' \max\Big\{ \Big\| D\overline W \Big\|_{C(B^{\bar N}(\rho)\times B^{\bar n}(\eta_1))}, \sup_{j}\Big\| DW^j\Big\|_{C(B^{\bar N}(\rho)\times B^{\bar n}(\eta_1))}\Big\}< \frac{1}{2} \eta_1.
$$
\end{proof}

\begin{proof}[\textbf{Proof of Lemma \ref{wujinzang}}]
We need to show there exists $0<\rho'\leq \rho''$ such that $w^j \to \bar w$ uniformly on $[0,1] \times B^{\bar N}(\rho')\times B^{\bar n}(\eta_0)$.

Let $\rho' \in (0,\rho'']$ be a small number, to be chosen later. We have
\begin{align*}
&\quad \Big\| w^j(\epsilon,t_B,u)-\overline w(\epsilon,t_B,u)\Big\|_{C([0,1]\times B^{\bar N}(\rho')\times B^{\bar n}(\eta_0))}\\
&=\Big\| \int_0^\epsilon \Big(\frac{d}{d\sigma} w^j(\sigma, t_B,u) - \frac{d}{d\sigma} \overline w(\sigma,t_B,u)\Big)\,d\sigma\Big\|_{C([0,1]\times B^{\bar N}(\rho')\times B^{\bar n}(\eta_0))}\\
&=\Big\| \int_0^\epsilon \Big(\frac{1}{\sigma} W^j(\sigma t_B, w^j(\sigma, t_B,u)) - \frac{1}{\sigma} \overline W(\sigma t_B, \overline w(\sigma, t_B, u))\Big)\,d\sigma\Big\|_{C([0,1]\times B^{\bar N}(\rho')\times B^{\bar n}(\eta_0))}\\
&\leq \Big\| \frac{1}{\sigma} W^j(\sigma t_B, w^j(\sigma, t_B,u)) - \frac{1}{\sigma} \overline W(\sigma t_B, \overline w(\sigma, t_B, u)) \Big\|_{C([0,1]\times B^{\bar N}(\rho')\times B^{\bar n}(\eta_0))}\\
&\leq \Big\| \frac{1}{\sigma} W^j(\sigma t_B, w^j(\sigma, t_B,u)) - \frac{1}{\sigma} \overline W(\sigma t_B, w^j(\sigma, t_B, u)) \Big\|_{C([0,1]\times B^{\bar N}(\rho')\times B^{\bar n}(\eta_0))} \\
&\quad + \Big\| \frac{1}{\sigma} \overline W(\sigma t_B, w^j(\sigma, t_B,u)) - \frac{1}{\sigma} \overline W(\sigma t_B, \overline w(\sigma, t_B, u)) \Big\|_{C([0,1]\times B^{\bar N}(\rho')\times B^{\bar n}(\eta_0))}\\
&\leq \Big\| |t_B| \big| D(W^j-\overline W)\big| \Big\|_{C(B^{\bar N}(\rho')\times B^{\bar n}(\frac{1}{2}\eta_1))} + \\
& \quad  \Big\| |t_B||D^2\overline W|\Big\|_{C(B^{\bar N}(\rho')\times B^{\bar n}(\frac{1}{2}\eta_1))} \Big\| w^j(\sigma, t_B, u)- \overline w(\sigma, t_B, u)\Big\|_{C([0,1]\times B^{\bar N}(\rho')\times B^{\bar n}(\eta_0))}\\
&\leq \rho' \Big\|W^j-\overline W\Big\|_{C^1(B^{\bar N}(\rho)\times B^{\bar n}(\frac{1}{2}\eta_1))}\\
&\quad  + \rho' \Big\|\overline W\Big\|_{C^2(B^{\bar N}(\rho)\times B^{\bar n}(\frac{1}{2}\eta_1))} \Big\| w^j(\epsilon, t_B, u) - \overline w(\epsilon, t_B,u)\Big\|_{C([0,1]\times B^{\bar N}(\rho')\times B^{\bar n}(\eta_0))}.
\end{align*}
Choose $0<\rho' \leq \rho''$ such that
$$
\rho' \Big\| \overline W \Big\|_{C^2(B^{\bar N}(\rho)\times B^{\bar n}(\frac{1}{2}\eta_1))} \leq \frac{1}{2},
$$
then as $j\to \infty$,
\begin{align*}
\Big\|  w^j(\epsilon, t_B, u) - \overline w(\epsilon, t_B,u)\Big\|_{C([0,1]\times B^{\bar N}(\rho')\times B^{\bar n}(\eta_0))} \leq 2\rho' \Big\| W^j - \overline W \Big\|_{C^1(B^{\bar N}(\rho)\times B^{\bar n}(\frac{1}{2}\eta_1))}\to 0.
\end{align*}
\end{proof}

\begin{proof}[\textbf{Proof of Lemma \ref{yangtianxiang}}]
We need to show for any fixed $\{\varsigma_k\}_{k\in \mathbb{N}^{\bar \nu}}\in  \overline{\mathcal{K}}(\bar N, \bar e, \bar a, \bar \nu)$, the operator
$$
T_{j,M} f(x) := \psi(x) \sum_{\substack{k\in \mathbb{N}^{ \bar \nu} \\ |k|\leq M}} \int f(\gamma_{\tilde \tau_j^{-1} \delta_j t}(x)) \delta_0(t_{B^C}) \varsigma_k^{(2^k)}(t_B)\,dt
$$
is bounded on $L^p(\mathbb{R}^n)$ uniformly in $j (\geq N_0),M$.

Let 
\begin{align*}
&e' =\{e_i^\mu: i\in B^C, \mu \in A^C\} = \{e_i^\mu: \bar N < i \leq N, \bar \nu < \mu \leq \nu\}, \quad e_i' = (e_i^{\bar \nu +1}, \ldots, e_i^{\nu}) \text{ for } \bar N < i \leq N,\\
&\tilde e =\{e_i^\mu: i\in B^C, 1\leq \mu \leq \nu\} = \{e_i^\mu: \bar N < i \leq N, 1 \leq \mu \leq \nu\}, \quad \tilde e_i = (e_i^{1}, \ldots, e_i^{\nu}) \text{ for } \bar N < i \leq N.
\end{align*}
Then for each $i\in B^C$, $e_i^\mu \neq 0$ for some $\mu \in A^C$, and for each $\mu \in A^C$, $e_i^\mu \neq 0$ for some $i\in B^C$, i.e., $e'$ satisfies (\ref{e}). If given $t_{B^C} \in \mathbb{R}^{N-\bar N}$ and $\delta \in [0,\infty)^{\nu- \bar \nu}$, we use dilations based on $e'$:
$$
\delta t_{B^C} = (\delta^{e_{\bar N+1}'} t_{\bar N +1}, \ldots, \delta^{e_N'} t_N),
$$
and denote $\varsigma^{(\delta)}(t_{B^C})= \delta^{ e_{\bar N+1}' + \cdots + e_{N}'} \varsigma(\delta t_{B^C})$. If given $t_{B^C} \in \mathbb{R}^{N-\bar N}$ and $\delta \in [0,\infty)^{\nu}$, we use dilations based on $\tilde e$:
$$
\delta t_{B^C} = (\delta^{\tilde e_{\bar N+1}} t_{\bar N +1}, \ldots, \delta^{\tilde e_N} t_N),
$$
and denote $\varsigma^{(\delta)}(t_{B^C})= \delta^{ \tilde e_{\bar N+1} + \cdots + \tilde e_{N}} \varsigma(\delta t_{B^C})$. The dilations based on $e$ and $\bar e$ are as before (see Section \ref{raj} and Subsection \ref{natural2}).
\begin{center}
\begin{tabular}{ |c|c c| } 
 \hline
   & dimension & parameter \\ 
   \hline
 $e$ & $N$ & $\nu$ \\ 
 $\bar e$ & $\bar N$ & $\bar \nu$ \\ 
 $e'$ & $N-\bar N$ & $\nu -\bar \nu$ \\
 $\tilde e$ & $N-\bar N$ & $\nu$\\
 \hline
\end{tabular}
\end{center}
For $1\leq \mu \leq \bar \nu$, let $t_B^\mu$ be the vector consisting of the coordinates $t_i$ of $t$ with $1\leq i \leq \bar N$ and $e_i^\mu \neq 0$, and let $\tilde t_{B}^\mu$ be the vector consisting of the coordinates $t_i$ of $t$ with $1\leq i \leq \bar N$ and $e_i^\mu = 0$. For $1\leq \mu \leq \nu$, let $t_{B^C}^\mu$ be the vector consisting of the coordinates $t_i$ of $t$ with $\bar N< i \leq N$ and $e_i^\mu \neq 0$.

Since $e'$ satisfies (\ref{e}), by Lemma \ref{ff}, the Dirac function $\delta_0(t_{B^C})$ can be written as
$$
\delta_0(t_{B^C}) = \sum_{s\in \mathbb{N}^{\nu - \bar \nu}} \tilde \varsigma_s^{(2^s)}(t_{B^C}),
$$
for some $\{\tilde \varsigma_s\}_{s\in \mathbb{N}^{\nu - \bar \nu}}\in \overline{\mathcal{K}}(N-\bar N, e', \bar a, \nu - \bar \nu)$. Associated to each $\{\varsigma_k\}_{k\in \mathbb{N}^{\bar \nu}} \in \overline{\mathcal{K}}(\bar N, \bar e, \bar a, \bar \nu)$, we define a $\{\bar \varsigma_{(k,s)} \}_{(k,s)\in \mathbb{N}^\nu}$ by letting $\bar \varsigma_{(k,s)}(t) = \varsigma_k(t_B) \cdot \tilde \varsigma_s(t_{B^C})$ for every $(k,s)\in \mathbb{N}^\nu$. Fix arbitrary $(k,s)\in \mathbb{N}^{\nu}$. If $(k,s)_\mu \neq 0$ for some $\mu \in A^C$, then $s_\mu \neq 0$. 
Thus $\int \tilde \varsigma_s(t_{B^C}) \,dt_{B^C}^\mu \equiv 0$, and 
$$
\int \bar \varsigma_{(k,s)}(t) \,dt^\mu \equiv  \varsigma_k(t_B) \int \tilde \varsigma_s(t_{B^C}) \,dt_{B^C}^\mu \equiv 0.
$$
If $(k,s)_\mu \neq 0$ for some $\mu \in A$, then $k_\mu \neq 0$. Thus $\int \varsigma_k(t_B) dt_B^\mu \equiv 0$, and
$$
\int \bar \varsigma_{(k,s)}(t) \,dt^\mu \equiv \int \tilde \varsigma_s(t_{B^C}) \int \varsigma_k(t_B) \,dt_B^\mu \,dt_{B^C}^\mu \equiv 0.
$$
And since $\{\bar \varsigma_{(k,s)} \}_{(k,s)\in \mathbb{N}^\nu} \subseteq C^\infty_c(\overline{B^N(\sqrt{2} \bar a)}) \subseteq C^\infty_c(B^N(a))$ is a bounded set, we have $\{\bar \varsigma_{(k,s)}\}_{(k,s)} \in \overline{\mathcal{K}}(N,e, \sqrt{2} \bar a, \nu)$, $\sum_{(k,s)\in \mathbb{N}^\nu} \bar \varsigma_{(k,s)}^{(2^k,2^s)} \in \mathcal{K}(N,e,a,\nu)$. Then the operator
\begin{align*}
T_{\{\varsigma_k\}_{k\in \mathbb{N}^{\bar \nu}}} f(x) :&= \psi(x) \sum_{(k,s)\in \mathbb{N}^\nu} \int f(\gamma_t(x)) \bar \varsigma_{(k,s)}^{(2^k,2^s)}(t)\,dt\\
& = \psi(x) \sum_{k\in \mathbb{N}^{\bar \nu}} \int f(\gamma_t(x)) \sum_{s\in \mathbb{N}^{\nu- \bar \nu}} 2^{(k,s)\cdot\sum_{i=1}^N e_i}  \varsigma_k(2^k t_B) \tilde \varsigma_s ((2^k,2^s) t_{B^C})\,dt\\
&=  \psi(x) \sum_{k\in \mathbb{N}^{\bar \nu}} \int f(\gamma_t(x)) \varsigma_k^{(2^k)}(t_B) \sum_{s\in \mathbb{N}^{\nu-  \bar \nu}}   2^{(k,0)\cdot\sum_{i=\bar N+1}^N e_i} \tilde \varsigma_s^{(2^s)} ((2^k,1) t_{B^C})\,dt\\
&=  \psi(x) \sum_{k\in \mathbb{N}^{\bar \nu}} \int f(\gamma_t(x)) \varsigma_k^{(2^k)}(t_B)  \delta_0^{(2^k,1)}(t_{B^C})\,dt\\
&=  \psi(x) \sum_{k\in \mathbb{N}^{\bar \nu}} \int f(\gamma_t(x)) \varsigma_k^{(2^k)}(t_B)  \delta_0(t_{B^C})\,dt.
\end{align*}
is bounded on $L^p(\mathbb{R}^n)$.

Then we have the linear map
\begin{align*}
\overline{\mathcal{K}}( \bar N, \bar e,\bar a,\bar \nu) \to \mathcal{B}(L^p(\mathbb{R}^n)), \quad \{\varsigma_k\}_{k\in \mathbb{N}^{\bar \nu}}  &\mapsto T_{\{\varsigma_k\}_{k\in \mathbb{N}^{\bar \nu}}}.
\end{align*}
Note this map is not one-to-one. $\overline{\mathcal{K}}( \bar N, \bar e,\bar a,\bar \nu)$ and $\mathcal{B}(L^p(\mathbb{R}^n))$ are both Fr\'echet spaces. If $\{\varsigma_k\}  \mapsto T_{\{\varsigma_k\}}$ is a closed map, then by the closed graph theorem, $\{\varsigma_k\}  \mapsto T_{\{\varsigma_k\}}$ is a bounded map.

Now we show $\{\varsigma_k\}  \mapsto T_{\{\varsigma_k\}}$ is a closed map. Suppose $\{\varsigma_k^j\}$ converges in $\overline{\mathcal{K}}( \bar N, \bar e,\bar a,\bar \nu)$ to $\{\varsigma_k^0\}$ as $j\to \infty$, and $T_{\{\varsigma_k^j\}} \to T_0$ in $\mathcal{B}(L^p(\mathbb{R}^n))$ as $j\to \infty$. We need to show $T_{\{\varsigma_k^0\}} = T_0$. Fix arbitrary $f\in C_c^\infty(\mathbb{R}^n)$. We first show $T_{\{\varsigma_k^j - \varsigma_k^0\}}f\to 0$ pointwise as $j\to \infty$. For each $k\in \mathbb{N}^{\bar \nu}$, choose $1\leq \mu(k) \leq \bar \nu$ such that $k_{\mu(k)} = |k|_\infty$. Fix arbitrary $x\in \text{~supp}\,\psi$. Then $f(\gamma_{(t_B,0)}(x))$ is a smooth function of $t_B$ on $B^{\bar N}(a)$, and its $C^1(B^{\bar N}(\bar a))$ norm is bounded by some constant $C_{f,x}$. Let 
$$
c'=\min \{e_i^\mu: 1\leq i \leq \bar N, 1\leq \mu \leq \bar \nu, e_i^\mu > 0\}, \quad  c= \frac{c'}{\sqrt{\bar \nu}}>0.
$$
Recall $t_B^{\mu(k)}$ is the vector consisting of the coordinates $t_i$ of $t$ with $1\leq i \leq \bar N$ and $e_i^{\mu(k)} \neq 0$, and $\tilde t_B^{\mu(k)}$ is the vector consisting of the coordinates $t_i$ of $t$ with $1\leq i \leq \bar N$ and $e_i^{\mu(k)} =0$.
Then 
\begin{align*}
|T_{\{\varsigma_k^j-\varsigma_k^0\}}f(x)| &\leq |\psi(x)| \sum_{k\in \mathbb{N}^{\bar \nu}} \Big| \int f(\gamma_{(t_B,0)}(x)) \big((\varsigma_k^{j})^{(2^k)}(t_B)-(\varsigma_k^{0})^{(2^k)}(t_B)\big)\,dt_B\Big|\\
&= \Big|\psi(x) \int f(\gamma_{(t_B,0)}(x)) (\varsigma_0^j(t_B) - \varsigma_0^0(t_B))\,dt_B \Big|\\
&\quad + |\psi(x)| \sum_{\substack{k\in \mathbb{N}^{\bar \nu} \\ k\neq 0}} \Big|\int f(\gamma_{(t_B^{\mu(k)}, \tilde t_{B}^{\mu(k)},0)}(x)) \big((\varsigma_k^{j})^{(2^k)}(t_B)-(\varsigma_k^{0})^{(2^k)}(t_B)\big)\,dt_B \\
&\qquad \qquad \qquad \quad - \int f(\gamma_{(0, \tilde t_{B}^{\mu(k)},0)}(x)) \big((\varsigma_k^{j})^{(2^k)}(t_B)-(\varsigma_k^{0})^{(2^k)}(t_B)\big)\,dt_B\Big| \\
&\leq C_{f,x}  |\psi(x)| \Big(\int \big|\varsigma_0^j(t_B) -\varsigma_0^0(t_B)\big| \,dt_B + \sum_{\substack{k\in \mathbb{N}^{\bar \nu} \\ k\neq 0}} \bar a 2^{-c' k_{\mu(k)}} \int \big|(\varsigma_k^{j})^{(2^k)}(t_B)-(\varsigma_k^{0})^{(2^k)}(t_B)\big| \,dt_B \Big)\\
&\lesssim  C_{f,x}  |\psi(x)| (\bar a^{\bar N} + \bar a^{\bar N+1} \sum_{\substack{k\in \mathbb{N}^{\bar \nu} \\ k\neq 0}} 2^{-c | k|}  \big) \cdot \big\|\{\varsigma_k^j - \varsigma_k^0\big\}_k\|_{0},
\end{align*}
where $\|\cdot \|_0$ is a semi-norm defined in Definition \ref{semi}.
Hence $T_{\{\varsigma_k^j\}} f\to T_{\{\varsigma_k^0\}}f$ pointwise as $j\to \infty$.

Since $T_{\{\varsigma_k^j\}} \to T_0$ in $\mathcal{B}(L^p(\mathbb{R}^n))$ as $j\to \infty$, $T_{\{\varsigma_k^j\}} f$ converges to $T_0 f$ in $L^p(\mathbb{R}^n)$ and thus converges in measure. So there exists a subsequence $T_{\{\varsigma_k^{j_l}\}}f$ convergent to $T_0f$ almost everywhere as $l\to \infty$. Then $T_{\{\varsigma_k^0\}}f(x) = T_0f(x)$ a.e.. Since $C_c^\infty(\mathbb{R}^n)$ is dense in $L^p(\mathbb{R}^n)\, (1<p<\infty)$, we have $T_{\{\varsigma_k^0\}} =T_0$. Therefore $\{\varsigma_k\} \mapsto T_{\{\varsigma_k\} }$ is a closed map, and thus $\{\varsigma_k\} \mapsto T_{\{\varsigma_k\} }$ is a bounded map.

Fix arbitrary
$\{\varsigma_k\}_{k\in \mathbb{N}^{\bar \nu}} \in \overline{\mathcal{K}}(\bar N, \bar e, \bar a, \bar \nu) $. Recall $\tilde \tau_j = (\tau_j^{b_1}, \ldots, \tau_j^{b_{ \nu}})$, $ \delta_j = (\delta_j^1, \ldots, \delta_j^{ \nu})$. Let $\bar \tau_j=(\tau_j^{b_1}, \ldots, \tau_j^{b_{\bar \nu}})$, $\bar \delta_j = (\delta_j^1, \ldots, \delta_j^{\bar \nu})$. Then $\bar \tau_j, \bar \delta_j, \bar \delta_j^{-1} \in (0,\infty)^{\bar \nu}$, and 
$$
T_{j,M}f(x) = \psi(x) \sum_{\substack{k\in \mathbb{N}^{\bar \nu} \\ |k|\leq M}} \int f(\gamma_t(x)) \delta_0(t_{B^C}) \varsigma_k^{(2^k \bar \tau_j \bar \delta_j^{-1})}(t_B)\,dt = T_{\big\{\varsigma_k^{(\bar \tau_j \bar \delta_j^{-1})}\big\}} f(x).
$$

\begin{example}
Before we start the technical proof that the kernel $\sum_k \varsigma_k^{(2^k \bar \tau_j \bar \delta_j^{-1})}$ corresponds to an element in $\overline{\mathcal{K}}(\bar N, \bar e, \bar a, \bar \nu)$ bounded uniformly in $j$, we illustrate the basic idea by a $1$-parameter example. 
Let $\{\varsigma_k(t)\}_{k\in \mathbb{N}}$ be a bounded sequence in $C_c^\infty(\overline{B^1(\bar a)})$, satisfying $\int \varsigma_k(t)\,dt =0$ for $k\neq 0$, i.e., $\{\varsigma_k\}_k \in \overline{\mathcal{K}}(1, \{1\}, \bar a, 1)$. Let $m_j \in \mathbb{N}_{>0}$ be such that $m_j \to \infty$ as $j\to \infty$, and let 
$$
K_j(t)=\sum_{k\in \mathbb{N}} \varsigma_k^{(2^{k+m_j+1})}(t) = \sum_{k\in \mathbb{N}} 2^{k+m_j +1} \varsigma_k(2^{k+m_j +1} t).
$$
The goal is to decompose the kernel $K_j(t)$ into $\sum_{l\in \mathbb{N}} \phi_{j,l}^{(2^l)}(t)$, such that $\big\{ \{\phi_{j,l}\}_l: j\in \mathbb{N} \big\}$ is bounded in $\overline{\mathcal{K}}(1, \{1\}, \bar a, 1)$. If we let $\phi_{j,l}(t) = \varsigma_l^{(2^{m_j +1})}(t)$ for every $l$, $\{\varsigma_l^{(2^{m_j+1})}(t)\}_{j,l\in \mathbb{N}}$ is not bounded in $C_c^\infty(\overline{B^1(\bar a)})$; if we let $\phi_{j,l}(t) = \varsigma_{l-m_j}^{(2)} (t)$ for every $l\geq m_j$ and $\phi_{j,l} \equiv 0$ for every $l<m_j$, $\int \phi_{j, m_j}(t)\,dt = \int 2 \varsigma_0(2t)\,dt$ does not necessarily vanish, i.e., does not necessarily satisfy the cancellation condition. These issues are resolved by letting
$$
\phi_{j,l}(t)=
\left\{
\begin{aligned}
&2 \varsigma_{l-m_j}(2t), & l>m_j,\\
&2\varsigma_0(2t) -\varsigma_0(t), & 0<l \leq m_j,\\
&2\varsigma_0(2t), & l=0.
\end{aligned}
\right.
$$
Note $\{\phi_{j,l}\}_{j,l}$ is bounded in $C_c^\infty(\overline{B^1(\bar a)})$, and $\int \phi_{j,l}(t)\,dt =0$ for $l\neq 0$. Thus $\big\{ \{\phi_{j,l}\}_l: j\in \mathbb{N} \big\}$ is bounded in $\overline{\mathcal{K}}(1, \{1\}, \bar a, 1)$. We have
\begin{align*}
\sum_{l\in \mathbb{N}} \phi_{j,l}^{(2^l)}(t) &= 2\varsigma_0 (2t) + \sum_{0<l \leq m_j} \big( \varsigma_0^{(2^{l+1})}(t) - \varsigma_0^{(2^{l})}(t) \big) + \sum_{l>m_j} \varsigma_{l-m_j}^{(2^{l+1})} (t)\\
&=  \varsigma_0^{(2^{m_j+1})}(t) + \sum_{k>0}  \varsigma_k^{(2^{m_j +k+1})}(t)=K_j(t).
\end{align*}
\end{example}

Now we come back to the general multi-parameter case. Recall for $j\geq N_0$, $\tilde \tau_j^{-1} \delta_j \in [0,\frac{1}{2}]^\nu$, thus for each $j \geq N_0$, there exists $m_j=(m_j^1, \ldots, m_j^{\bar \nu}) \in \mathbb{N}^{\bar \nu}$, such that for $1\leq \mu \leq \bar \nu$,
$$
2^{m_j^\mu+1} \leq (\tau_j)^{b_\mu} (\delta_j^\mu)^{-1} < 2^{m_j^\mu +2}.
$$
Hence $\big\{\varsigma_k^{(2^{-p-m_j} \bar \tau_j \bar \delta_j^{-1})}: j\geq N_0, k\in \mathbb{N}^{\bar \nu}, p\in \{0,1\}^{\bar \nu}\big\}$ is bounded in $C_c^\infty(\overline{B^{\bar N}(\bar a)})$. 
For any $a,b\in \mathbb{Z}^{\bar \nu}$, let $a_\mu$ denote the $\mu$-th coordinate of $a$, let $a_+$ denote the vector obtained by taking the positive part of each coordinate of $a$, and let $a \vee b$ denote the vector obtained by taking the larger one of each coordinate of $a$ and $b$.
For every $j \geq N_0$, $l\in \mathbb{N}^{\bar \nu}$, define
\begin{equation}\label{redecompose}
\tilde \phi_{j,l}(t_B) =  \sum_{\substack{p\in \{0,1\}^{\bar \nu} \\ l-p\in \mathbb{N}^{\bar \nu} \\ p\leq (m_j +1-l)_+}} (-1)^{|p|_1}\, \varsigma_{(l-m_j)_+}^{(2^{-p-m_j} \bar \tau_j \bar \delta_j^{-1})}(t_B),
\end{equation}
where $|p|_1$ denotes the $l^1$ norm of the vector $p$. Then $\{\tilde \phi_{j,l} \}_{\substack{j\geq N_0\\ l\in \mathbb{N}^{\bar \nu}}}$ is bounded in $C_c^\infty(\overline{B^{\bar N}(\bar a)})$.

Fix $j\geq N_0$, $l\in \mathbb{N}^{\bar \nu}$ and $1\leq \mu \leq \bar \nu$ such that $l_\mu \neq 0$. We wish to show
\begin{equation}\label{cancellation}
\int \tilde \phi_{j,l} (t_B)\,dt_B^\mu \equiv 0.
\end{equation}
If $l_\mu > (m_j)\mu$, then $\big( (l-m_j)_+ \big)_\mu >0$, and thus $\int \varsigma_{(l-m_j)_+}(t_B)\,dt_B^\mu \equiv 0$, which implies (\ref{cancellation}). If $l_\mu \leq (m_j)_\mu$, then $\big((m_j+1-l)_+\big)_\mu \geq 1$, and since $l_\mu \geq 1$, the sum in (\ref{redecompose}) contains two terms: one for $p_\mu =1$ and one for $p_\mu =0$. These two terms have opposite signs, and therefore the sum of their integrals against $t_B^\mu$ is $0$. (\ref{cancellation}) follows.

For every $j\geq N_0, M\in \mathbb{N}, l\in \mathbb{N}^{\bar \nu}$, define
$$
\phi_{j,M,l} =
\left\{
\begin{aligned}
\tilde \phi_{j,l}, \quad |(l-m_j)_+|\leq M,\\
0, \quad |(l-m_j)_+|> M.
\end{aligned}
\right.
$$
Then
$$
\Big\{ \{ \phi_{j,M,l}\}_{l\in \mathbb{N}^{\bar \nu}}: j\geq N_0, M\in \mathbb{N} \Big\}
$$ 
is bounded in $\overline{\mathcal{K}}(\bar N, \bar e, \bar a, \bar \nu)$.

Fix $j \geq N_0, M\in \mathbb{N}$. We have
\begin{align*}
\sum_{l\in \mathbb{N}^{\bar \nu}} \phi_{j,M,l}^{(2^l)} &=\sum_{l:\,|(l-m_j)_+|\leq M} \tilde \phi_{j,l}^{(2^l)} =\sum_{l:\,|(l-m_j)_+|\leq M}\, \sum_{\substack{p:\,p\in \{0,1\}^{\bar \nu} \\ l-p\in \mathbb{N}^{\bar \nu} \\ p\leq (m_j +1-l)_+}} (-1)^{|p|_1}\, \varsigma_{(l-m_j)_+}^{(2^{l-p-m_j} \bar \tau_j \bar \delta_j^{-1})}\\
&= \sum_{\substack{k:\,k- m_j \in \mathbb{N}^{\bar \nu}\\ |k-m_j|\leq M}} \,\, \sum_{l:\,l \vee m_j = k} \, \sum_{\substack{p:\,p\in \{0,1\}^{\bar \nu} \\ l-p\in \mathbb{N}^{\bar \nu} \\ p\leq (m_j +1-l)_+}} (-1)^{|p|_1} \big( \varsigma_{k-m_j}^{(2^{-m_j} \bar \tau_j \bar \delta_j^{-1})} \big)^{(2^{l-p})} \\
&= \sum_{\substack{k:\,k- m_j \in \mathbb{N}^{\bar \nu} \\ |k-m_j| \leq M}} \varsigma_{k-m_j}^{(2^{k-m_j} \bar \tau_j \bar \delta_j^{-1})},
\end{align*}
where in the second line we let $k := (l-m_j)_+ + m_j = l\vee m_j$ and classify the $l$ according to the value of $k$, which results in a telescoping sum over $l,p$, and the last equality is verified by counting how many times each term $(\varsigma_{k-m_j}^{(2^{-m_j}\bar \tau_j \bar \delta_j^{-1})} )^{(2^{l'})}$ appears in the telescoping sum and with what sign. In other words, for fixed $k$, the $l':=l- p \leq k$ in the telescoping sum vary only in coordinates $\{1\leq \mu \leq \bar \nu: l_\mu \leq (m_j)_\mu\}$, and if $l'$ has at least one coordinate less than that of $k$, by denoting $\nu_0$ the number of coordinates of $l'$ less than that of $k$, the coefficient of $\varsigma_{k-m_j}^{(2^{l'-m_j} \bar \tau_j \bar \delta_j^{-1})}$ in the telescoping sum is equal to $$
\sum_{\mu=0}^{\nu_0} (-1)^\mu \binom{\nu_0}{\mu} =0,
$$
whereas the coefficient of $\varsigma_{k-m_j}^{(2^{k-m_j} \bar \tau_j \bar \delta_j^{-1})}$ is $1$.

Therefore the corresponding operator
\begin{align*}
\psi(x) \sum_{\substack{k:\, k-m_j \in \mathbb{N}^{\bar \nu}\\ |k-m_j|\leq M}} \int f(\gamma_t(x)) \delta_0(t_{B^C}) \varsigma_{k-m_j}^{(2^{k-m_j} \bar \tau_j \bar \delta_j^{-1})} (t_B)\,dt &= \psi(x) \sum_{\substack{k:\, k-m_j \in \mathbb{N}^{\bar \nu}\\ |k-m_j|\leq M}} \int f(\gamma_{\tilde \tau_j^{-1} \delta_j t}(x)) \delta_0(t_{B^C}) \varsigma_{k-m_j}^{(2^{k-m_j})} (t_B)\,dt\\
&= \psi(x) \sum_{\substack{ k \in \mathbb{N}^{\bar \nu}\\ |k|\leq M}} \int f(\gamma_{\tilde \tau_j^{-1} \delta_j t}(x)) \delta_0(t_{B^C}) \varsigma_{k}^{(2^{k})} (t_B)\,dt
\end{align*}
is bounded on $L^p(\mathbb{R}^n)$ uniformly in $j(\geq N_0),M$.
\end{proof}


\begin{proof}[\textbf{Proof of Lemma \ref{modu}}]
We need to show there exists $N_2 \geq N_1$ such that for $j \geq N_2$, there exists $0< \sigma_j\leq \bar \sigma$, such that $E_j$ is a diffeomorphism from $B^{\bar n}(\frac{3}{4}\eta_1) \times B^{n-\bar n}(\sigma_j)$ onto its image, and that for $(u,v)\in B^{\bar n}(\frac{3}{4}\eta_1) \times B^{n-\bar n}(\sigma_j)$,
$$
|\det DE_j(u,v)| \approx |\det DE_j(0,0)|>0.
$$

Recall the normal space to the submanifold $\Phi_j(B^{\bar n}(\eta_1))$ at $x_j$ has an orthonormal basis $\{V_{\bar n+1}, \ldots, V_n\}$. By a rotation of the coordinate axes of $\mathbb{R}^n$, we can assume
$$
V_i = (\underbrace{0, \ldots, 0}_{i-1}, 1, \underbrace{0, \ldots, 0}_{n-i}) \text{ for } \bar n+1 \leq i \leq n.
$$
Let $\mathbf{V}$ be the $n\times (n-\bar n)$ matrix with column vectors $V_{\bar n+1} ,\ldots, V_n$, i.e.
$$
\begin{bmatrix}
0 & 0 &\cdots & 0\\
\cdots & \cdots & \cdots & \cdots\\
0 & 0 &\cdots & 0\\
1 & 0 &\cdots & 0\\
0 & 1 & \cdots & 0\\
\cdots & \cdots & \cdots & \cdots\\
0 & 0 &\cdots & 1
\end{bmatrix}.
$$
Note $\mathbf{V}$ and these rotations of $\mathbb{R}^n$ depend on $j$, but since rotations have Jacobian equal to $1$, the following estimates are still uniform in $j$.

We use notations in Definition \ref{indoor}. By Theorem \ref{16}, for every $j$, $\Phi_j$ is a diffeomorphism from $B^{\bar n}(\eta_1)$ onto its image, and on $B^{\bar n}(\eta_1)$,
\begin{align*}
\big| \det \mathbf{Y}^j_{J} \big| \approx 1,\quad \big| \det_{\bar n \times \bar n} D\Phi_j \big| \approx \big| \det_{\bar n \times \bar n} D\Phi_j(0) \big| >0.
\end{align*}
There exists $I_j \in \mathcal{I}(\bar n, n)$ such that
$$
\Big| \det (\tau_j^{-1} \delta_j \mathbf{X})_{I_j \times J} (x_j) \Big| = \Big| \det_{\bar n \times \bar n} (\tau_j^{-1} \delta_j \mathbf{X})_{J} (x_j) \Big|_\infty = \Big| \det_{\bar n \times \bar n} (\tau_j^{-1} \delta_j \mathbf{X}) (x_j) \Big|_\infty>0.
$$

Note for every $j$, $\tau_j^{-1} \delta_j \mathbf{X} = \{ \tau_j^{-\theta_l} \delta_j^{d_l} X_l\}_{l=1}^Q$ satisfies (\ref{commutator}) and (\ref{coefficient}). Then even if we replace $\mathbf{X}$ by $\tau_j^{-1} \delta_j \mathbf{X}$, the same set of $m$-admissible constants for $m\in \mathbb{N}$ still works. Thus we can apply Lemma \ref{7} with $\mathbf{X}$ replaced by $\tau_j^{-1} \delta_j \mathbf{X}$, with $x=x_j, \delta=(1,\ldots,1)$, $\epsilon =1$, (note $\widetilde{\xi_1}(1) = \xi_1$). Then for every $j$ and every $u\in B^{\bar n}(\eta_1)$, by Lemma \ref{7} (ii),
\begin{align*}
\big| \det (\tau_j^{-1} \delta_j \mathbf{X})_{I_j \times J} \big( \Phi_j(u) \big) \big| &\approx \big| \det_{\bar n\times \bar n} (\tau_j^{-1} \delta_j \mathbf{X}) \big(\Phi_j(u)\big) \big| \approx \big| \det_{\bar n\times \bar n} (\tau_j^{-1} \delta_j \mathbf{X})_{J} \big( \Phi_j(u) \big) \big|,
\end{align*}
and by Lemma \ref{7} (i),
\begin{align*}
\big| \det (\tau_j^{-1} \delta_j \mathbf{X})_{I_j \times J} \big( \Phi_j(u) \big) \big| &\approx \big| \det_{\bar n \times \bar n} (\tau_j^{-1} \delta_j \mathbf{X}) \big(\Phi_j(u)\big) \big| \approx \big| \det_{\bar n\times \bar n} (\tau_j^{-1} \delta_j \mathbf{X}) (x_j) \big| \approx \big| \det (\tau_j^{-1} \delta_j \mathbf{X})_{I_j \times J} (x_j) \big|.
\end{align*}
Since the tangent vectors $(\tau_j^{-1} \delta_j \mathbf{X})_J(x_j)$ to the submanifold $\Phi_j(B^{\bar n}(\eta_1))$ at $x_j$ are orthogonal to $\mathbf{V}$, we have $I_j = \{1, \ldots, \bar n\}$, denoted as $I$. Without loss of generality, by reordering $(\mathbf{X}, \mathbf{d})$, we can assume $J = (1, \ldots, \bar n)$.

We first show that for $j \geq N_1$, there exists $0< \tilde \sigma_j\leq \bar \sigma$, for $(u,v)\in B^{\bar n}(\frac{3}{4}\eta_1) \times B^{n-\bar n}(\tilde \sigma_j)$,
\begin{equation}\label{mujianping}
|\det DE_j(u,v)| \approx |\det DE_j(0,0)|>0.
\end{equation}
By compactness of $\overline{B^{\bar n}(\frac{3}{4}\eta_1)}\times \{0\}$ and continuity of $\big| \det DE_j\big|$, it suffices to show that (\ref{mujianping}) holds on $\overline{B^{\bar n}(\frac{3}{4}\eta_1)}\times \{0\}$.

Recall for $j \geq N_1$, 
\begin{equation}\label{apple}
\|H_j\|_{C^2(B^{\bar n}(\eta_1)\times B^{n-\bar n}(\bar \sigma))} \leq C'.
\end{equation}
Denote $\tau_j^{-1} \delta_j u := \{\tau_j^{-\theta_l} \delta_j^{d_l} u_l\}_{l\in J} = (\tau_j^{-\theta_1} \delta_j^{d_1}u_1, \ldots, \tau_j^{-\theta_{\bar n}} \delta_j^{d_{\bar n}}u_{\bar n})$, and $\tau_j^{-1} \delta_j(u,v) := (\tau_j^{-1} \delta_j u, v)$. Then $E_j(u,v) = H_j(\tau_j^{-1} \delta_j (u, v))$. Denote the Frobenius norm for a matrix $A=(a_{i,j})$ as $\|A\|_*=\big(\sum_{i,j}|a_{i,j}|^2 \big)^{\frac{1}{2}}$. For $j \geq N_1$ and $u\in B^{\bar n}(\eta_1)$,
\begin{equation}\label{over}
\begin{aligned}
\big\| D_v E_j(u,0) -\mathbf{V} \big\|_*  &= \big\| D_v E_j(u,0) - D_v E_j(0,0) \big\|_* = \big\| D_v H_j(\tau_j^{-1} \delta_j u, 0) - D_vH_j (0,0) \big\|_*\\
& \leq C' \sup_{u\in B^{\bar n}(\eta_1)} |\tau_j^{-1} \delta_j u| \lesssim C' \max_{1\leq i \leq \bar n} \tau_j^{-\theta_i} \delta_j^{d_i} \to 0, \quad \text{as~} j\to \infty.
\end{aligned}
\end{equation}

For every $j$ and every $u\in B^{\bar n}(\eta_1)$, $(\tau_j^{-1} \delta_j \mathbf{X})_{J}(\Phi_j(u)) = D\Phi_j(u) \cdot \mathbf{Y}_{J}^j(u)$, thus for every $I'\in \mathcal{I}(\bar n, n)$, 
$$
\big| \det (\tau_j^{-1} \delta_j \mathbf{X})_{I'\times J}(\Phi_j(u)) \big| = \big| \det (D\Phi_j)_{I'} (u) \big| \cdot \big| \det \mathbf{Y}^j_{J} (u)\big| \approx \big| \det (D\Phi_j)_{I'} (u) \big|,
$$
and hence
\begin{equation}\label{hope}
\begin{aligned}
\big| \det  (D\Phi_j)_{I} (u) \big| &\approx \big| \det (\tau_j^{-1} \delta_j \mathbf{X})_{I\times J} (\Phi_j(u)) \big| \approx \big| \det_{\bar n\times \bar n} (\tau_j^{-1} \delta_j \mathbf{X})_{J} (\Phi_j(u)) \big|\\
& \approx \big| \det_{\bar n\times \bar n} D\Phi_j(u)\big| \approx \big| \det (\tau_j^{-1} \delta_j \mathbf{X})_{I\times J} (x_j) \big| \approx \big| \det (D\Phi_j)_{I}(0) \big|>0.
\end{aligned}
\end{equation}
By (\ref{over}) and (\ref{hope}), there exists $N_2 \geq N_1$ such that for $j \geq N_2$ and $u\in B^{\bar n}(\eta_1)$,
\begin{equation}\label{ryan}
\begin{aligned}
\big| \det DE_j(u,0) \big| &= \big| \det  \big[ D\Phi_j(u), D_v E_j(u,0)\big] \big| \approx \big| \det \big[ D\Phi_j(u), \mathbf{V}\big] \big| = \big| \det (D\Phi_j)_{I}(u)\big|\\
&\approx \big| \det (D\Phi_j)_{I}(0)\big| = \big| \det \big[ D\Phi_j(0), \mathbf{V}\big] \big| = \big| \det DE_j (0,0) \big|>0.
\end{aligned}
\end{equation}
Hence for $j \geq N_2$, there exists $0< \tilde \sigma_j\leq \bar \sigma$, such that for $(u,v)\in B^{\bar n}(\frac{3}{4}\eta_1) \times B^{n-\bar n}(\tilde \sigma_j)$, $|\det DE_j(u,v)| \approx |\det DE_j(0,0)|>0$.

Next we show each $E_j (j\geq N_2)$ is a diffeomorphism from $B^{\bar n}(\frac{3}{4}\eta_1)\times B^{n-\bar n}(\sigma_j)$ onto its image, for some $0<\sigma_j \leq \tilde \sigma_j$. Note $E_j\big|_{\overline{B^{\bar n}(\frac{3}{4}\eta_1)}\times \{0\}}$ is a diffeomorphism onto its image. Thus $E_j\big|_{\overline{B^{\bar n}(\frac{3}{4}\eta_1)}\times \{0\}}$ is injective. By compactness of $\overline{B^{\bar n}(\frac{3}{4}\eta_1)}\times \{0\}$, we have
$$
\overline{B^{\bar n}(\frac{3}{4}\eta_1)} \times \{0\}= \bigcap_{\epsilon>0} \Big(\overline{B^{\bar n}(\frac{3}{4}\eta_1)}\times B^{n-\bar n}(\epsilon)\Big).
$$
Suppose by contradiction that for every $\epsilon>0$, $E_j$ is not injective on $\overline{B^{\bar n}(\frac{3}{4}\eta_1)}\times B^{n-\bar n}(\epsilon)$. Then there exist sequences $\epsilon_k \to 0$, $a_k \neq b_k \in \overline{B^{\bar n}(\frac{3}{4}\eta_1)}\times B^{n-\bar n}(\epsilon_k)$ such that $E_j(a_k) = E_j(b_k)$. Bounded sequences $a_k, b_k$ have convergent subsequences $a_{k_l}\to a_0, b_{k_l}\to b_0$ for some $a_0, b_0 \in \overline{B^{\bar n}(\frac{3}{4}\eta_1)}\times \{0\}$. Thus $E_j(a_0) = E_j(b_0)$. But $E_j$ is injective on $\overline{B^{\bar n}(\frac{3}{4}\eta_1)} \times \{0\}$. So $a_0= b_0$. Then $E_j$ is not injective on any neighborhood of $a_0$, which contradicts the fact that $DE_j(a_0)$ is invertible. Hence there exists $0<\sigma_j \leq \tilde \sigma_j$, such that $E_j$ is injective on $\overline{B^{\bar n}(\frac{3}{4}\eta_1)}\times B^{n-\bar n}(\sigma_j)$, and thus $E_j$ is a diffeomorphism from $B^{\bar n}(\frac{3}{4}\eta_1)\times B^{n-\bar n}(\sigma_j)$ onto its image for $j\geq N_2$.
\end{proof}

\begin{proof}[\textbf{Proof of Lemma \ref{strada}}]
We need to show there exists $C_0\geq 1$, such that for $j \geq N_2$, there exists $0<\epsilon_j \leq \frac{1}{2C_0} \sigma_j$, such that for $0<\epsilon \leq \epsilon_j$ and $t\in B^N(\rho')$, $E_j^{-1} \circ \gamma_{\tilde \tau_j^{-1}\delta_j t} \circ E_j$ maps from $B^{\bar n}(\eta_0)$ $\times B^{n-\bar n}(\epsilon)$ into $B^{\bar n}(\frac{3}{4}\eta_1)\times B^{n-\bar n}(C_0\epsilon)$.

By Corollary \ref{VV}, for every $j \geq N_2$, $\gamma_{\tilde \tau_j^{-1} \delta_j t} \circ \Phi_j(u) \in \Phi_j(B^{\bar n}(\frac{1}{2}\eta_1))$ for $u\in B^{\bar n}(\eta_0)$ and $t\in B^N(\rho')$. Fix arbitrary $j \geq N_2$. Since functions $\gamma_t(x), E_j(u, v)$ have finite $C^1$ norms, as $v\to 0$ we have
\begin{equation}\label{sally}
\gamma_{\tilde \tau_j^{-1} \delta_j t} \circ E_j(u,v) \to \gamma_{\tilde \tau_j^{-1} \delta_j t} \circ E_j(u,0) \in E_j(B^{\bar n}(\frac{1}{2}\eta_1)\times \{0\}), \text{ uniformly in } t\in B^N(\rho'), u\in B^{\bar n}(\eta_0).
\end{equation}
So there exists $0<\bar \epsilon_j \leq \frac{1}{2C_0}\sigma_j$, such that $E_j^{-1}\circ \gamma_{\tilde \tau_j^{-1} \delta_j t}\circ E_j (u,v)$ is well-defined for $t\in B^N(\rho'), (u,v)\in B^{\bar n}(\eta_0)\times B^{n-\bar n}(\bar \epsilon_j)$, and takes values in $B^{\bar n}(\frac{3}{4}\eta_1) \times B^{n-\bar n}(\sigma_j)$.

As in the proof of the previous lemma, we assume $I=J=\{1, \ldots, \bar n\}$. 
Denote $I_0 = \{1, \ldots, n\}$, $I^C =J^C= \{\bar n +1, \ldots, n\}$. It suffices to show there exists $C_0\geq 1$, such that for $j \geq N_2$, there exists $0<\epsilon_j \leq \bar \epsilon_j$, such that for $0< \epsilon \leq \epsilon_j$ and $t\in B^N(\rho')$, on $B^{\bar n}(\eta_0)\times B^{n-\bar n}(\epsilon)$, 
$$
\big| \big( E_j^{-1} \circ \gamma_{\tilde \tau_j^{-1} \delta_j t} \circ E_j (u,v) \big)_{I^C} \big| \leq C_0 \epsilon.
$$

For $j \geq N_2$, let $\epsilon_j\in (0,\bar \epsilon_j]$ be a small number, to be chosen later. For $0<\epsilon \leq \epsilon_j, t\in B^N(\rho')$, and $(u,v) \in B^{\bar n}(\eta_0)\times B^{n-\bar n}(\epsilon)$,
\begin{align*}
\big| \big( E_j^{-1} \circ \gamma_{\tilde \tau_j^{-1} \delta_j t} \circ E_j (u,v) \big)_{I^C} \big|&= \big| \big( E_j^{-1} \circ \gamma_{\tilde \tau_j^{-1} \delta_j t} \circ E_j (u,v) \big)_{I^C} - \big( E_j^{-1} \circ \gamma_{\tilde \tau_j^{-1} \delta_j t} \circ E_j (u,0) \big)_{I^C} \big|\\
&\leq \sup_{(u,v)\in B^{\bar n}(\eta_0)\times B^{n-\bar n}(\epsilon)} \big\| D_v \big( (E_j^{-1}\circ \gamma_{\tilde \tau_j^{-1} \delta_j t} \circ E_j )_{I^C}\big)(u,v) \big\|_* \cdot \epsilon\\
&= \sup_{B^{\bar n}(\eta_0)\times B^{n-\bar n}(\epsilon)} \big\|  \big( D(E_j^{-1}\circ \gamma_{\tilde \tau_j^{-1} \delta_j t} \circ E_j) \big)_{I^C\times J^C}(u,v) \big\|_* \cdot \epsilon.
\end{align*}
It suffices to show that there exists $C_0\geq 1$, for $j\geq N_2$,
$$
\sup_{B^{\bar n}(\eta_0)\times B^{n-\bar n}(\epsilon_j)} \big\|  \big( D(E_j^{-1}\circ \gamma_{\tilde \tau_j^{-1} \delta_j t} \circ E_j) \big)_{I^C\times J^C}(u,v) \big\|_* \leq C_0.
$$

By (\ref{apple}), $\sup_{B^{\bar n}(\eta_1)\times B^{n-\bar n}(\bar \sigma)} \big\| DH_j\big\|_* \lesssim C'$ for $j\geq N_2$. So for $j \geq N_2$, $t\in B^N(\rho')$,
\begin{align*}
&\big\| \big( DH_j(\tau_j^{-1} \delta_j(u,v))\big)_{I_0\times J^C}\big\|_* \lesssim C', \text{~for~} (u,v)\in B^{\bar n}(\eta_1)\times B^{n-\bar n}(\bar \epsilon_j),\\
&\big\|D\gamma_{\tilde \tau_j^{-1} \delta_j t} \big\|_* \leq \| \gamma_t(x) \|_{C^1(B^N(\rho)\times \Omega)}, \text{~on~} E_j(B^{\bar n}(\eta_1) \times B^{n-\bar n}(\bar \epsilon_j)),\\
&\big\| \big( DH_j(0,0)^{-1}\big)_{I^C\times I_0} \big\|_* = \left\|\begin{bmatrix}
\big(\mathbf{X}_{I\times J}(x_j)\big)^{-1} & 0\\
0 & I_{n-\bar n}
\end{bmatrix}_{I^C\times I_0}\right\|_* =\big\| \big[0,I_{n-\bar n}\big]\big\|_* \lesssim 1.
\end{align*}

For an $n\times n$ matrix $A$, $S_j A$ denotes the matrix obtained by dilating the $i$-th column vector of $A$ by $\tau_j^{-\theta_i} \delta_j^{d_i}$ for $1\leq i \leq \bar n$, and $S^*_j A$ denotes the matrix obtained by dilating the $i$-th row of $A$ by $\tau_j^{\theta_i} \delta_j^{-d_i}$ for $1\leq i \leq \bar n$. Then
$$
DE_j(u,v) = S_j \big( DH_j(\tau_j^{-1} \delta_j(u,v))\big), \quad \big(DE_j(u,v) \big)^{-1} = S_j^*\big(DH_j(\tau_j^{-1} \delta_j(u,v))^{-1}\big).
$$
Denote $i^C = I_0\backslash \{i\}, j^C = I_0\backslash \{j\}$. By (\ref{hope}) and (\ref{ryan}), for $j \geq N_2$ and $u\in B^{\bar n}(\eta_1)$,
\begin{align*}
& \quad \big| \det DH_j(\tau_j^{-1} \delta_j (u, 0)) \big| = \big| \det S_j^{-1} \big( DE_j(u,0) \big) \big| = \tau_j^{\sum_{i=1}^{\bar n} \theta_i} \delta_j^{-\sum_{i=1}^{\bar n} d_i} \big| \det  DE_j(u,0) \big| \\
&\approx \tau_j^{\sum_{i=1}^{\bar n} \theta_i} \delta_j^{-\sum_{i=1}^{\bar n} d_i} \big| \det  DE_j(0,0) \big| = \big| \det S_j^{-1} \big( DE_j(0,0) \big) \big| = \big| \det DH_j( 0, 0) \big|>0,
\end{align*}
and
\begin{align*}
&\quad \sup_{(i,j)\in I_0\times J^C} \big| \det \big( DH_j(\tau_j^{-1} \delta_j (u, 0)) \big)_{i^C\times j^C} \big|  = \tau_j^{\sum_{i=1}^{\bar n} \theta_i} \delta_j^{-\sum_{i=1}^{\bar n} d_i} \sup_{(i,j)\in I_0\times J^C}  \big| \det \big(DE_j(u,0) \big)_{i^C\times j^C} \big|\\
&\lesssim \tau_j^{\sum_{i=1}^{\bar n} \theta_i} \delta_j^{-\sum_{i=1}^{\bar n} d_i} \big| \det \big(D\Phi_j(u)\big)_{I} \big| \approx \tau_j^{\sum_{i=1}^{\bar n} \theta_i} \delta_j^{-\sum_{i=1}^{\bar n} d_i} \big| \det \big(D\Phi_j(0)\big)_{I} \big|\\
&= \tau_j^{\sum_{i=1}^{\bar n} \theta_i} \delta_j^{-\sum_{i=1}^{\bar n} d_i} \sup_{(i,j)\in I_0\times J^C}  \big| \det \big(DE_j(0,0) \big)_{i^C\times j^C} \big|\\
&= \sup_{(i,j)\in I_0\times J^C} \big| \det \big( DH_j(0, 0) \big)_{i^C\times j^C} \big|.
\end{align*}
Therefore for $j \geq N_2$ and $u\in B^{\bar n}(\eta_1)$,
\begin{align*}
&\quad \big\| \big( \big(DH_j(\tau_j^{-1} \delta_j (u, 0))\big)^{-1}\big)_{I^C\times I_0} \big\|_* = \left\| \frac{\big( \text{adj}\, DH_j(\tau_j^{-1} \delta_j (u,0)) \big)_{I^C\times I_0}}{\det DH_j(\tau_j^{-1} \delta_j (u, 0))} \right\|_* \\
& \lesssim \sup_{(i,j)\in I_0\times J^C} \frac{\big| \det \big(DH_j (\tau_j^{-1} \delta_j (u, 0))\big)_{i^C\times j^C} \big|}{\big|\det DH_j(\tau_j^{-1} \delta_j (u, 0))\big|} \lesssim \sup_{(i,j)\in I_0\times J^C} \frac{\big| \det \big(DH_j (0, 0)\big)_{i^C\times j^C} \big|}{\big|\det DH_j(0, 0)\big|}\\
&\leq \left\| \frac{\big( \text{adj}\, DH_j(0,0) \big)_{I^C\times I_0}}{\det DH_j(0, 0)} \right\|_*  =\big\| \big( DH_j(0, 0)^{-1}\big)_{I^C\times I_0} \big\|_* \lesssim 1,
\end{align*}
where $\text{adj}\,(\cdot)$ means the adjoint of a matrix, i.e., the transpose of the cofactor matrix. Then by shrinking $\bar \epsilon_j$, for $(u,v) \in B^{\bar n}(\frac{3}{4}\eta_1) \times B^{n-\bar n}(\bar \epsilon_j)$,
$$
\big\| \big( DH_j(\tau_j^{-1} \delta_j (u,v))^{-1}\big)_{I^C\times I_0} \big\|_* \lesssim 1.
$$

By (\ref{sally}), for $j \geq N_2$, there exists $0<\epsilon_j \leq \bar \epsilon_j$ such that for $t\in B^N(\rho')$ and $(u,v)\in B^{\bar n}(\eta_0)\times B^{n-\bar n}(\epsilon_j)$, $ E_j^{-1} \circ \gamma_{\tilde \tau_j^{-1} \delta_j t} \circ E_j(u,v)\in B^{\bar n}(\frac{3}{4}\eta_1) \times B^{n-\bar n}(\bar \epsilon_j)$, and thus
$$
\big\| \big(DH_j\big(\tau_j^{-1} \delta_j (E_j^{-1}\circ \gamma_{\tilde \tau_j^{-1} \delta_j t} \circ E_j(u,v)) \big)^{-1}\big)_{I^C\times I_0} \big\|_* \lesssim 1,
$$
which implies
\begin{align*}
&\quad \big\| \big( D(E_j^{-1} \circ \gamma_{\tilde \tau_j^{-1} \delta_j t} \circ E_j)\big)_{I^C \times J^C}(u,v) \big\|_* \\
&= \Big\| \Big( DE_j\big( E_j^{-1}\circ \gamma_{\tilde \tau_j^{-1} \delta_j t} \circ E_j(u,v) \big)\Big)^{-1}_{I^C \times I_0} \cdot D\gamma_{\tilde \tau_j^{-1} \delta_j t} \big( E_j(u,v) \big) \cdot \big(DE_j(u,v) \big)_{I_0 \times J^C}\Big\|_* \\
&= \Big\| S_j^*\Big(DH_j\big(\tau_j^{-1} \delta_j (E_j^{-1}\circ \gamma_{\tilde \tau_j^{-1} \delta_j t} \circ E_j(u,v)) \big)^{-1}\Big)_{I^C\times I_0} \cdot D\gamma_{\tilde \tau_j^{-1} \delta_j t} \big( E_j(u,v) \big) \cdot  S_j \big( DH_j(\tau_j^{-1} \delta_j(u,v)) \big)_{I_0\times J^C}\Big\|_* \\
&\leq  \big\| \Big(DH_j\big(\tau_j^{-1} \delta_j (E_j^{-1}\circ \gamma_{\tilde \tau_j^{-1} \delta_j t} \circ E_j(u,v)) \big)^{-1}\Big)_{I^C\times I_0} \big\|_* \cdot \\
& \qquad \qquad \qquad \qquad \qquad \qquad \big\| D\gamma_{\tilde \tau_j^{-1} \delta_j t} \big( E_j(u,v) \big) \big\|_* \cdot \big\| \big( DH_j(\tau_j^{-1} \delta_j(u,v))\big)_{I_0\times J^C}\big\|_* \\
&\lesssim C' \| \gamma_t(x) \|_{C^1(B^N(\rho)\times \Omega)}.
\end{align*}
Hence there exists $C_0\geq 1$ such that for $j \geq N_2$, $t\in B^N(\rho')$,
$$
\sup_{B^{\bar n}(\eta_0)\times B^{n-\bar n}(\epsilon_j)} \big\|  \big( D(E_j^{-1}\circ \gamma_{\tilde \tau_j^{-1} \delta_j t} \circ E_j) \big)_{I^C\times J^C}(u,v) \big\|_* \leq C_0.
$$
\end{proof}

\bibliographystyle{amsalpha}
\bibliography{ref}

\newcommand{\noop}[1]{}
\providecommand{\bysame}{\leavevmode\hbox to3em{\hrulefill}\thinspace}
\providecommand{\MR}{\relax\ifhmode\unskip\space\fi MR }
\providecommand{\MRhref}[2]{%
  \href{http://www.ams.org/mathscinet-getitem?mr=#1}{#2}
}
\providecommand{\href}[2]{#2}
\begin{thebibliography}{CNSW99}

\bibitem[CDSS20]{CHRISTSTREET}
Michael Christ, Spyridon Dendrinos, Betsy Stovall, and Brian Street,
  \emph{Endpoint {L}ebesgue estimates for weighted averages on polynomial
  curves}, Amer. J. Math. \textbf{142} (2020), no.~6, 1661--1731. \MR{4176542}

\bibitem[CHKY08]{CHKY08}
Yong-Kum Cho, Sunggeum Hong, Joonil Kim, and Chan~Woo Yang,
  \emph{Multiparameter singular integrals and maximal operators along flat
  surfaces}, Rev. Mat. Iberoam. \textbf{24} (2008), no.~3, 1047--1073.
  \MR{2490209}

\bibitem[CHKY09]{TRIPLE}
\bysame, \emph{Triple {H}ilbert transforms along polynomial surfaces}, Integral
  Equations Operator Theory \textbf{65} (2009), no.~4, 485--528. \MR{2576306}

\bibitem[CNSW99]{CNSW}
Michael Christ, Alexander Nagel, Elias~M. Stein, and Stephen Wainger,
  \emph{Singular and maximal {R}adon transforms: analysis and geometry}, Ann.
  of Math. (2) \textbf{150} (1999), no.~2, 489--577. \MR{1726701}

\bibitem[CWW00]{CWW00}
Anthony Carbery, Stephen Wainger, and James Wright, \emph{Double {H}ilbert
  transforms along polynomial surfaces in {$\mathbb{R}^3$}}, Duke Math. J.
  \textbf{101} (2000), no.~3, 499--513. \MR{1740686}

\bibitem[CWW06]{CWW06}
\bysame, \emph{Singular integrals and the {N}ewton diagram}, Collect. Math.
  (2006), no.~Vol. Extra, 171--194. \MR{2264209}

\bibitem[CWW09]{9}
\bysame, \emph{Triple {H}ilbert transforms along polynomial surfaces in
  {$\mathbb{R}^4$}}, Rev. Mat. Iberoam. \textbf{25} (2009), no.~2, 471--519.
  \MR{2554163}

\bibitem[Die60]{18}
J.~Dieudonn\'{e}, \emph{Foundations of modern analysis}, Pure and Applied
  Mathematics, Vol. X, Academic Press, New York-London, 1960. \MR{0120319}

\bibitem[Fab67]{FABES}
Eugene~B. Fabes, \emph{Singular integrals and partial differential equations of
  parabolic type}, Studia Math. \textbf{28} (1966/67), 81--131. \MR{213744}

\bibitem[Gal79]{GAL79}
Andr\'{e} Galligo, \emph{Th\'{e}or\`eme de division et stabilit\'{e} en
  g\'{e}om\'{e}trie analytique locale}, Ann. Inst. Fourier (Grenoble)
  \textbf{29} (1979), no.~2, vii, 107--184. \MR{539695}

\bibitem[Gre01a]{GREENBLATT2}
Michael Greenblatt, \emph{Boundedness of singular {R}adon transforms on {$L^p$}
  spaces under a finite-type condition}, Amer. J. Math. \textbf{123} (2001),
  no.~6, 1009--1053. \MR{1867310}

\bibitem[Gre01b]{GREENBLATT1}
\bysame, \emph{A method for proving {$L^p$}-boundedness of singular {R}adon
  transforms in codimension 1}, Duke Math. J. \textbf{108} (2001), no.~2,
  363--393. \MR{1833395}

\bibitem[Gre07a]{GREENBLATT}
\bysame, \emph{An analogue to a theorem of {F}efferman and {P}hong for
  averaging operators along curves with singular fractional integral kernel},
  Geom. Funct. Anal. \textbf{17} (2007), no.~4, 1106--1138. \MR{2373012}

\bibitem[Gre07b]{GREENBLATT6}
\bysame, \emph{A {$T(1)$} theorem for singular {R}adon transforms}, Math. Ann.
  \textbf{339} (2007), no.~3, 599--626. \MR{2336061}

\bibitem[Gre19a]{GREENBLATT5}
\bysame, \emph{{$L^p$} {S}obolev regularity for a class of {R}adon and
  {R}adon-like transforms of various codimension}, J. Fourier Anal. Appl.
  \textbf{25} (2019), no.~4, 1987--2003. \MR{3977144}

\bibitem[Gre19b]{GREENBLATT3}
\bysame, \emph{{$L^p$} {S}obolev regularity of averaging operators over
  hypersurfaces and the {N}ewton polyhedron}, J. Funct. Anal. \textbf{276}
  (2019), no.~5, 1510--1527. \MR{3912783}

\bibitem[Gre20]{GREENBLATT4}
\bysame, \emph{Smoothing theorems for {R}adon transforms over hypersurfaces and
  related operators}, Forum Math. \textbf{32} (2020), no.~6, 1637--1647.
  \MR{4168708}

\bibitem[Hu20]{SPARSE}
Bingyang Hu, \emph{Sparse domination of singular {R}adon transform}, J. Math.
  Pures Appl. (9) \textbf{139} (2020), 235--316. \MR{4108413}

\bibitem[Kim15]{KIM}
Joonil Kim, \emph{Multiple {H}ilbert transforms associated with polynomials},
  Mem. Amer. Math. Soc. \textbf{237} (2015), no.~1119, v+120. \MR{3400384}

\bibitem[Lob70]{LOB70}
Claude Lobry, \emph{Contr\^{o}labilit\'{e} des syst\`emes non lin\'{e}aires},
  SIAM J. Control \textbf{8} (1970), 573--605. \MR{0271979}

\bibitem[NRS01]{NRS01}
Alexander Nagel, Fulvio Ricci, and Elias~M. Stein, \emph{Singular integrals
  with flag kernels and analysis on quadratic {CR} manifolds}, J. Funct. Anal.
  \textbf{181} (2001), no.~1, 29--118. \MR{1818111}

\bibitem[NSW85]{NSW85}
Alexander Nagel, Elias~M. Stein, and Stephen Wainger, \emph{Balls and metrics
  defined by vector fields. {I}. {B}asic properties}, Acta Math. \textbf{155}
  (1985), no.~1-2, 103--147. \MR{793239}

\bibitem[NW77]{17}
Alexander Nagel and Stephen Wainger, \emph{{$L^{2}$} boundedness of {H}ilbert
  transforms along surfaces and convolution operators homogeneous with respect
  to a multiple parameter group}, Amer. J. Math. \textbf{99} (1977), no.~4,
  761--785. \MR{450901}

\bibitem[Pat08]{d=3}
Sanjay Patel, \emph{Double {H}ilbert transforms along polynomial surfaces in
  {$\mathbb{R}^3$}}, Glasg. Math. J. \textbf{50} (2008), no.~3, 395--428.
  \MR{2451738}

\bibitem[Pat09]{d=1}
\bysame, \emph{{$L^p$} estimates for a double {H}ilbert transform}, Acta Sci.
  Math. (Szeged) \textbf{75} (2009), no.~1-2, 241--264. \MR{2533413}

\bibitem[PY08]{8}
Malabika Pramanik and Chan~Woo Yang, \emph{Double {H}ilbert transform along
  real-analytic surfaces in {$\mathbb{R}^{d+2}$}}, J. Lond. Math. Soc. (2)
  \textbf{77} (2008), no.~2, 363--386. \MR{2400397}

\bibitem[She18]{ENDPOINT}
Jiawei Shen, \emph{Hardy space theory and endpoint estimates for
  multi-parameter singular {R}adon transforms}, Ph.D. thesis, Wayne State
  University, 2018.

\bibitem[SS12]{ANALYTIC}
Elias~M. Stein and Brian Street, \emph{Multi-parameter singular {R}adon
  transforms {III}: {R}eal analytic surfaces}, Adv. Math. \textbf{229} (2012),
  no.~4, 2210--2238. \MR{2880220}

\bibitem[SS13]{LP}
\bysame, \emph{Multi-parameter singular {R}adon transforms {II}: {T}he {$L^p$}
  theory}, Adv. Math. \textbf{248} (2013), 736--783. \MR{3107526}

\bibitem[SS18]{BETSY}
Betsy Stovall and Brian Street, \emph{Coordinates adapted to vector fields:
  canonical coordinates}, Geom. Funct. Anal. \textbf{28} (2018), no.~6,
  1780--1862. \MR{3881835}

\bibitem[Sto11]{BETSY11}
Betsy Stovall, \emph{{$L^p$} improving multilinear {R}adon-like transforms},
  Rev. Mat. Iberoam. \textbf{27} (2011), no.~3, 1059--1085. \MR{2895344}

\bibitem[Sto14]{STOVALL}
\bysame, \emph{Uniform {$L^p$}-improving for weighted averages on curves},
  Anal. PDE \textbf{7} (2014), no.~5, 1109--1136. \MR{3265961}

\bibitem[Str11]{CARNOT}
Brian Street, \emph{Multi-parameter {C}arnot-{C}arath\'{e}odory balls and the
  theorem of {F}robenius}, Rev. Mat. Iberoam. \textbf{27} (2011), no.~2,
  645--732. \MR{2848534}

\bibitem[Str12]{L2}
\bysame, \emph{Multi-parameter singular {R}adon transforms {I}: {T}he {$L^2$}
  theory}, J. Anal. Math. \textbf{116} (2012), 83--162. \MR{2892618}

\bibitem[Str17]{BRIAN17}
\bysame, \emph{Sobolev spaces associated to singular and fractional {R}adon
  transforms}, Rev. Mat. Iberoam. \textbf{33} (2017), no.~2, 633--748.
  \MR{3651019}

\bibitem[Str21a]{STREET20}
Brian Street, \emph{Coordinates adapted to vector fields {II}: Sharp results},
  Amer. Jour. Math., in press, arXiv:1808.04159, 2021.

\bibitem[Str21b]{COORDINATES}
\bysame, \emph{Coordinates adapted to vector fields {III}: Real analyticity},
  Asian J. Math., in press, arXiv:1808.04635, 2021.

\bibitem[SW03]{SW03}
Andreas Seeger and Stephen Wainger, \emph{Singular {R}adon transforms and
  maximal functions under convexity assumptions}, Rev. Mat. Iberoamericana
  \textbf{19} (2003), no.~3, 1019--1044. \MR{2053571}

\bibitem[TW03]{TW03}
Terence Tao and James Wright, \emph{{$L^p$} improving bounds for averages along
  curves}, J. Amer. Math. Soc. \textbf{16} (2003), no.~3, 605--638.
  \MR{1969206}

\bibitem[Vit17]{DOUBLE}
Marco Vitturi, \emph{Double {H}ilbert transforms along surfaces in the
  {H}eisenberg group}, Ph.D. thesis, The University of Edinburgh, 2017.

\bibitem[VW20]{12}
Marco Vitturi and James Wright, \emph{Multiparameter singular integrals on the
  {H}eisenberg group: uniform estimates}, Trans. Amer. Math. Soc. \textbf{373}
  (2020), no.~8, 5439--5465. \MR{4127882}

\bibitem[ZS75]{ZS75}
Oscar Zariski and Pierre Samuel, \emph{Commutative algebra. {V}ol. {II}},
  Springer-Verlag, New York-Heidelberg, 1975, Reprint of the 1960 edition,
  Graduate Texts in Mathematics, Vol. 29. \MR{0389876}

\end{thebibliography}

\vspace{2em}
\noindent
Lingxiao Zhang\\
Department of Mathematics\\
University of Wisconsin-Madison\\
480 Lincoln Dr., Madison, WI, 53706, USA\\
email: \href{mailto:lzhang395@wisc.edu}{lzhang395@wisc.edu}\\
MSC: primary 42B20; secondary 44A12, 32B99, 53C17

\end{document}